\documentclass[12pt,reqno]{amsart}

\usepackage[cp1251]{inputenc}
\usepackage{amscd,amsmath,amsxtra,amssymb,amstext,latexsym,amsthm}
\usepackage[all]{xy}
\usepackage{color}
\usepackage{array,graphicx}
\usepackage[all]{xy}

\def\leq{\leqslant}
\def\geq{\geqslant}

\newtheorem{thm}{Theorem}
\newtheorem{lem}
{Lemma}
\newtheorem{prop}
{Proposition}
\newtheorem{claim}
{Statement}
\newtheorem{df}
{Definition}
{Corollary}
\newtheorem{rem}
{Remark}
{Question}
\newtheorem{ex-thm}{Theorem-Example}

{\catcode`\@=11
\gdef\n@te#1#2{\leavevmode\vadjust{%
 {\setbox\z@\hbox to\z@{\strut#1}%
  \setbox\z@\hbox{\raise\dp\strutbox\box\z@}\ht\z@=\z@\dp\z@=\z@%
  #2\box\z@}}}
\gdef\leftnote#1{\n@te{\hss#1\quad}{}}
\gdef\rightnote#1{\n@te{\quad\kern-\leftskip#1\hss}{\moveright\hsize}}
\gdef\?{\FN@\qumark}
\gdef\qumark{\ifx\next"\DN@"##1"{\leftnote{\rm##1}}\else
 \DN@{\leftnote{\rm??}}\fi{\rm??}\next@}}

\begin{document}
\baselineskip=13.7pt plus 2pt

\title[On irreducible germs of generic morphisms]{On irreducible germs of generic morphisms}

\author[Vik.S. Kulikov]{Vik.S. Kulikov}

\address{Steklov Mathematical Institute of Russian Academy of Sciences, Moscow, Russia}
 \email{kulikov@mi-ras.ru}

\dedicatory{} \subjclass{}

\keywords{}

\maketitle

\def\st{{\sf st}}

\quad \qquad \qquad

\begin{abstract} 
The article examines a  set of irreducible germs $f_P:U_P\to V_p$ of 
generic morphisms $f:S\to\mathbb P^2$ to the projective plane whose branch curve germs $B_P\subset V_p$ have singularities 
equisingular deformation equivalent to singularities given by equations $x^{k_1}-y^{k_2}=0$ with coprime $k_1,k_2\in\mathbb N$.

\end{abstract}

\def\st{{\sf st}}

\setcounter{section}{-1}

\section{Introduction} 
\subsection{} Let $S$ be a non-singular irreducible projective surface defined over the field of complex numbers $\mathbb C$ and $f:S\to\mathbb P^2$ be a finite morphism to the projective plane $\mathbb P^2$ branched in a curve $B_f\subset\mathbb P^2$. 
Two finite morphisms $f_i:S_i\to \mathbb P^2$, $i=1,2$, are {\it equivalent} if there is an isomorphism $\varphi: S_1\to S_2$ such that $f_1=f_2\circ\varphi$.

Denote by $\mathcal F_{(2)}$ a set of the finite morphisms $f:S\to\mathbb P^2$  satisfying the following  
{\it generality conditions}:
\begin{itemize}
\item[(i)] {\it the branch curve $B_f\subset\mathbb P^2$ of $f$ is irreducible,}
\item[({ii})] {\it $f$ is ramified with multiplicity $2$ at a generic point of its ramification locus $R_f\subset S$,}
\item[({iii})] $f_{\mid R_f}:R_f\to B_f$ {\it is a bi-rational morphism,}
\end{itemize}
and call them the {\it generic morphisms to the projective plane}. Historically, a  generic morphism $f:S\to\mathbb P^2$ has been referred to as a  {\it generic cover of the projective plane} if in addition it satisfies the following condition:
\begin{itemize}
\item[({ iv})] {\it the singular points of $B_f$ are the ordinary nodes and ordinary cusps only.}
\end{itemize}

A morphism $f\in\mathcal F_{(2)}$ induces a {\it monodromy homomorphism}
$f_*:\pi_1(\mathbb P^2\setminus B_f,q)\to \mathbb S_{\deg f}$  to the symmetric group $\mathbb S_{\deg f}$ acting  on the fibre $f^{-1}(q)=\{ q_1,\dots,q_{\deg f}\}$. Its image $G_f:=f_*(\pi_1(\mathbb P^2\setminus B_f,q))\subset \mathbb S_{\deg f}$  is called the {\it monodromy group} of  $f$.

Let $p$ be a point of the curve $B_f\subset \mathbb P^2$. It is well known that the {\it local fundamental group} $\pi_1^{loc}(B_f,p):=\pi_1(V_p\setminus B_f)$ does not depend on $V_p$, where $V_p\subset \mathbb P^2$ is  a sufficiently small complex analytical neighborhood of the point $p$ bi-holomorphic to a ball of radius  $r\ll 1$ centered at $p$. The group  $\pi_1^{loc}(B_f,p)$ is generated by so called {\it geometric generators}, i.e., elements represented by positively oriented little loops in $V_p\setminus B_f$ around the irreducible components of the  curve germ $B_f$.

The image $G_{f,p}:=\text{im}\, f_{p*}:=\text{im}\, f_*\circ \iota_*$ is called the {\it local monodromy group} of $f$ at the point $p$, where $$\iota_*:\pi_1^{loc}(B_f,p)=\pi_1(V_p\setminus B_f,\widetilde q)\to \pi_1(\mathbb P^2\setminus B_f,q)$$ is a {\it local monodromy homomorphism} at $p$ defined uniquely up to conjugation in $\pi_1(\mathbb P^2\setminus B_f,q)$ by  imbedding $\iota: V_p\hookrightarrow \mathbb P^2$. 

Let the inverse image $f^{-1}(p)$ of a point $p\in \text{Sing}\, B_f$ consist of $k$ points,
$$f^{-1}(p)=\{ P_1,\dots,P_k\},$$  Then $f^{-1}(V_p)=\bigsqcup_{j=1}^kU_{P_j}\subset S$ is a disjoint union of $k$  connected neighborhoods $U_{P_j}$ of the points $P_j\in f^{-1}(p)$. We call $f_{P_j}:=f_{\mid U_{P_j}}:U_{P_j}\to V_p$ a {\it germ} of $f$ at the point $P_j$ and a germ $f_{{P_j}}$ is called {\it non-trivial} if $\deg f_{{P_j}}\geq 2$. Denote by $B_{P_j}\subset B_f\cap V_p$ the branch curve germ of a non-trivial germ $f_{{P_j}}:U_{P_j}\to V_p$ and denote by $R_{P_j}:=R_f\cap U_{P_j}$ the ramification curve germ of the germ $f_{P_j}:U_{P_j}\to V_p$. The inverse image $f_{{P_j}}^{-1}(B_{P_j})$ is the union of two curve germs, $f_{{P_j}}^{-1}(B_{P_j})=R_{P_j}\cup C_{P_j}\subset U_{P_j},$
where $C_{P_j}\subset U_{P_j}$ is the additional curve germ to the germ $R_{P_j}$ in $f_{P_j}^{-1}(B_{P_j})$. 

The imbeddings $\iota_j: V_p\setminus B_f\hookrightarrow V_p\setminus B_{P_j}$ define a homomorphism 
$$\prod f_{P_j*}\circ\iota_{j*}: \pi_1(V_p\setminus B_f)\to\prod G_{f_{P_j},p}=\prod\mathbb S_{\deg f_{P_j}}$$
into the direct product $\prod G_{f_{P_j},p}$ of the monodromy groups $G_{f_{P_j},p}$ of germs $f_{P_j}:U_{P_j}\to V_p$, where the direct product is taken over all $j$ for which $\deg f_{P_j}\geq 2$, and it is easy to see that the local monodromy homomorphism $f_{p*}:\pi_1^{loc}(B_f,p)\to \mathbb S_{\deg f}$ factorises into the composition of the homomorphism $\prod f_{P_j*}\circ\iota_{j*}$ and an imbedding $\prod\mathbb S_{\deg f_{P_j}}\hookrightarrow \mathbb S_{\deg f_p}$.

By Riemann - Stein Theorem (\cite{St}, \cite{G-R}), a non-trivial germ $f_{P_j}: U_{P_j}\to V_p$ of a generic morphism is defined uniquely by its monodromy homomorphism 
$f_{P_j*}:\pi_1^{loc}(B_{P_j},p)\to \mathbb S_{\deg f_{P_j}}$ that satisfies the following {\it generality conditions} (see Proposition \ref{prop10}):
\begin{itemize}
\item[$(\text{\bf i})$] $f_{P_j*}(\gamma)\in \mathbb S_{\deg f_{P_j}}$ are transpositions for  geometric generators $\gamma$ of  $\pi_1^{loc}(B_{P_j},p)$,
\item[$(\text{\bf ii})$] $f_{P_j*}$ is an epimorphism.
\end{itemize}

Denote by $\mathcal G_{(2)}$ the set of irreducible non-trivial germs $f_P:U_P\to V_p$ of generic morphisms to the projective plane.
A germ $\{ f_P:U_P\to V_p\}\in \mathcal G_{(2)}$ is said to be {\it quasi-generic} if it satisfies the {\it maximum degree condition}
\begin{itemize}
\item[$(\text{\bf iii})$] {\it   the degree $\deg f_{{P}}= \mu_p(B_{P})+1$, where $\mu_p(B_{P})$ is the multiplicity  at $p$ of the branch   curve germ  $B_{P}\subset V_p$ of the germ $f_P$.}
\end{itemize}

\begin{df} \label{def} A finite morphism $f\in \mathcal F_{(2)}$ is called a {\it quasi-generic cover} if  all its non-trivial germs  $f_{{P}}: U_{P}\to V_p$ 
are quasi-generic.
\end{df}

Denote by $\mathcal G_{(2),Q}$ the subset of $\mathcal G_{(2)}$ consisting  of quasi-generic germs and $\mathcal F_{(2),Q}$ the subset of $\mathcal F_{(2)}$ consisting of the quasi-generic covers.

Let $\{ f_P:U_P\to V_p\}\in \mathcal G_{2}$ and  $W$ be the normalization of the fibre product 
$$U_1\times _{V_p}U_2=\{ \, (x,y)\in U_1\times U_2\, \,  \mid \, \, f_1(x)=f_2(y)\, \, \} $$  
over $V_p$ of two copies $f_i:=f_P:U_i:=U_P\to V_p$, $i=1,2$, of the germ $f_P$ and let $R_i,C_i\subset U_i$ be curve germs corresponding to $R_P,C_P\subset U_P$. Denote by 
$$g_{1}:W\to U_1,\quad g_{2}:W\to U_2,\quad g_{1,2}:W\to V_p$$ the corresponding natural mappings.
In Subsection 2.3.3 of \cite{K-Q}, it was proved that $W=W'\bigsqcup W''$ is a disjoint union of two connected neighborhoods such that $g_{i\mid W'}: W'\to U_{i}$ are biholomorphic mappings and $g_{1,2\mid W'}: W'\to V_p$ coincides with $f_{1}$ and $f_{2}$.  The point $W''\cap g_{1,2}^{-1}(p)$ can be a singular point of $W''$. Let $\nu: \widetilde W\to W''$ be a resolution of singularity. Denote by $h_i=g_i\circ \nu$. We have
$$h_1^*(R_1)=\overline R+\overline C+E,$$
where $\overline R=h_1^{-1}(R_1)\cap h_2^{-1}(R_2)$, $\overline C=h_1^{-1}(R_1)\cap h_2^{-1}(C_2)$, and $E$ is a divisor supported 
in $\nu^{-1}(W''\cap g_{1,2}^{-1}(p))$. 

A germ $f_{P}: U_{P}\to V_p$ is said to have {\it extra property} if either $\deg f_{P}=2$, or the divisor $E$ can be represented as a sum $E=  E_{\overline R}+E_{\overline C}$ of two effective divisors  $E_{\overline R}$ and $E_{\overline C}$ such that
$$(\overline R+E_{\overline R},\overline C+E_{\overline C})_{\widetilde W}\leq  2 \delta_{P}(R_P)+\mu_{p}(B_P)-1$$ 
if $\deg f_{P}>2$, where $(Div_1,Div_2)_{\widetilde W}$ is the intersection number of divisors $Div_1$ and $Div_2$ in $\widetilde W$ and $\delta_P(R_P)$ is the $\delta$-invariant of the singularity of curve germ $R_P\subset U_P$ at $P$.

A quasi-generic cover $f:S\to\mathbb P^2$ is said to be {\it extra-quasi-generic} if over each point $p\in \text{Sing}\, B_f$, all its non-trivial germs have extra-property. Denote by $\mathcal F_{(2),E}\subset\mathcal F_{(2),Q}$ the set of extra-quasi-generic covers. 

\subsection{} The collection
$$\Gamma_f=\{ f_{p*}: \pi_1^{loc}(B_f,p)\to G_{f,p} \mid p\in Sing\, B_f\}$$
of local monodromy homomorphisms $f_{p*}: \pi_1^{loc}(B_f,p)\to G_{f,p}$  considered up to internal automorphisms of the groups $G_{f,p}$, is called the {\it local monodromy data} of a generic morphism $f:S\to\mathbb P^2$ and the pair $pas(f)=(B_f,\Gamma_f)$ is called the {\it passport} of $f$.

Chisini's Conjecture (\cite{Ch}) claims that two generic covers $f:S_i\to \mathbb P^2$, $i=1,2$, branched in a curve $B=B_{f_1}=B_{f_2}$ are equivalent if $\max(\deg f_1,\deg f_2)\geq 5$. In \cite{K1} -- \cite{N}, for morphisms belonging to some subsets of $\mathcal F_{(2)}$, it was proved  
theorems similar to Chisini's Conjecture (in \cite{K7} and \cite{K-Q}, called {\it Chisini Theorems}). 
Namely, a statement  is called a {\it Chisini Theorem  for morphisms from a subset $\mathcal M$ of $\mathcal F_{(2)}$} if it states that there is a constant $\frak{d}=\frak{d}(\mathcal M)\in\mathbb N$ such that if $f_1$ and $f_2\in \mathcal M$ satisfy conditions: $pas(f_1)=pas(f_2)$ and $\max(\deg f_1,\deg f_2)\geq \frak{d}$, then the morphisms $f_1$ and $f_2$ are equivalent.
In particular, in \cite{K-Q}, it was proved the following 
\begin{thm} \label{main} Chisini Theorem  with constant $\frak d=12$ is true for the extra-quasi-generic covers of the projective plane. \end{thm}

\subsection{} Let $S\subset \mathbb F_1*\mathcal C_2$ be a semigroup generated by elements $x_1$ and $x_2$ in the free product 
$\mathbb F_1*\mathcal C_2=\langle x_1,x_2\mid x_2^2=1\rangle$ of the free group $\mathbb F_1$ and the cyclic group $\mathcal C_2$ of second order.
Denote 
$$\begin{array}{rcl} \mathcal P=  \{ (k_1,k_2)\in \mathbb N^2\mid  
k_1,k_2\,\, \text{are coprime}\} & \simeq & \mathbb Q_+=\{ r=\frac{k_1}{k_2}\in \mathbb Q\mid r>0\}, \\
\mathcal O=  \{ (k_1,k_2)\in \mathcal P\mid  k_1\geq k_2\} & \simeq & \mathbb Q_{\geq 1}=\{ r=\frac{k_1}{k_2}\in \mathbb Q\mid r\geq 1\}. \end{array}$$

The semigroup $S$ acts on $\mathcal P$  (and, accordingly, on the set of positive rational numbers $\mathbb Q_+$) as follows:
$$ x_1((k_1,k_2))=(k_1+k_2,k_2),\qquad x_2((k_1,k_2))=(k_2,k_1). $$

The set of orbits under the action of $x_2$ on $\mathcal P$ is the set of pairs $\{ (k_1,k_2),(k_2,k_1)\}$. For shortness, denote an orbit 
$\{ (k_1,k_2),(k_2,k_1)\}$, where $k_1\geq k_2$, by $\{k_1,k_2\}$ and identify the set of orbits with the set $\mathcal O$. 

Define two actions $\alpha$ and $\beta$ on the set $\mathcal O$, 
$$\alpha: \{ k_{1},k_{2}\}\longmapsto \{ k_{1}+k_{2},k_{2}\}, \qquad \beta: \{ k_{1},k_{2}\}\longmapsto \{ k_{1}+k_{2},k_{1}\}$$ 

In Subsection \ref{sect1.3}, it will be shown that $\mathcal P$ is the orbit of $(1,1)\in\mathcal P$ under the action of $S$ on $\mathcal P$ and the set $\mathcal O$ is the orbit of $\{ 1,1\}\in \mathcal O$ under the action of the free semigroup $\mathcal S$ generated by $\alpha$ and $\beta$. 

\subsection{} 
Denote by $B_{k_1,k_2}\subset \mathbb C^2$ a curve given by equation $x^{k_1}-y^{k_2}=0$. Note that the linear transformation of the coordinates $x\leftrightarrow y$ determines the isomorphism of the pairs $(\mathbb C^2,B_{k_1,k_2})$ and $(\mathbb C^2,B_{k_2,k_1})$. We say that the singularity at $p\in V_p$ of a curve germ $(B,p)\subset (V,p)$ is $\mathcal P_{\{ k_1,k_2\}}$-{\it simplest} (where $\{ k_1,k_2\} \in \mathcal O$) if $(B,p)$ is equivisingular deformation equivalent to the curve germ whose Puiseux  series is either $y=x^{\frac{k_1}{k_2}}$ or  $y=x^{\frac{k_2}{k_1}}$  
and we say that the singularity of a curve germ $B$ is  {\it Puiseux simplest} (or simply $\mathcal P$-simplest) 
if it is one of $\mathcal P_{\{ k_1,k_2\}}$-simplest singularities. (Of course, the germs of curves given by $y=x^k$, $k\in \mathbb N$, are nonsingular and consequently they are isomorphic to each other. But, to simplify the formulation of the results (see also Remark \ref{rem-s} in Subsection \ref{sect-4.2}), we will assume that these curves germs are "singular"\, and will distinguish these types of singularities.)

\subsection{} The aim of the article is to describe a set $\mathcal G_{(2),\mathcal P}$ of the germs $\{ f_P:U_P\to V_p\}\in\mathcal G_{(2)}$ whose branch curve germs $B_P\subset V_p$ are irreducuble and have singularities of the $\mathcal P$-simplest types. 

In the first four sections of the article, some well-known facts are recalled, as well as the statements that are used in Section \ref{sect5} in proving the main results of the article are proved.

In Section \ref{sect1}, in order to describe the resolutions of singularities of $\mathcal P$-simplest types, we investigate a set $\mathcal D=\{ (k_1,k_2,q_1, q_2)\in \mathbb Z^4 \} $ of the solutions of equation 
\begin{equation}\label{sys1} x_1x_2-x_1x_4-x_2x_3=1,\end{equation}
and its subset $$\mathcal D_{\mathcal P}=\{ (k_1,k_2,q_1, q_2)\in \mathcal D\mid k_1>0,\,\, k_2>0,\,\,  0\leq q_1< k_1,\,\, 0\leq q_2< k_2\}.$$ 
Note that $(1,1,0,0)\in \mathcal D_{\mathcal P}$ 
and it is easy to check that if $(k_1,k_2,q_1,0)\in\mathcal D_{\mathcal P}$ in which $q_1\neq 0$, then $k_2=1$ and $q_1=k_1-1>0$. Also it is easy to see that 
the solutions $(k_1,k_2,q_1,q_2)\in\mathcal D_{\mathcal P}$ of equation (\ref{sys1}) have the following properties
\begin{itemize} 
\item[$1)$]  $k_1$ and $k_2$ are coprime, 
\item[$2)$] $k_1$ and $q_1$ (resp., $k_2$ and $q_2$) are coprime if $q_1\neq 0$ (resp., $q_2\neq 0$).
\end{itemize} 

Denote by $H \subset \text{GL}(4,\mathbb Z)$ the image of the group $\mathcal F_1*\mathcal C_2$ under a homomorphism $\varphi: \mathcal F_1*\mathcal C_2\to  \text{GL}(4,\mathbb Z)$ sending $x_1$ to $h_1$ and $x_2$ to $h_2$, where
$$ h_1= \left( \begin{array}{cccc} 
1 & 1 & 0 & 0 \\
0 & 1 & 0 & 0 \\ 
0 & 1 & 1 & -1 \\ 
0 & 0 & 0 & 1 
\end{array}\right), \qquad 
 h_2= \left( \begin{array}{cccc} 
0 & 1 & 0 & 0 \\
1 & 0 & 0 & 0 \\ 
0 & 0 & 0 & 1 \\ 
0 & 0 & 1 & 0 
\end{array}\right),   
$$
and denote by $S_{\mathcal D}$  the image of the semigroup $S$. Denote also by $\mathcal O_{\mathcal D}$ the set of orbits of elements of $\mathcal D_{\mathcal P}$ under the action of $h_2$ on $\mathcal D_{\mathcal P}$. For shortness, we denote the  orbit of $s=(k_1,k_2,q_1,q_2)\in\mathcal D_{\mathcal P}$, where $k_1\geq k_2$, by $\{ \frac{k_1}{q_1},\frac{k_2}{q_2}\}\in \mathcal O_{\mathcal D}$.

Let $\text{\rm pr}: \mathcal D_{\mathcal P}\to\mathcal P$, $\text{pr}: \mathcal O_{\mathcal D}\to\mathcal O$, and 
$\text{pr}_1: \mathcal O_{\mathcal D}\to\mathcal O$ be the projections given as follows: 
$$\begin{array}{rcl} \text{\rm pr}: (k_1,k_2,q_1,q_2)\in\mathcal D_{\mathcal P} & \mapsto &  (k_1,k_2)\in \mathcal P, \\ 
\text{\rm pr}: \{ \frac{k_1}{q_1},\frac{k_2}{q_2}\}\in\mathcal O_{\mathcal D} & \mapsto & \{k_1,k_2\}\in \mathcal O, \\ 
\text{\rm pr}_1: \{ \frac{k_1}{q_1},\frac{k_2}{q_2}\}\in\mathcal O_{\mathcal D} &\mapsto & \{ k_1,q_1\}\in \mathcal O.\end{array}$$

Define two actions on the set $\mathcal O_{\mathcal D}$ (denoted again by $\alpha$ and $\beta$) as follows:
$$\alpha: \{ \frac{k_{1}}{q_1},\frac{k_{2}}{q_2}\}\mapsto \{ \frac{k_{1}+k_{2}}{k_2+q_1-q_2},\frac{k_{2}}{q_2}\}, \quad   
\beta: \{ \frac{k_{1}}{q_1},\frac{k_{2}}{q_2}\}\mapsto \{ \frac{k_{1}+k_{2}}{k_1+q_2-q_1},\frac{k_{1}}{q_1}\}$$ 
and denote by $\mathcal S_{\mathcal D}$ a semigroup freely generated by the actions $\alpha$ and $\beta$.

The main result of Section \ref{sect1} is the following 
\begin{thm} \label{main1} $(1)$ The group $H$ leaves invariant the set $\mathcal D$.
\newline
$(2)$ The set $\mathcal D_{\mathcal P}$ is the orbit of $(1,1,0,0)$ under the action of $S_{\mathcal D}$  on $\mathcal D$. 
\newline
$(3)$ The set $\mathcal O_{\mathcal D}$ is the orbit of $\{ \frac{1}{0},\frac{1}{0}\}$ under the action of $\mathcal S_{\mathcal D}$  on 
$\mathcal O_{\mathcal D}$. 
\newline
$(4)$ The projections 
$\text{\rm pr}: \mathcal D_{\mathcal P}\to \mathcal P$, $\text{\rm pr}: \mathcal O_{\mathcal D}\to \mathcal O$, and 
$\text{\rm pr}_1: \mathcal O_{\mathcal D}\to \mathcal O$ are bijections.
\end{thm}

The set of branch curve germs $B_P\subset V_p$ of germs $\{ f_P:U_P\to V_p\}\in\mathcal G_{(2),\mathcal P}$ is closely connected with  the subset
$$ \begin{array}{rl}
D_{\mathcal P}= \{ & \widetilde s=(k_1,k_2,q_1,q_2,q_3,q_4,m_1,m_2) \in\mathbb Z^8\mid \\ 
& \min(k_1,k_2)\geq 2,\,\,   
m_1\geq 0,\, m_2\geq 0, \,\,   
0< q_1<k_1, \\
&  0< q_2<k_2, \,\,\, 0< q_3<k_1+k_2,\,\,\,  0< q_4< k_1k_2\,\,\,  \} 
\end{array} $$
of the set of solutions of the system of equations  (\ref{sys1}), (\ref{syst2}), (\ref{syst3}): 
$$ \phantom{aaaaaaaaaaaaaaaaaaaaaaaa} x_1x_2-x_1x_4-x_2x_3=1 \phantom{aaaaaaaaaaaaaaaaaaaaaaaa} (1) $$ 
\begin{equation} \label{syst2} x_1x_2(x_1+x_2)-x_1x_2x_5- x_6(x_1+x_2)=1 \end{equation}
\begin{equation}\label{syst3} x_6=x_1x_7+x_3=x_2x_8+x_4  \end{equation}

In Section \ref{sect2}, we prove the following
\begin{thm} \label{mainD} The projection 
$$\text{\rm pr}_{D_{\mathcal P},\mathcal D_{\mathcal P}} :\widetilde s=(k_1,k_2,q_1,q_2,q_3,q_4,m_1,m_2)\in D_{\mathcal P}\longmapsto
\overline s=(k_1,k_2,q_1,q_2)\in\mathcal D_{\mathcal P, 0}$$
is a bijection, where $\mathcal D_{\mathcal P, 0}=\{ \overline s=(k_1,k_2,q_1,q_2)\in\mathcal D\mid \min(k_1,k_2)\geq 2\}$.
\end{thm}

In Section \ref{sect3}, we remind facts connected with expansions of rational numbers $\frac{k}{q}>1$ as continued fractions, and in Section \ref{sect4}, we investigate properties of resolution of singularities of curve germs having $\mathcal P$-simplest singularity types.

Denote by $\mathcal G_{(2),\mathcal O}$ a set of germs $\{ f_{P}:U_{P}\to V_p\}\in  \mathcal G_{(2),\mathcal P}$ branched in curve 
germs $B_{P}$ which are equisingular deformation equivalent to curve germs given by equations $x^{k_1k_2}-y^{k_1+k_2}=0$, 
$\{ k_1,k_2\}\in \mathcal O$, $k_2>1$, and such that $\deg f_{P}=\mu_p(B_{P})=k_1+k_2$, 
denote  $\mathcal G_{(2),\overline{\mathcal D}}= \mathcal G_{(2),\mathcal D}\cup\mathcal G_{(2),\mathbb N}$, where $\mathcal G_{(2),\mathcal D}$ is a set of germs $\{ f_{P}:U_{P}\to V_p\}\in  \mathcal G_{(2),\mathcal P}$ branched in curve germs $B_{P}$ which are equisingular deformation equivalent to curve germs given by equations $x^{k_1(k_2+1)}-y^{k_2}=0$, where $\{ k_1,k_2\}\in \mathcal D_{\mathcal P,0}$, and such that $\deg f_{P}=\mu_p(B_{P})=k_2+1$, and 
$\mathcal G_{(2),\mathbb N}$ is a set of germs $\{ f_{P}:U_{P}\to V_p\}\in  \mathcal G_{(2),\mathcal P}$ branched in curve germs $B_{P}$ which are equisingular deformation equivalent to  germs given by equations $x^{k+1}-y^{k}=0$,  $k\in \mathbb N$, and such that $\deg f_{P}=\mu_p(B_{P})=k+1$.

The main result of the article, proved in Section \ref{sect5}, is the following two theorems.

\begin{thm} \label{main0} $(\bf 1)$ For each $\{ k_1,k_2\}\in \mathcal O$, $k_2>1$, there is a single $($up to equivalence$)$ germ 
$\{ f_{P}:U_{P}\to V_p\}\in \mathcal G_{(2),\mathcal O}$  branched in  a curve germ $B_{P}$ which is equisingular deformation equivalent to a germ given by equation $x^{k_1k_2}-y^{k_1+k_2}=0$. \newline 
$(\bf 2)$ For each $\overline s=(k_1,k_2,q_1,q_2)\in\mathcal D_{\mathcal P,0}$ and for each $(k_1,k_2)$ where $k_1=1$ and $k_2\in \mathbb N$, there is a
single $($up to equivalence$)$ germ 
$\{ f_{P}:U_{P}\to V_p\}\in \mathcal G_{(2),\overline{\mathcal D}}$  branched in  a curve germ $B_{P}$ which is equisingular deformation equivalent to a germ given by equation $x^{k_1(k_2+1)}-y^{k_2}=0$. \newline    
$(\bf 3)$ The set $\mathcal G_{(2),\mathcal P}$ is the disjoint  union of two subsets, 
$\mathcal G_{(2),\mathcal P}=\mathcal G_{(2),\mathcal O}\cup\mathcal G_{(2),\overline{\mathcal D}}$.
\end{thm}

\begin{thm} \label{main00} 
Chisini Theorem  with constant $\frak d=12$ is true for the generic morphisms to the projective plane whose all irreducible germs belong to 
$\mathcal G_{(2),\mathbb N}$.
\end{thm}

In view of Theorems \ref{main0} and \ref{main00}, the following question is of interest: {\it Are there any Chisini Theorems for any subsets of the set of generic morphisms to the projective plane whose irreducible germs belong to the set  $\mathcal G_{(2),\mathcal O}\cup\mathcal G_{(2),\mathcal D}$}? \\

Further  we use freely the notation introduced above.

\section{Proof of Theorem \ref{main1}}\label {sect1}

\subsection{Euclid's additive algorithm}\label{sect1.3} Consider the set $\mathcal O$ of orbits   under the action of $x_2$ on $\mathcal P$. 

In this subsection, we recall Euclid's additive algorithm for finding the greatest  common divisor of natural numbers $k_{0,1}=k_1$ and $k_{0,2}=k_2$. It acts on the set $\mathcal O$ and at  $i$th step of the algorithm, we make either a replacement $\varepsilon_{1,i}$ or a replacement $\varepsilon_{2,i}$, where
\begin{equation} \label{epsilon} \begin{array}{lr} \varepsilon_{1,i} :\{ k_{i-1,1},k_{i-1,2}\} \mapsto \{ k_{i,1},k_{i,2}\}= \{ k_{i-1,1}-k_{i-1,2},k_{i-1,2}\} \quad  &
\text{\rm if}\,\, k_{i-1,1}\geq 2k_{i-1,2}, \\ 
\varepsilon_{2,i} :\{ k_{i-1,1},k_{i-1,2}\} \mapsto \{ k_{i,1},k_{i,2}\}=\{ k_{i-1,2},k_{i-1,1}-k_{i-1,2}\} \quad  & \text{\rm if}\,\, k_{i-1,1}\leq 2k_{i-1,2}.\end{array} \end{equation} 
These steps are done until we get a pair $(k_{n,1},k_{n,2})$ in which $k_{n,1}=k_{n,2}$.  

Denote by $\text{\rm n}_{\mathcal E}(k_1,k_2):=n$ the number of steps of Euclid's additive algorithm applied to an orbit $\{ k_1,k_2\}$ to obtain the orbit $\{ \text{g.c.d.}(k_1,k_2),\text{g.c.d.}(k_1,k_2)\}$. 

In the case when $(k_1,k_2)\in\mathcal P$, we have $(k_{n,1},k_{n,2})=(1,1)$. 

The inverse replacements $\varepsilon_1^{-1}:=\varepsilon_{1,i}^{-1}$ and $\varepsilon_2^{-1}:=\varepsilon_{2,i}^{-1}$ are 
\begin{equation}\label{alpha-beta} \alpha:=\varepsilon_{1}^{-1}: \{ k_{1},k_{2}\}\mapsto \{ k_{1}+k_{2},k_{2}\}, \quad   
\beta:=\varepsilon_{2}^{-1}: \{ k_{1},k_{2}\}\mapsto \{ k_{1}+k_{2},k_{1}\}. \end{equation} 

The replacements $\alpha$ and $\beta $ generate a semigroup $\mathcal S$ acting on the set $\mathcal O$.  
It is easy to see that $\mathcal O$ is the orbit of $\{ 1,1\}$ under the action of $\mathcal S$. 

If we do one more step of Euclid's additive algorithm, then we obtain the pair $\{ 1,0\}$ from $\{ 1,1\}$. We denote by $\overline{\mathcal O}=\mathcal O\cup\{\{ 1,0\}\}$ and associate an oriented graph $\Gamma_{\overline{\mathcal O}}$ with the set $\overline{\mathcal O}$  as follows. 
\begin{itemize} 
\item[$(\text{\rm g}_1)$] The vertices $v$ of $\Gamma_{\overline{\mathcal O}}$ correspond one-to-one to the orbits $\{ k_1,k_2\}\in \overline{\mathcal O}$, and if this does not lead to a misunderstanding, we denote the vertices in the same way as their corresponding orbits, i.e. by $\{ k_1,k_2\}$.
\item[$(\text{\rm g}_2)$] The vertex $\{ 1,0\}$ is of valence $1$, the vertex $\{ 1,1\}$ is of valence $2$, and all other vertices are of valence $3$. 
\item[$(\text{\rm g}_3)$] The vertex $\{ 1,0\}$ is the end of an edge labeled by $\varepsilon_2$ and whose start vertex is $\{ 1,1\}$, the vertex  $\{ 1,1\}$ is the end of an edge labeled by $\varepsilon_1$ and whose start vertex is $\{ 2,1\}$, and all over edges  
    are labeled  either by the replacement $\varepsilon_1$, or by the replacement $\varepsilon_2$. 
\item[$(\text{\rm g}_4)$] Each vertex $\{ k_1,k_2\}$ is the start vertex of the only one edge labeled  either by $\varepsilon_1$ if $k_1>2k_2$, 
or by $\varepsilon_2$ if $k_1<2k_2$.
\item[$(\text{\rm g}_5)$] Each vertex $\{ k_1,k_2\}$, $k_1>k_2\geq 1$, is the end of two edges having different labels. 
\end{itemize}
 
A vertex $\{ k_1,k_2\}$ is said to be of {\it level} $\text{\rm ln}(\{ k_1,k_2\})=n_{\mathcal E}(k_1,k_2)+1$. The vertex $\{ 1,0\}$ is called the root of $\Gamma_{\overline{\mathcal O}}$. Note that $\text{\rm ln}(\{ k_1,k_2\})$ is the length of the single path along the edges of $\Gamma_{\overline{\mathcal O}}$ connecting the vertex $\{ k_1,k_2\}$ with $\{ 1,0\}$.  

\subsection{Proof of statement $(1)$ of Theorem \ref{main1}} 
It is easy to see that $h_2(s)\in\mathcal D$ for $s=(k_1,k_2,q_1,q_2)\in\mathcal D$ and the linear transformation $h_1$ acts as follows: 
\begin{equation} \label{equx} h_1: s=(k_1,k_2,q_1,q_2) \longmapsto h_1(s)=(k_1+k_2,k_2, k_2+q_1-q_2,q_2).\end{equation}  
We have $ (k_1+k_2)k_2-(k_1+k_2)q_2-k_2(k_2+q_1-q_2)=k_1k_2-k_1q_2-k_2q_1=1$ and therefore $h_1(\mathcal D)\subset \mathcal D$.
Similarly, it is easy to see that 
\begin{equation}\label{equy} h_1^{-1}: s=(k_1,k_2,q_1,q_2) \longmapsto h_1^{-1}(s)=(k_1-k_2,k_2,-k_2+q_1+q_2,q_2)\end{equation}  
and 
$ (k_1-k_2)k_2-(k_1-k_2)q_2-k_2(-k_2+q_1+q_2)=k_1k_2-k_1q_2-k_2q_1=1.$
Therefore the group $H$ leaves invariant the set $\mathcal D$. \qed

\subsection{Proof of statement $(2)$ of Theorem \ref{main1}}
First of all, note that if $k_1=k_2=k$ in $ s=(k_1,k_2,q_1,q_2)\in\mathcal D_{\mathcal P}$, then 
$k(k-q_1-q_2)=1$, since $s=(k_1,k_2,q_1,q_2)$ is a solution of equation (\ref{sys1}). Consequently, 
$k=1$ and $q_1=q_2=0$, i.e., $s=(1,1,0,0)$, since $0\leq q_1,q_2<k$. 

Let us show that if $k_1\neq k_2$ in $s=(k_1,k_2,q_1,q_2)\in\mathcal D_{\mathcal P}$,  then 
\begin{equation} \label{ine} \min(k_1,k_2)\leq q_1+q_2. \end{equation} 
Indeed, if, for example, $k_1>k_2$ and we assume the opposite, then $k_1>1$ and $$1=k_1k_2-k_1q_2-k_2q_1>k_1k_2-k_1q_2-k_1q_1=k_1(k_2-q_1-q_2)>k_1> 1,$$
i.e., we obtain a contradiction. 

Now, to complete the proof of statement $(2)$, it suffices to note that if $k_1> k_2$ in $s=(k_1,k_2,q_1,q_2)\in\mathcal D_{\mathcal P}$, then it follows from (\ref{equy}) that applying $h_1^{-1}$ to $s$, we get  $h_1^{-1}(s)=(k_1',k_2', q_1',q_2')\in\mathcal D_{\mathcal P}$ for which $\max(k_1',k_2')<\max(k_1,k_2)$. Therefore, after several steps of applying transformations $h_2$ and $h_1^{-1}$ to an element $\overline s \in\mathcal D_{\mathcal P}$, we will get $s_0=(1,1,0,0)$. \qed

\subsection{Proof of statement $(3)$ of Theorem \ref{main1}} 
Consider the set $\mathcal O_{\mathcal D}$ of orbits of the action $h_2$ on the set $\mathcal D_{\mathcal P}$.
For shortness, we denote an element $s=(k_1,k_2,q_1,q_2)\in \mathcal D_{\mathcal P}$ by $(\frac{k_1}{q_1},\frac{k_2}{q_2})$ and an orbit 
$\{ (\frac{k_1}{q_1},\frac{k_2}{q_2}),(\frac{k_2}{q_2},\frac{k_1}{q_1})\}$   under the action of $h_2$ on $\mathcal D_{\mathcal P}$
by $\{ \frac{k_1}{q_1},\frac{k_2}{q_2}\}$, where  $k_1\geq k_2$.   

The linear transformations $\alpha=h_1$ and $\beta=h_1h_2$ define two actions  on the set $\mathcal O_{\mathcal D}$ of orbits  
$\{ \frac{k_1}{q_1},\frac{k_2}{q_2}\}=\{ (\frac{k_1}{q_1},\frac{k_2}{q_2}),(\frac{k_2}{q_2},\frac{k_1}{q_1})\}$ denoted again by $\alpha$ and $\beta$, 
$$ \alpha: \{ \frac{k_{1}}{q_1},\frac{k_{2}}{q_2}\}\mapsto \{ \frac{k_{1}+k_{2}}{k_2+q_1-q_2},\frac{k_{2}}{q_2}\}, \quad   
\beta: \{ \frac{k_{1}}{q_1},\frac{k_{2}}{q_2}\}\mapsto \{ \frac{k_{1}+k_{2}}{k_1+q_2-q_1},\frac{k_{1}}{q_1}\}$$
and define two inverse transformations:
$$\varepsilon_{1} :\{ \frac{k_1}{q_1},\frac{k_2}{q_2}\} \mapsto  \{ \frac{k_1-k_2}{q_1+q_2-k_2},\frac{k_2}{q_2}\} \quad  
\text{\rm if}\,\, k_{1}> 2k_{2},$$  
$$\varepsilon_{2} :\{ \frac{k_1}{q_1},\frac{k_2}{q_2}\} \mapsto \{ \frac{k_2}{q_2},\frac{k_1-k_2}{q_1+q_2-k_2}\} \quad 
\text{\rm if}\,\, k_{1}< 2k_{2}.$$ 
According to inequality (\ref{ine}), the transformations $\varepsilon_1$ and $\varepsilon_2$ are correctly defined.

It is easy to see that 
\begin{equation} \label{com} \alpha (\text{\rm pr}(\{ \frac{k_{1}}{q_1},\frac{k_{2}}{q_2}\}))=\text{\rm pr}(\alpha(\{ \frac{k_{1}}{q_1},\frac{k_{2}}{q_2}\})), \quad 
\beta (\text{\rm pr}(\{ \frac{k_{1}}{q_1},\frac{k_{2}}{q_2}\}))=\text{\rm pr}(\beta(\{ \frac{k_{1}}{q_1},\frac{k_{2}}{q_2}\})),\end{equation} 
where $\text{\rm pr}(\{ \frac{k_{1}}{q_1},\frac{k_{2}}{q_2}\})=\{ k_1,k_2\}$ is the projection defined in Introduction.

Let us add one more element signed by $\{ \frac{1}{0},\frac{0}{-1}\}$ to the set $\mathcal O_{\mathcal D}$ and denote the obtained set by 
$\overline{\mathcal O}_{\mathcal D}=\mathcal O_{\mathcal D}\cup\{ \{ \frac{1}{0},\frac{0}{-1}\}\}$.
 
We associate an oriented graph $\Gamma_{\overline{\mathcal O}_{\mathcal D}}$ with the set $\overline{\mathcal O}_{\mathcal D}$  as follows. 
\begin{itemize} 
\item[$(\text{\rm \bf g}_1)$] The vertices $v$ of $\Gamma_{\overline{\mathcal O}_{\mathcal D}}$ correspond one-to-one to the orbits 
$\{ \frac{k_1}{q_1},\frac{k_2}{q_2}\}\in \overline{\mathcal O}_{\mathcal D}$, and if this does not lead to a misunderstanding, we denote the vertices in the same way as their corresponding orbits, i.e. by $\{ \frac{k_1}{q_1},\frac{k_2}{q_2}\}$.
\item[$(\text{\rm\bf g}_2)$] The vertex $\{ \frac{1}{0},\frac{0}{-1}\}$ is of valence $1$, the vertex $\{ \frac{1}{0},\frac{1}{0}\}$ is of valence $2$, and all other vertices are of valence $3$. 
\item[$(\text{\rm \bf g}_3)$] The vertex $\{ \frac{1}{0},\frac{0}{-1}\}$ is the end of an edge labeled by $\varepsilon_2$ and whose start vertex is 
$\{ \frac{1}{0},\frac{1}{0}\}$, the vertex  $\{ \frac{1}{0},\frac{1}{0}\}$ is the end of an edge labeled by $\varepsilon_1$ and whose start vertex is 
$\{ \frac{2}{1},\frac{1}{0}\}$, and all over edges are labeled  either by the replacement $\varepsilon_1$, or by the replacement $\varepsilon_2$. 
\item[$(\text{\rm \bf g}_4)$] Each vertex $\{ \frac{k_1}{q_1},\frac{k_2}{q_2}\}$ is the start vertex of the only one edge labeled  either by $\varepsilon_1$ if $k_1>2k_2$, or by $\varepsilon_2$ if $k_1<2k_2$.
\item[$(\text{\rm g}_5)$] Each vertex $\{ \frac{k_1}{q_1},\frac{k_2}{q_2}\}$, $k_1>k_2\geq 1$, is the end of two vertex having different labels.    
\end{itemize}
 
A vertex $\{ \frac{k_1}{q_1},\frac{k_2}{q_2}\}$ is said to be of {\it level} $\text{\rm ln}(\{ \frac{k_1}{q_1},\frac{k_2}{q_2}\})=n_{\mathcal E}(k_1,k_2)+1$. The vertex $\{ \frac{1}{0},\frac{0}{-1}\}$ is called the root of $\Gamma_{\overline{\mathcal O}_{\mathcal D}}$. Note that 
$\text{\rm ln}(\{ \frac{k_1}{q_1},\frac{k_2}{q_2}\})$ is the length of the single path along the edges of $\Gamma_{\overline{\mathcal O}_{\mathcal D}}$ connecting the vertex $\{ \frac{k_1}{q_1},\frac{k_2}{q_2}\}$ with $\{ \frac{1}{0},\frac{0}{-1}\}$.  

A part of the graph $\Gamma_{\overline{\mathcal O}_{\mathcal D}}$ containing the vertices  whose level $\leq 5$ looks as follows.   

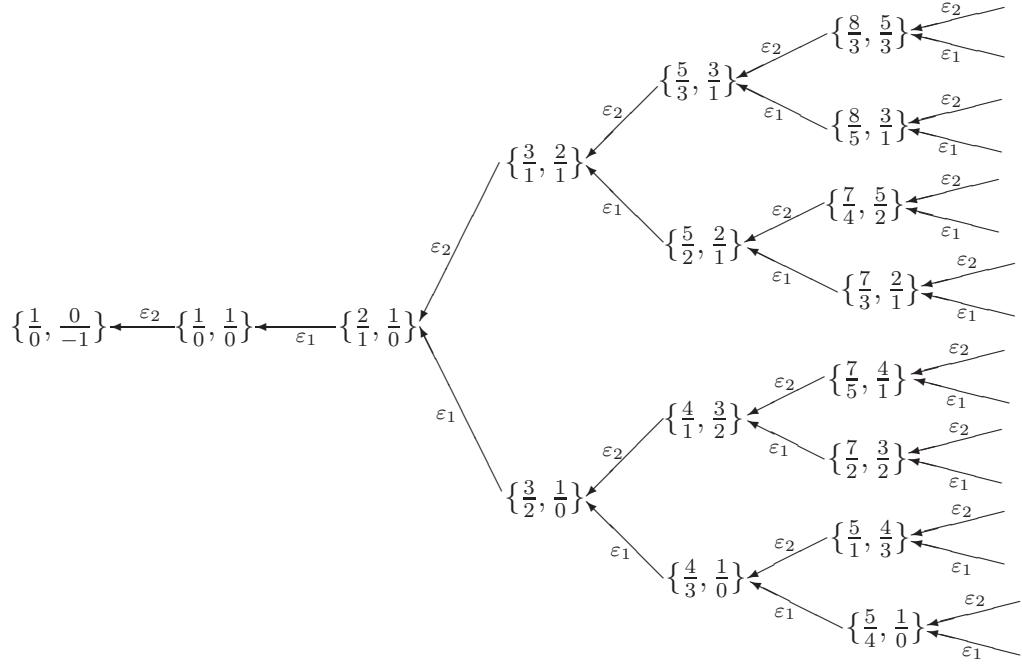
\begin{picture}(400,270)
\put(24,122){\mbox{$\{ \frac{1}{0},\frac{0}{-1}\}$}}
\put(86,125){\vector(-1,0){24}}
\put(67,130){\mbox{$\phantom{a}_{\varepsilon_2}$}}
\put(86,122){\mbox{$\{ \frac{1}{0},\frac{1}{0}\}$}}
\put(147,125){\vector(-1,0){30}}
\put(126,121){\mbox{$\phantom{a}_{\varepsilon_1}$}}
\put(148,122){\mbox{$\{ \frac{2}{1},\frac{1}{0}\}$}}
\put(209,187){\vector(-1,-2){30}}
\put(177,155){\mbox{$\phantom{a}_{\varepsilon_2}$}}
\put(210,63){\vector(-1,2){31}}
\put(179,92){\mbox{$\phantom{a}_{\varepsilon_1}$}}
\put(211,184){\mbox{$\{ \frac{3}{1},\frac{2}{1}\}$}}
\put(269,216){\vector(-1,-1){27}}
\put(242,207){\mbox{$\phantom{a}_{\varepsilon_2}$}}
\put(271,157){\vector(-1,1){29}}
\put(242,170){\mbox{$\phantom{a}_{\varepsilon_1}$}}
\put(211,57){\mbox{$\{ \frac{3}{2},\frac{1}{0}\}$}}
\put(271,90){\vector(-1,-1){29}}
\put(242,78){\mbox{$\phantom{a}_{\varepsilon_2}$}}
\put(271,30){\vector(-1,1){29}}
\put(245,40){\mbox{$\phantom{a}_{\varepsilon_1}$}}
\put(269,215){\mbox{$\{ \frac{5}{3},\frac{3}{1}\}$}}
\put(333,236){\vector(-2,-1){34}}
\put(302,231){\mbox{$\phantom{a}_{\varepsilon_2}$}}
\put(333,200){\vector(-2,1){34}}
\put(303,206){\mbox{$\phantom{a}_{\varepsilon_1}$}}
\put(334,233){\mbox{$\{ \frac{8}{3},\frac{5}{3}\}$}}
\put(400,246){\vector(-4,-1){35}}
\put(370,246){\mbox{$\phantom{a}_{\varepsilon_2}$}}
\put(400,227){\vector(-4,1){35}}
\put(370,228){\mbox{$\phantom{a}_{\varepsilon_1}$}}
\put(334,198){\mbox{$\{ \frac{8}{5},\frac{3}{1}\}$}}
\put(399,211){\vector(-4,-1){35}}
\put(370,211){\mbox{$\phantom{a}_{\varepsilon_2}$}}
\put(399,191){\vector(-4,1){35}}
\put(369,193){\mbox{$\phantom{a}_{\varepsilon_1}$}}
\put(271,153){\mbox{$\{ \frac{5}{2},\frac{2}{1}\}$}}
\put(332,172){\vector(-2,-1){30}}
\put(306,169){\mbox{$\phantom{a}_{\varepsilon_2}$}}
\put(337,138){\vector(-2,1){34}}
\put(306,144){\mbox{$\phantom{a}_{\varepsilon_1}$}}
\put(332,168){\mbox{$\{ \frac{7}{4},\frac{5}{2}\}$}}
\put(398,181){\vector(-4,-1){35}}
\put(370,181){\mbox{$\phantom{a}_{\varepsilon_2}$}}
\put(398,161){\vector(-4,1){35}}
\put(370,161){\mbox{$\phantom{a}_{\varepsilon_1}$}}
\put(338,136){\mbox{$\{ \frac{7}{3},\frac{2}{1}\}$}}
\put(404,149){\vector(-4,-1){35}}
\put(376,150){\mbox{$\phantom{a}_{\varepsilon_2}$}}
\put(404,129){\vector(-4,1){35}}
\put(376,130){\mbox{$\phantom{a}_{\varepsilon_1}$}}
\put(333,102){\mbox{$\{ \frac{7}{5},\frac{4}{1}\}$}}
\put(400,116){\vector(-4,-1){35}}
\put(373,117){\mbox{${\phantom{a}_{\varepsilon_2}}$}}
\put(402,96){\vector(-4,1){36}}
\put(372,97){\mbox{$\phantom{a}_{\varepsilon_1}$}}
\put(333,72){\mbox{$\{ \frac{7}{2},\frac{3}{2}\}$}}
\put(399,86){\vector(-4,-1){35}}
\put(373,87){\mbox{$\phantom{a}_{\varepsilon_2}$}}
\put(399,66){\vector(-4,1){35}}
\put(373,66){\mbox{$\phantom{a}_{\varepsilon_1}$}}
\put(334,42){\mbox{$\{ \frac{5}{1},\frac{4}{3}\}$}}
\put(400,55){\vector(-4,-1){35}}
\put(374,56){\mbox{$\phantom{a}_{\varepsilon_2}$}}
\put(400,35){\vector(-4,1){35}}
\put(374,34){\mbox{$\phantom{a}_{\varepsilon_1}$}}
\put(340,8){\mbox{$\{ \frac{5}{4},\frac{1}{0}\}$}}
\put(406,21){\vector(-4,-1){35}}
\put(379,22){\mbox{$\phantom{a}_{\varepsilon_2}$}}
\put(406,1){\vector(-4,1){35}}
\put(378,2){\mbox{$\phantom{a}_{\varepsilon_1}$}}
\put(271,87){\mbox{$\{ \frac{4}{1},\frac{3}{2}\}$}}
\put(332,106){\vector(-2,-1){29}}
\put(307,104){\mbox{$\phantom{a}_{\varepsilon_2}$}}
\put(332,75){\vector(-2,1){29}}
\put(305,79){\mbox{$\phantom{a}_{\varepsilon_1}$}}
\put(272,27){\mbox{$\{ \frac{4}{3},\frac{1}{0}\}$}}
\put(333,45){\vector(-2,-1){30}}
\put(307,43){\mbox{$\phantom{a}_{\varepsilon_2}$}}
\put(339,11){\vector(-2,1){35}}
\put(307,17){\mbox{$\phantom{a}_{\varepsilon_1}$}}
\end{picture}
  
Denote by 
$\text{\rm pr}_{\Gamma_{\overline{\mathcal O}_{\mathcal D}},\Gamma_{\overline{\mathcal O}}}:  \Gamma_{\overline{\mathcal O}_{\mathcal D}}\to 
\Gamma_{\overline{\mathcal O}}$ a mapping sending the vertices $\{ \frac{k_1}{q_1},\frac{k_2}{q_2}\}$  
of $\Gamma_{\overline{\mathcal O}_{\mathcal D}}$ to the vertices $\{ k_1,k_2\}$ of $\Gamma_{\overline{\mathcal O}}$  and sending the edges,  connecting vertices of $\Gamma_{\overline{\mathcal O}_{\mathcal D}}$ and labeled by  $\varepsilon_1$ and $\varepsilon_2$,  respectively, to the edges of $\Gamma_{{\overline{\mathcal O}}}$ having the same labels  $\varepsilon_1$ and $\varepsilon_2$ and connecting the images of the vertices.  
 
\begin{claim} \label{pro1} 
The mapping $\text{\rm pr}_{\Gamma_{\overline{\mathcal O}_{\mathcal D}},\Gamma_{\overline{\mathcal O}}}:  \Gamma_{\overline{\mathcal O}_{\mathcal D}}\to \Gamma_{\overline{\mathcal O}}$ is an isomorphism of oriented graphs.
\end{claim}
\proof  It is easy to see that for each $\{ k_1,k_2\}\in\overline{\mathcal O}$, starting from the roots of graphs $\Gamma_{\overline{\mathcal O}_{\mathcal D}}$ and $\Gamma_{\overline{\mathcal O}}$ and moving along the edges  having the same labels, we can reach simultaneously the vertex 
$\{ k_1,k_2\}\in \Gamma_{\overline{\mathcal O}}$ and a vertex $\{ \frac{k_1}{q_1},\frac{k_2}{q_2}\}\in \Gamma_{\overline{\mathcal O}_{\mathcal D}}$. 
Therefore the mapping $\text{\rm pr}_{\Gamma_{\overline{\mathcal O}_{\mathcal D}},\Gamma_{\overline{\mathcal O}}}$ is a cover. 

Assume that there are two different orbits $\{ \frac{k_1}{q_1},\frac{k_2}{q_2}\}\in\mathcal O_{\mathcal D}$ and $\{ \frac{k_1}{q'_1},\frac{k_2}{q'_2}\}\in\mathcal O_{\mathcal D}$ such that 
$\text{\rm\ pr}_{\Gamma_{\overline{\mathcal O}_{\mathcal D}},\Gamma_{\overline{\mathcal O}}}(\{ \frac{k_1}{q_1},\frac{k_2}{q_2}\})=
\text{\rm\ pr}_{\Gamma_{\overline{\mathcal O}_{\mathcal D}},\Gamma_{\overline{\mathcal O}}}(\{ \frac{k_1}{q'_1},\frac{k_2}{q'_2}\})$. 
It follows from (\ref{sys1}) that
$$ k_1(q_2-q_2')+k_2(q_1-q'_1)=0.$$
Therefore, there must be $n\in\mathbb Z$ such that $q_1-q'_1=nk_1$ and $q_2-q'_2=-nk_2$, since $k_1$ and $k_2$ are coprime. On the other hand, we have inequalities: $0\leq q_i,q'_i<k_i$ for $i=1,2$. Therefore $|q_1-q'_1|<k_1$ and $|q_2-q'_2|<k_2$, and hence $n=0$. Consequently, the mapping 
$\text{\rm pr}_{\Gamma_{\overline{\mathcal O}_{\mathcal D}},\Gamma_{\overline{\mathcal O}}}:  \Gamma_{\overline{\mathcal O}_{\mathcal D}}\to 
\Gamma_{\overline{\mathcal O}}$ is an injection. \qed \\

Now, Statement $(3)$ of Theorem \ref{main1} follows from Statement \ref{pro1} and definition of transformations $\varepsilon_1$ and $\varepsilon_2$, since we can reach any vertex $\{ \frac{k_1}{q_1},\frac{k_2}{q_2}\}\in \Gamma_{\overline{\mathcal O}_{\mathcal D}}$ moving along the path in 
$\Gamma_{\overline{\mathcal O}_{\mathcal P}}$ connecting the vertices $\{ \frac{1}{0},\frac{0}{-1}\}$ and $\{ \frac{k_1}{q_1},\frac{k_2}{q_2}\}$. \\

\subsection{Proof of statement $(4)$ of Theorem \ref{main1}} 
The statements that projections $\text{\rm pr}: \mathcal D_{\mathcal P}\to \mathcal P$ and $\text{\rm pr}: \mathcal O_{\mathcal D}\to \mathcal O$ are bijections directly follow from Statement \ref{pro1}. 

To prove that the projection $\text{\rm pr}_1: \mathcal O_{\mathcal D}\to \mathcal O$ is a surjection, note that 
since $k_1$ and $q_1$ are coprime positive integers, then there is a uniquely defined positive integer $k_2<k_1$ such that 
$q_1k_2+1=0 (\text{\rm mod}\, k_1)$. In addition, there is a unique natural number $n<k_2$ such that $q_1 k_2+1=k_1 n$. A direct check shows that
$\overline s=(k_1,k_2,q_1,k_2-n)\in \mathcal D_{\mathcal P}$ and hence, $\{ \frac{k_1}{q_1},\frac{k_2}{k_2-n}\}\in \text{\rm pr}_1^{-1}(\{ \frac{k_1}{q_1}\})$. 

To prove that the projection $\text{\rm pr}_1: \mathcal O_{\mathcal D}\to \mathcal O$ is an injection, assume that there are two different orbits 
$\{ \frac{k_1}{q_1},\frac{k_2}{q_2}\}, \{ \frac{k_1}{q_1}, \frac{k'_2}{q'_2}\}  \in\mathcal O_{\mathcal D}$. 
Then it follows from definition of the sets $\mathcal D_{\mathcal P}$ and $\mathcal O_{\mathcal D}$ that 
\begin{equation}\label{cont1} k_1(q_2-q'_2)-(k_1-q_1)(k_2-k'_2)=0\end{equation} 
and if $k_2=k'_2$ then $q_2=q'_2$. Therefore, without loss of generality, we can assume that $k_2>k'_2$ and hence, $0< k:=k_2-k'_2<k_1$, since $0<k_2,k'_2\leq k_1$. Let $d=\text{g.c.d.}(k,q)$, where $q=q_2-q'_2$, and let  $k:=da_1$ and $q=da_2$.

It follows from (\ref{cont1}) that $k_1a_2-(k_1-q_1)a_1=0$, where $0<a_1<k_1$ and $\text{g.c.d.}(a_1,a_2)=1$. But, in this case $k_1=a_1d_1$, where $d_1>1$, and, consequently, $\text{g.c.d.}(k_1,k_1-q_1)>1$ which contradicts the fact that $k_1$  and $q_1$ are coprime numbers.  Therefore the projection 
$\text{\rm pr}_1: \mathcal O_{\mathcal D}\to \mathcal O$ is an injection and Statement $(4)$ of Theorem \ref{main1} is proved. 

\begin{lem}\label{lem-ineq} The following inequalities  
$ k_2\leq q_1+q_2<k_1$
hold for all orbits $\{ \frac{k_1}{q_1},\frac{k_2}{q_2}\}  \in\mathcal O_{\mathcal D}\setminus \{\{ \frac{1}{0},\frac{1}{0}\}\}$.
\end{lem}
\proof For all $\{ \frac{k_1}{q_1},\frac{k_2}{q_2}\}  \in\mathcal O_{\mathcal D}\setminus \{\{ \frac{1}{0},\frac{1}{0}\}\}$, the inequalities $1\leq k_2<k_1$ hold. Therefore 
$$k_1(k_2-q_2-q_1)<k_1k_2-k_1q_2-k_2q_1=1<k_2(k_1-q_2-q_1)$$
and hence $k_2-q_2-q_1\leq 0$ and $k_1-q_2-q_1>0$. \qed 

\section{Proof of Theorem \ref{mainD} }\label{sect2}
\subsection{One more system of Diophantine equations}\label{sect2.1}  
Consider a subset
$$ \begin{array}{rrrl}
\widetilde D_{\mathcal P}= \{ & \widetilde s & = & (k_1,k_2,q_1,q_2)\times (a_1,a_2) \in\mathbb Z^4\times \mathbb Z^2\mid \min(k_1,k_2)\geq 2,  \\ 
&  & & 0< q_1<k_1,\,\,  0< q_2<k_2,    0< a_1\leq k_1,\,\,\,  0< a_2\leq k_2\,\,\,  \} 
\end{array} $$
of the set of solutions of the system of equations  (\ref{sys1}), (\ref{syst4}): 
$$ \phantom{aaaaaaaaaaaaaaaaaaaaaaaa} x_1x_2-x_1x_4-x_2x_3=1 \phantom{aaaaaaaaaaaaaaaaaaaaaaaa} (1) $$ 
\begin{equation}\label{syst4} x_1y_2-x_2y_1=x_3-x_4.  \end{equation}
Denote by $\text{\rm pr}_{4}: (k_1,k_2,q_1,q_2)\times(a_1,a_2)\mapsto (k_1,k_2,q_1,q_2)$ 
the restriction on $\widetilde D_{\mathcal P}$ of the projection $\text{pr}:\mathbb Z^4\times\mathbb Z^2\to \mathbb Z^4$ to the first factor. 
\begin{prop} \label{proj1} The projection $\text{pr}_{4}: \widetilde D_{\mathcal P}\to \mathcal D_{\mathcal P}$ is a bejection.
\end{prop} 
\proof Assume that two different solutions $\overline s=s\times (a_1,a_2)$ and $\overline s'=s\times (a'_1,a'_2)$ belong to 
 $\text{pr}_{4}^{-1}(s)$ for some $ s=(k_1,k_2,q_1,q_2)\in\mathcal D_{\mathcal P}$. Then 
$k_1a_2-k_2a_1=q_1-q_2$ and $k_1a'_2-k_2a'_1=q_1-q_2$. Therefore 
\begin{equation}\label{cont} k_1(a_2-a'_2)-k_2(a_1-a'_1)=0\end{equation} 
and, without loss of generality, we can assume that $0< b_2:=a_2-a'_2<k_2$, since $0<a_2,a'_2\leq k_2$. Denote $b_1=a_1-a'_1$ and let  
$d=\text{g.c.d.}(b_1,b_2)$, $b_1:=dc_1$ and $b_2=dc_2$.

It follows from (\ref{cont}) that $k_1c_2-k_2c_1=0$, where $0<c_2<k_2$. But, $k_1$ and $k_2$ are coprime integers. Therefore equality (\ref{cont}) is impossible if $a_2\neq a_2'$ and hence, $\text{pr}_{4}: \widetilde D_{\mathcal P}\to \mathcal D_{\mathcal P}$ is an injection. 

Without loss of generality, we can assume that $k_1>k_2$. Then $|q_1-q_2|<k_1$, since $0<q_2<k_2<k_1$ and $q_1<k_1$. 

There are three possibilities: either $q_1=q_2:=q$, or $0<q_1-q_2:=b_1<k_1$, or $0<q_2-q_1:=b_2<k_2$ (since $q_2<k_2$). 

In the first case, it is easy to see that $\widetilde s=(k_1,k_2,q,q)\times(k_1,k_2)\in \widetilde D_{\mathcal P}$.

In the second case, there is a natural number $a_1<k_1$ such that $k_2a_1 +b_1=0 (\text{\rm mod}\, k_1)$, i.e., there exists also a natural number 
$a_2$ such that $k_1a_2=k_2a_1+b_1$, and it is easy to see that $a_2$ must be less then $k_2$. Therefore 
$$\text{\rm pr}_4^{-1}((k_1,k_2,q_1,q_2))=\widetilde s=(k_1,k_2,q_1,q_2)\times(a_1,a_2)\in \widetilde D_{\mathcal P}.$$

In the third case, there is a natural number $a_2<k_2$ such that $k_1a_2 +b_2=0 (\text{\rm mod}\, k_2)$, i.e., there exists also a natural number 
$a_1$ such that $k_2a_1=k_1a_2+b_2$, and it is easy to see that $a_1$ must be less then $k_1$. Therefore 
$$\phantom{aaaaaaaaaaaaa} \text{\rm pr}_4^{-1}((k_1,k_2,q_1,q_2))=\widetilde s=(k_1,k_2,q_1,q_2)\times(a_1,a_2)\in \widetilde D_{\mathcal P}.\phantom{aaaaaaa} \qed$$

\subsection{Description of the set $D_{\mathcal P}$}\label{sect2.2} 
To find the last four coordinates $(q_3,q_4,m_1,m_2)$ of $\overline s=(k_1,k_2,q_1,q_2,q_3,q_4,m_1,m_2)\in \text{pr}_{\widetilde D_{\mathcal P},
\mathcal D_{\mathcal P}}^{-1}(s)$ whose  first four coordinates of $s=(k_1,k_2,q_1,q_2)\in \mathcal D_{\mathcal P,0}$ are known,      
we need to resolve a system of Diophantine equations obtained after substitution   $k_1,k_2,q_1,q_2$ in (\ref{sys1}), (\ref{syst2}), (\ref{syst3}) 
instead of variables $x_1,\dots,x_4$. After this substitution, we obtain  that $(q_3,q_4,m_1,m_2)$ is a solution of the following system of equalities and inequalities  
 
\begin{equation} \label{sys3} k_1k_2(k_1+k_2)-k_1k_2x_5- x_6(k_1+k_2)=1 \end{equation}
\begin{equation}\label{sys4} 0<x_6=k_1x_7+q_1=k_2x_8+q_2<k_1k_2,  \end{equation}
\begin{equation}\label{sys5}   
0< x_5<k_1+k_2,\quad   0\leq x_7,\quad 0\leq x_8\end{equation}
the coefficients of which  are related by the following equality  
\begin{equation} \label{sys2}
 k_1k_2-k_1q_2-k_2q_1=1 \end{equation}
 and satisfy the following inequalities
\begin{equation}\label{sys6} 
0< q_1<k_1,\quad 0< q_2<k_2, \quad \min(k_1,k_2)\geq 2. \end{equation} 

Theorem \ref{mainD} is a consequence of  Proposition \ref{proj1} and the following
\begin{prop}\label{prop-m1} The correspondence 
$$\widetilde s=(k_1,k_2,q_1,q_2)\times (a_1,a_2)\,\, \longmapsto \,\, \overline s=(k_1,k_2,q_1,q_2,q_3,q_4,m_1,m_2)$$
between the elements $\widetilde s\in {\widetilde D_{\mathcal P}}$ and the elements  $\overline s\in D_{\mathcal P}$, where 
\begin{equation}\label{cor} q_3=a_1+a_2-1,\,\,  q_4=k_2(k_1-a_1)+q_1,,\,\, m_1=k_1-a_1,\,\  m_2=k_2-a_2,\end{equation}
is one-to-one. 
\end{prop}
\proof 
It follows from inequalities (\ref{sys4}), (\ref{sys5}), and (\ref{sys6}) that 
\begin{equation}\label{x7,x8} 0\leq x_7< k_2,\qquad 0\leq x_8< k_1. \end{equation}

By (\ref{sys4}), we can express $q_1$ as $q_1=k_2x_8-k_1x_7+q_2$ and substitute this expression into equality (\ref{sys2}). We obtain an equation
$$\begin{array}{c} k_2(k_1-(k_2x_8-k_1x_7+q_2))-1-k_1q_2 \stackrel{(\ref{sys2})}{=} \\ 
k_2k_1-k_2(k_2x_8-k_1x_7+q_2)-k_1k_2+k_2q_1=0
\end{array}$$
that is equivalent to equation  
\begin{equation} \label{ms2} k_1x_7-k_2x_8=q_2-q_1, \end{equation}
since $k_2\neq 0$. 

Applying equalities (\ref{sys4}) (and after that equality (\ref{sys2})) to (\ref{sys3}) we obtain 
$$\begin{array}{c} k_1k_2(k_1+k_2)-x_5k_1k_2- x_6(k_1+k_2)= \\
k_1k_2(k_1+k_2)-x_5k_1k_2- (k_2x_8+q_2)k_1-(k_1x_7+q_1)k_2= \\
k_1k_2(k_1+k_2)-x_5k_1k_2- k_1k_2x_8-k_1k_2x_7 -k_1q_2- k_2q_1 \stackrel{(\ref{sys2})}{=} \\
k_1k_2(k_1+k_2)-x_5k_1k_2- k_1k_2x_8-k_1k_2x_7 -k_1k_2 +1= \\
k_1k_2[(k_1+k_2)-x_5- x_8-x_7 -1] +1=1.
\end{array} $$
Consequently, 
\begin{equation} \label{ms1} x_7+x_8=k_1+k_2-x_5-1, 
\end{equation} 
since $k_1k_2\neq 0$. 

Let's replace variables:
$$ x_7=k_2-y_2,\qquad x_8=k_1-y_1. $$
Applying (\ref{ms2}) and (\ref{x7,x8}), we obtain
\begin{equation}\label{ms0} k_1y_2-k_2y_1=q_1-q_2,  \quad 0<y_1\leq k_1,\,\, 0<y_2\leq k_2.\end{equation}
It follows from (\ref{ms1}) that 
\begin{equation} x_5=y_1+y_2-1.\end{equation}

Let $(a_1,a_2)$ be a solution of equation (\ref{ms0}). Then $s\times (a_1,a_2)$ is a solution of system (\ref{sys1}), (\ref{syst4}), where $s=(k_1,k_2,q_1,q_2)$ is a solution of equation (\ref{sys1}). Then if we put 
$$ m_1=k_1-a_1,\,\, m_2=k_2-a_2,\,\, q_3=a_1+a_2-1,\,\,q_4=k_1(k_1-a_1)+q_1$$
and if we perform the calculations made above in reverse order, we get that  
$$\widetilde s=(k_1,k_2,q_1,q_2,q_3,q_4,m_1,m_2)\in \text{pr}_{\widetilde D_{\mathcal P},\mathcal D_{\mathcal P}}^{-1}(s).$$ 
is a solution of system of equations $(\ref{sys1})$, (\ref{syst2}), $(\ref{syst3})$. \qed \\

Note that the following equality
\begin{equation} \label{Delta} m_1+m_2+q_3=k_1+k_2-1 \end{equation}
holds for $\widetilde s=(k_1,k_2,q_1,q_2,q_3,q_4,m_1,m_2)\in \widetilde D_{\mathcal P}$.

\section{Continued fractions and graphs of curves}\label{sect3}
\subsection{Expansions in continued fractions}
Let $\mathbb Z[u_1,\dots,u_n,\dots ]$ be the ring of polynomials in  variables $u_1,\dots,u_n,\dots$  with coefficients in $\mathbb Z$ and
  $\mathcal R=\mathbb Q(u_1,\dots,u_n,\dots)$ be its field of quotients. Consider  rational functions
\begin{equation}\label{H-J}
R_n(u_1,\dots, u_n)=u_n-\frac{1}{u_{n-1}-\frac{1}{u_{n-2}- \frac{1}{\dots -\frac{1}{u_1}}}}\in \mathcal R.
\end{equation}
By induction, it is easy to check that 
\begin{equation} \label{R} R_n(u_1,\dots,u_n)=\frac{D_n(u_1,\dots,u_n)}{D_{n-1}(u_1,\dots,u_{n-1})}, 
\end{equation} where
the polynomials $D_n(u_1,\dots,u_n)$  are given recurrently:
$$ D_{-1}:=0,\,\, D_0:=1,\,\, D_1:=D_1(u_1)=u_1,$$
$$D_2:=D_2(u_1,u_2)= u_2D_1-D_0=u_1u_2-1,$$
$$D_3:=D_3(u_1,u_2,u_3)=u_3D_2-D_1=u_1u_2u_3-u_1-u_3,$$
$$\dots\dots\dots$$
\begin{equation}\label{rec} D_n:=D_n(u_1,\dots,u_n)=u_nD_{n-1}-D_{n-2}.
\end{equation}

For $n\in \mathbb N$, denote by $M_{n,+}(u_1,\dots, u_n)$ and $M_{n,-}(u_1,\dots, u_n)$ the following matrices
$$ M_{n,\pm}:=M_{n,\pm}(u_1,\dots,u_n)= \left( \begin{array}{ccccccc}
u_1 & \pm 1   & 0 & \dots  & \dots     & \dots & 0 \\
\pm 1   & u_2 & \pm 1   & \ddots & \ddots & \ddots     &  \vdots  \\
0   & \pm 1   & u_3 & \ddots & \ddots & \ddots & \vdots  \\
\vdots & \ddots & \ddots & \ddots &  \ddots & \ddots & \vdots \\
\vdots & \ddots & \ddots & \ddots & u_{n-2} & \pm 1  & 0 \\
\vdots & \ddots & \ddots & \ddots & \pm 1  & u_{n-1} & \pm 1 \\ 
0 & \dots & \dots & \dots & 0 & \pm 1  & u_{n} 
\end{array} \right) $$ 
and $ M^-_{n,\pm}:=M_{n,\pm}(-u_1,\dots,-u_n)$.

\begin{lem}\label{prop-1} $D_n=\det M_{n,+}=\det M_{n,-}$ and $\det  M^-_{n,\pm}=(-1)^nD_n$.
\end{lem}
\noindent {\it Proof} by induction on $n$ applying equality (\ref{rec}). \qed 

Permuting the rows and columns of $M_{n,+}$ and applying Lemma \ref{prop-1}, we obtain that 
\begin{equation}\label{M} D_n(u_1,u_2,\dots,u_n)=D_n(u_n,u_{n-1},\dots, u_1).
\end{equation}

\begin{lem} \label{pos-def} For $w_1,\dots, w_n\in\mathbb Z$, the matrices $M_{n,\pm}(w_1,\dots,w_n)$ are positive definite if and only if 
the matrices $M_{n,\pm}(-w_1,\dots,-w_n)$ are negative definite. 

If $w_i\geq 2$ for $1\leq i\leq n$, then $M_{i,\pm}(w_1,\dots,w_i)$ are positive definite matrix and $0<d_i=\det M_{i,\pm}(w_1,\dots,w_i)<d_{i+1}=\det M_{i+1,\pm}(w_1,\dots,w_{i+1})$. 
\end{lem}
\noindent {\it Proof} directly follows from recurrence relation (\ref{rec}), Lemma \ref{prop-1} and Sylvester's criterion of positive definiteness of matrices. \qed \\

Denote  $D_{n_0,\text{\rm lt},i}:=D_{n_0-i-1}(u_1,\dots, u_{n_0-1-i})$, where $1\leq n_0\leq n$ and $0\leq i\leq n_0$ (recall that $D_0=1$ and $D_{-1}=0$), and  $D_{n_0,\text{rt},i}:=D_{n-n_0-i}(u_{n_0+1+i},\dots,u_n)$, where $i=0,\dots, n-n_0$.

Decomposing  the determinant of $M_{n,+}$ by its $n_0$th row, it follows from Lemma \ref{prop-1}, that 
\begin{equation}\label{eq0+}  D_n=  u_{n_0}D_{n_0,\text{\rm lt},0}D_{n_0,\text{\rm rt},0}-  
 D_{n_0,\text{\rm rt},0}D_{\text{\rm lt,1}}-    D_{n_0,\text{\rm lt},0}D_{\text{\rm rt},1}.\end{equation}

\subsection{Weighted chains}\label{sect3.2} A graph $\Gamma$ depicted in Figure 1, is called a {\it weighted chain of length 
$\text{\rm ln}(\Gamma)=\text{\rm n}(\Gamma)-1\geq 0$}, were $\text{\rm n}(\Gamma)$ is the number of its vertices  
and $w_i\in\mathbb Z$ are the weights of the vertices $v_i\in \Gamma$, $1\leq i\leq n(\Gamma)$. 

\begin{picture}(300,60)
\put(99,37){$w_1$}\put(105,30){\circle*{4}}\put(100,18){$v_{1}$}
\put(145,30){\circle*{4}}\put(141,37){$w_2$}
\put(140,19){$v_2$}\put(108,30){\line(1,0){40}}
\put(145,30){\line(1,0){40}}\put(215,30){\line(1,0){40}} \put(255,30){\circle*{4}}
\put(195,30){$.\, .\, .$}
\put(243,19){$v_{n(\Gamma)-1}$} \put(243,37){$w_{n(\Gamma)-1}$} \put(257,30){\line(1,0){40}} \put(299,30){\circle*{4}}
\put(294,37){$w_{n(\Gamma)}$} \put(295,18){$v_{n(\Gamma)}$}
\put(190,0){$\text{Fig.}\, 1$}
\end{picture} \vspace{0.5cm}

\noindent If $\text{\rm n}(\Gamma)>1$, then the chain $\Gamma$ has two vertices of valence $1$. We choose one of them, denote it by $v_{\text{or}(\Gamma)}$, call it the {\it origin}  of chain $\Gamma$, and number the vertices of $\Gamma$ consistently starting from the origin ($v_{\text{or}(\Gamma)}=v_1$). The vertex $v_{\text{\rm n}(\Gamma)}$ is called the {\it end} of $\Gamma$ (if $\text{\rm n}(\Gamma)=1$, then at the same time, the vertex $v_1$ is the origin and the end of $\Gamma$). If the origin is chosen then $\Gamma$ is called an {\it oriented chain}. If $\Gamma$ is an oriented chain, when by $\Gamma^{\text{\rm inv}}$ is denoted the oriented chain $\Gamma$ with orientation defined by the choice of $v_{\text{\rm n}}(\Gamma)$ as the origin.   

Denote $M(\Gamma):=M_{\text{\rm n}(\Gamma),+}(w_1,\dots,w_{\text{\rm n}(\Gamma)}),$
\begin{equation}\label{not} \begin{array}{c} D(\Gamma) :=\det M(\Gamma)=D_{\text{\rm n}(\Gamma)}(w_1,\dots,w_{\text{\rm n}(\Gamma)}), \\  
D_{\text{lt},i}(\Gamma):=D_{\text{\rm n}(\Gamma)-i}(w_{1},\dots, w_{\text{\rm n}(\Gamma)-i}), \\ D_{\text{rt},i}(\Gamma):=D_{\text{\rm n}(\Gamma)-i}(w_{i+1},\dots, w_{\text{\rm n}(\Gamma)}), \\ 
d(\Gamma) :=|D(\Gamma)|,\quad d_{\text{rt},i}(\Gamma):=|D_{\text{rt},i}(\Gamma)|, \quad d_{\text{lt},i}(\Gamma):=|D_{\text{lt},i}(\Gamma)|.
\end{array} \end{equation} 

It follows from (\ref{M}) and Lemma \ref{prop-1}, that
\begin{equation}\label{equ} D(\Gamma)=w_1D_{\text{rt},1}(\Gamma)-D_{\text{rt},2}(\Gamma)= w_{\text{\rm n}(\Gamma)}D_{\text{lt},1}(\Gamma)-D_{\text{lt},2}(\Gamma).
\end{equation} 

Let $\Gamma_1$ and $\Gamma_2$ be two oriented chains. Denote by $\Gamma_1\oplus\Gamma_2$ the oriented chain obtained as follows: we add the edge to the disjoint union of $\Gamma_1$ and $\Gamma_2$  connecting the end of $\Gamma_1$ with the origin of $\Gamma_2$ and by definition, the origin of $\Gamma_1\subset\Gamma_1\oplus\Gamma_2$ is the origin of $\Gamma_1\oplus\Gamma_2$. 

Let an oriented graph $\Gamma$ depicted in Fig. 2 be represented as a sum of three oriented graphs, $\Gamma=\Gamma_1\oplus\Gamma_2\oplus\Gamma_3$, where $\Gamma_2$ consists of a single vertex $v_{n_0}$ having weight $w_{n_0}$, $\text{n}(\Gamma_1)=n_0-1$ and $\text{n}(\Gamma_3)=n-n_0$. The subgraph $\Gamma_1\subset \Gamma$ is called the {\it left part} of the graph $\Gamma$, the subgraph $\Gamma_3\subset \Gamma$ is called its {\it right part}, and the vertex $v_{n_0}$ is called the {\it center} of $\Gamma$.

\begin{picture}(300,65)
\put(45,30){\circle*{3}}
\put(40,35){$\mbox{$w$}_{1}$}
\put(42,20){$\mbox{$v$}_{1}$}
\put(45,30){\line(1,0){40}}
\put(85,30){\circle*{3}}
\put(80,35){$\mbox{$w$}_{2}$}
\put(82,20){$\mbox{$v$}_{2}$}
\put(85,30){\line(1,0){40}}
\put(132,30){$\dots$}
\put(150,30){\line(1,0){40}}
\put(186,35){$\mbox{$w_{n_0}$}$}
\put(190,30){\circle*{3}}
\put(187,20){$\mbox{$v_{n_0}$}$}
\put(192,30){\line(1,0){40}}
\put(235,30){$\dots$}
\put(335,30){\circle*{2}}
\put(287,35){$\mbox{$w_{n-1}$}$}
\put(295,30){\circle*{3}}
\put(290,20){\mbox{$v_{n-1}$}}
\put(255,30){\line(1,0){40}}
\put(327,35){\mbox{$w_{n}$}}
\put(335,30){\circle*{3}}
\put(332,20){\mbox{$v_n$}}
\put(295,30){\line(1,0){40}}
\put(180,-5){$\text{Fig.}\, 2$}
\end{picture} \vspace{0.8cm}
 
If we denote by $d(\Gamma):=D_n(w_1,\dots,w_n)$, $d_{\text{\rm lt},i}(\Gamma):=D_{n_0-i-1}(w_1,\dots, w_{n_0-1-i})$,  
and  $d_{\text{rt},i}(\Gamma):=D_{n-n_0-i}(w_{n_0+1+i},\dots,w_n)$, then by (\ref{eq0+}), we have equality

\begin{equation}\label{eq0}  d(\Gamma)=  w_{n_0}d_{\text{\rm lt},0}(\Gamma)d_{\text{\rm rt},0}(\Gamma)-  
 d_{\text{\rm rt},0}(\Gamma)d_{\text{\rm lt,1}}(\Gamma)- d_{\text{\rm lt},0}(\Gamma)d_{\text{\rm rt},1}(\Gamma).\end{equation} 

\subsection{Weighted oriented graphs $\Gamma_{\frac{k}{q}}$ and $\Gamma_{\{ \frac{k_1}{q_1},\frac{k_2}{q_2}\}}$}
Let $(k,q)\in \mathbb N^2$ be a pair of coprime numbers, $k>q$, and 
\begin{equation}\label{con-fr}
\frac{k}{q}=w_1-\frac{1}{w_{2}-\frac{1}{w_{3}- \frac{1}{\dots -\frac{1}{w_{n(k,q)}}}}}
\end{equation}
expansion as continued fraction. Let us associate  a weighted oriented chain $\Gamma_{\frac{k}{q}}$ (depicted in Figure 1 in which $n=n(k,q)$  and weights $w_i$ are defined by (\ref{con-fr})). Note that 
\begin{equation}\label{fra} k= d(\Gamma_{\frac{k}{q}})=d_{\text{rt},0}(\Gamma_{\frac{k}{q}})\,\, \text{and}\,\, q=d_{\text{\rm rt},1}(\Gamma_{\frac{k}{q}}).\end{equation}

The following statement is obvious.
\begin{claim} \label{chai} Let $\Gamma$ be an oriented weighted graph depicted in Fig. 1 and whose all weights $w_i\geq 2$, and let 
$d(\Gamma)=k$, $d_{\text{ln},1}(\Gamma)=q$. Then $\Gamma$ and $\Gamma_{\frac{k}{q}}$ are isomorphic as weighted oriented chains.  
\end{claim}

Similarly, let us associate a graph $\Gamma_{\{ \frac{k_1}{q_1},\frac{k_2}{q_2}\}}$ with an orbit $\{ \frac{k_1}{q_1},\frac{k_2}{q_2}\}\in \mathcal O_{\mathcal D}$ as follows. If $\frac{k_2}{q_2}\neq \frac{1}{0}$, then $\Gamma_{\{ \frac{k_1}{q_1},\frac{k_2}{q_2}\}}=\Gamma^{\text{\rm inv}}_{\frac{k_1}{q_1}}\oplus \Gamma_{\{ \frac{1}{0},\frac{1}{0}\}}\oplus\Gamma_{\frac{k_2}{q_2}}$, where $\Gamma_{\{ \frac{1}{0},\frac{1}{0}\}}$ is the graph consisting of a single vertex whose weight $w=1$. If $\frac{k_2}{q_2}=\frac{1}{0}$ and $k_1>1$, then 
$\Gamma_{\{ \frac{k_1}{k_1-1},\frac{1}{0}\}}=\Gamma_{\frac{k_1}{k_1-1}}\oplus \Gamma_{\{ \frac{1}{0},\frac{1}{0}\}}$. 

If we put $n_0=\text{\rm n}(\Gamma_{\frac{k_1}{q_1}})+1$ and $n=\text{\rm n}(\Gamma_{\{ \frac{k_1}{q_1},\frac{k_2}{q_2}\}})=
\text{\rm n}(\Gamma_{\frac{k_1}{q_1}})+ \text{\rm n}(\Gamma_{\frac{k_2}{q_2}})+1$, then the graph $\Gamma_{\{ \frac{k_1}{q_1},\frac{k_2}{q_2}\}}$ is depicted in Fig. 2 in which the vertex $v_{n_0}$ has weight $w_{n_0}=1$. 

\begin{claim} \label{state} For  $\{ \frac{k_1}{q_1},\frac{k_2}{q_2}\}\in\mathcal O_{\mathcal D}$, we have equality
$$d_{\text{\rm lt},0}d_{\text{\rm rt},0}-d_{\text{\rm rt},0}d_{\text{\rm lt,1}}-   
 d_{\text{\rm lt},0}d_{\text{\rm rt},1}=1,$$
 where $d_{\text{\rm lt},i}:=d_{\text{\rm lt},i}(\Gamma_{\{ \frac{k_1}{q_1},\frac{k_2}{q_2}\}})$ and $d_{\text{\rm lt},i}:=d_{\text{\rm lt},i}(\Gamma_{\{ \frac{k_1}{q_1},\frac{k_2}{q_2}\}})$ for $i=0,1$. 
\end{claim}
\noindent {\it Proof} follows from equality (\ref{eq0}), equality (\ref{fra}), definition (\ref{not}), and definition  of the set $\mathcal D_{\mathcal P}$. \qed  

\subsection{Chains of curves} Let $\overline E_n=E_1\cup \dots \cup E_n\subset X$ be a chain of smooth  irreducible curves lying in  a smooth surface $X$, that is, the intersection numbers $(E_i,E_{i+1})_X=1$ for $i=1,\dots, n-1$ and $(E_i,E_j)_X=0$ if $|i-j|>1$. The  weighted graph $\Gamma(\overline E_n)$ of the curve $\overline E_n$ is depicted in Figure 1, where the vertices $v_i$ are in one-to-one correspondence with the curves $E_i$ and they have weights $w_i=-(E_i^2)_X$. Put $d(\overline E_n):=d(\Gamma(\overline E_n))$.

Let $I(\overline E_n)=((E_i,E_j)_X)$ be the intersection matrix of the irreducible components of $\overline E_n$.  
Then $I(\overline E_n)=M^-_{n,+}$, where $M_{n,+}=M(\Gamma(\overline E_n))$, and 
\begin{equation}\label{eq6} \det I(\overline E_n)=(-1)^{n}d(\Gamma(\overline E_n)). 
\end{equation}

Let $V\subset X$ be a "tubular" neighborhood of $\overline E_n$ (here we use definition of "tubular" neighborhoods given in \cite{K-L}).
\begin{prop}\label{prop3} {\rm (Theorem 4 in \cite{K-L})} If $\overline E_n$ is a chain of rational curves, then 
$$\pi_1(V\setminus \overline E_n)=\mathbb Z/d(\overline E_n)\mathbb Z=\langle x_1\,\,|\,\, x_1^{d(\overline E_n)}=1\rangle=\langle x_{n}\,\,|\,\, x_{n}^{d(\overline E_n)}=1\rangle,$$
where the element $x_1$ $($resp., $x_{n})$ is represented by a simple loop around the curve $E_1$ $($resp., $E_{n})$. \end{prop} 

\section{ Resolution of of singular points of $\mathcal P$-simplest singularity type}\label{sect4}
\subsection{Chains of $\sigma$-processes.} \label{sect-4.1} Let  
$\overline{\sigma}_n =\sigma_1\circ \dots \circ \sigma_n: (V_n,p_n)\to (V_0,p_0)$ be a sequence of $n$ $\sigma$-processes $\sigma_i:(V_i,p_i)\to (V_{i-1},p_{i-1})$ with centers at the points $p_{i-1}$. For $1\leq i\leq j\leq n$ denote by $E_i\subset V_j$ the proper inverse image in $V_j$ of the exceptional curve of $\sigma$-process $\sigma_i:V_i\to V_{i-1}$ and $\overline E_j=E_1\cup\dots \cup E_j\subset V_j$.  

A sequence $\overline{\sigma}_n: V_n\to V_0$ is called a  {\it  chain of $\sigma$-processes} if it satisfies the following conditions:
\begin{itemize}
\item[(1)] for  $1\leq j\leq n-1$, the center $p_j\in E_j\subset V_j$, 
\item[(2)] for  $3\leq j\leq n-1$, if the curve $E_j\subset V_j$ intersects with two other curves $E_{i_1}$ and $E_{i_2}$, $1\leq i_1,i_2<j$, 
then $p_j\in V_j$ is either $E_j\cap E_{i_1}$ or $E_j\cap E_{i_2}$. 
\end{itemize}  

A chain of $\sigma$-processes $\overline{\sigma}_n: V_n\to V_0$, $n\geq 3$, is {\it non-degenerate} if it satisfies the following additional condition:
\begin{itemize}
\item[($n_1$)] $E_n$ intersects in $V_n$ with some curve $E_{i}$, were $i\leq n-2$.
\end{itemize}

It is easy to see that a chain $\overline{\sigma}_n: V_n\to V_0$ of $\sigma$-processes is  {\it degenerate} (i.e., not non-degenerate) if and only if 
\begin{itemize}
\item[$(\text{\bf \rm d}_1$)] 
$(E_i^2)_{V_n}=-2$ for $1\leq i\leq n-1$ and $(E_n^2)_{V_n}=-1$.   
\end{itemize}

The  weighted graph $\Gamma(\overline E_n)$ of the exceptional curve $\overline E_n \subset V_n$ is the chain depicted in Figure 2 
for which there is a bijection $\beta_E:\{ 1,\dots, n\}\to \{ 1,\dots, n\}$ satisfying the following conditions:  

\begin{itemize} \item[$(\text{\bf \rm c}_1$)] the vertices $v_j\in\Gamma(\overline E_n)$ correspond to the curves $E_{\beta_E(j)}\subset \overline E_n\subset V_n$,  
\item[$(\text{\bf \rm c}_2$)]  $\beta_E(1)=1$ and   $\beta_E^{-1}(n)=n_0$, 
\item[$(\text{\bf \rm c}_3$)] the weight  $w_{n_0}=-(E_{\beta_E(n_0)}^2)_{V_n}=1$,
\item[$(\text{\bf \rm c}_4$)] for $1\leq j\leq n-1$ the curves $E_{\beta_E({j})}$ and $E_{\beta_E({j+1})}$ have a common point,
\item[$(\text{\bf \rm c}_5$)] for  $j\neq n_0$, the weights $w_j=(E_{\beta_E(j)}^2)_{V_n}\geq 2$, 
\item[$(\text{\bf \rm c}_6)$] if $n\neq n_0$ then either $w_{n_0-1}=2$ and $w_{n_0+1}>2$ or $w_{n_0-1}>2$ and $w_{n_0+1}=2$, 
\item[$(\text{\bf \rm c}_7)$] if $n=n_0$ then $w_{n_j}=2$ for $j=1,\dots, n-1$.
\end{itemize}

The intersection matrix  $I(\overline E_n)=((E_{\beta_E(i)},E_{\beta_E(j)})_{V_n})$ of the irreducible components of the exceptional curve $\overline E_n$ of a chain  $\overline{\sigma}_n$ of $\sigma$-processes is negative definite and $\det I(\overline E_n)=(-1)^n$. The converse statement is also true. Namely, by Mumford's criterion (\cite{Mu}), if the intersection matrix  $I(\overline E_n)=((E_{\beta_E(i)},E_{\beta_E(j)})_{V_n})$ of the irreducible components of a chain $\overline E_n$ of rational curves is negative definite and $|\det I(\overline E_n)|=1$, then  $\overline E_n$ is an exceptional curve of a chain  $\overline{\sigma}_n$ of $\sigma$-processes. 

We call $E_{n,\text{lt}}=\cup_{i=1}^{n_0-1}E_{\beta_E(i)}\subset  V_n$ the {\it left part} and $E_{n,\text{rt}}=\cup_{i=n_0+1}^nE_{\beta_E(i)}\subset V_n$ the {\it right part} of $\overline E_n\subset V_n$. If $\overline{\sigma}_n$ is a non-degenerate chain, then $n_0<n$ and   
$$\Gamma(\overline E_n)=\Gamma(E_{n,\text{lt}})\oplus \Gamma_{\{ \frac{1}{0},\frac{1}{0}\}}\oplus \Gamma(E_{n,\text{rt}})=
\Gamma^{inv}_{\frac{k_1}{q_1}}\oplus \Gamma_{\{ \frac{1}{0},\frac{1}{0}\}}\oplus \Gamma_{\frac{k_2}{q_2}},$$ 
where $k_1=d_{\text{lt},0}(\overline E_n)$, $q_1=d_{\text{lt},1}(\overline E_n)$ and $k_2=d_{\text{rt},0}(\overline E_n)$, 
$q_2=d_{\text{rt},1}(\overline E_n)$.

If $\overline{\sigma}_n$ is a degenerate chain and $n>1$, then $n_0=n$ and 
$$\Gamma(\overline E_n)=\Gamma(E_{n,\text{lt}})\oplus \Gamma_{\{ \frac{1}{0},\frac{1}{0}\}}=
\Gamma^{inv}_{\frac{k_1}{q_1}}\oplus \Gamma_{\{ \frac{1}{0},\frac{1}{0}\}},$$ 
where $k_1=d_{\text{lt},0}(\overline E_n)$, $q_1=d_{\text{lt},1}(\overline E_n)$. Finally, if $n=1$, then 
$\Gamma(\overline E_1)=\Gamma_{\{ \frac{1}{0},\frac{1}{0}\}}$.

Define a map $\overline s:\mathcal E\to \mathbb Z^4$ from the set $\mathcal E$ of the chains $\overline{\sigma}_n$ of $\sigma$-processes as follows: 
$\overline s_{\overline{\sigma}_n}=(d_{\text{lt},0}(\overline E_n),d_{\text{rt},0}(\overline E_n),d_{\text{lt},1}(\overline E_n),
d_{\text{rt},1}(\overline E_n))$ if $\overline{\sigma}_n$ is a non-degenerate chain, 
$\overline s_{\overline{\sigma}_n}=(d_{\text{lt},0}(\overline E_n),1,d_{\text{lt},1}(\overline E_n),0)$ if $\overline{\sigma}_n$ is a degenerate chain, $n>1$, and $\overline s_{\overline{\sigma}_n}=(1,1,0,0)$ if $n=1$.

\begin{claim} \label{prop4} The images $\overline s_{\overline{\sigma}_n}\in \mathcal D_{\mathcal P}$ for all $\overline{\sigma}_n\in\mathcal E$.\end{claim}
\proof The intersection matrix  $I(\overline E_n)$ is negative definite  and 
$\det I(\overline E_n)=(-1)^n$. It follows from Lemma \ref{pos-def} that the matrix $M(\Gamma(\overline E_n))$ is positive definite and  $d(\overline E_n)=D(\Gamma(\overline E_n))=1$. Therefore Statement \ref{prop4} follows from equality (\ref{eq0}). \qed \\

Using notation used in the beginning of this Subsection, for $i=1,\dots,n$, denote $p_{i,\text{\rm lt}}=E_{i,\text{\rm ln}}\cap E_i\in V_i$ and $p_{\text{i,\rm rt}}=E_{\text{i,\rm rn}}\cap E_i\in V_i$, where $E_i$ is the exceptional curve of $\sigma_i:V_i\to V_{i-1}$ (if $\overline E_i\subset V_i$ is a degenerate chain, then $p_{\text{i,\rm rt}}$ is a point in $E_i$ not coinciding with $p_{i,\text{\rm lt}}$). 

\begin{prop} \label{left}  Let $\overline E_n\subset V_n$ be the exceptional curve of a chain $\overline{\sigma}_n$ of $\sigma$-processes, $n\geq 2$.  Then 
$d_{\text{\rm lt},0}(\overline E_n)>d_{\text{\rm rt},0}(\overline E_n)$ and therefore $\{ \frac{d_{\text{\rm lt},0}(\overline E_n)}{d_{\text{\rm lt},1}(\overline E_n)},\frac{d_{\text{\rm rt},0}(\overline E_n)}{d_{\text{\rm rt},1}(\overline E_n)}\}\in \mathcal O$. 

Conversely, for each $\{ \frac{k_1}{q_1},\frac{k_2}{q_2}\}\in \mathcal O$ there is a chain $\overline{\sigma}_n$ of $\sigma$-processes such that 
$k_1=d_{\text{\rm lt},0}(\overline E_n)$, $q_1=d_{\text{\rm lt},1}(\overline E_n)$ and $k_2=d_{\text{\rm rt},0}(\overline E_n)$, 
$q_2=d_{\text{\rm rt},1}(\overline E_n)$.
\end{prop} 

\proof Assume that for $m\leq n$ inequality $d_{\text{\rm lt},0}(\overline E_m)>d_{\text{\rm rt},0}(\overline E_m)$ is correct.  Consider a chain $\overline{\sigma}_{n+1}$ of $\sigma$-processes, $\overline{\sigma}_{n+1}=\overline{\sigma}_n\circ\sigma_{n+1}:V_{n+1}\to V_0$, where $\sigma_{n+1}:V_{n+1}\to V_n$ and $\overline{\sigma}_n: V_n\to V_0$. There are two cases:  $\sigma_{n+1}$ is either the blowup of the point $p_{n,lt}=E_{n,lt}\cap E_n\subset V_n$, or the blowup of the point $p_{n,rt}=E_{n,rt}\cap E_n\subset V_n$.   \

Let's compare graphs $\Gamma(E_{n,\text{\rm lt}})$ and $\Gamma(E_{n,\text{\rm rt}})$ of the chain  $\overline E_n$ with graphs $\Gamma(E_{n+1,\text{\rm lt}})$ $\Gamma(E_{n+1,\text{\rm rt}})$ of the chain $\overline E_{n+1}$. It is easy to see that in the first case  
the graphs $\Gamma(E_{n,\text{\rm lt}})$ and $\Gamma(E_{n+1,\text{\rm lt}})$ are isomorphic as oriented (non-weighted) graphs. The weights $w_i$ of the vertices $v_i$ of $\Gamma(E_{n,\text{\rm lt}})$ and the weights $\widetilde w_i$ of the vertices $v_i$ of $\Gamma(E_{n+1,\text{\rm lt}})$ are equal to each other for $i=1,\dots, n_0-2$ and $\widetilde w_{n_0-1}=w_{n_0-1}+1$. Therefore 
\begin{equation}\label{d_0-1} \begin{array}{ll} d_{\text{\rm lt},0}(\overline E_{n+1})= & \widetilde w_{n_0-1}d_{\text{\rm lt},1}(\overline E_n)-d_{\text{\rm lt},2}(\overline E_n)= (w_{n_0-1}+1)d_{\text{\rm lt},1}(\overline E_n)-d_{\text{\rm lt},2}(\overline E_n)=  \\ 
& (w_{n_0-1}d_{\text{\rm lt},1}(\overline E_n)-d_{\text{\rm lt},2}(\overline E_n))+d_{\text{\rm lt},1}(\overline E_n)=
 d_{\text{\rm lt},0}(\overline E_n)+d_{\text{\rm lt},1}(\overline E_n).\end{array}\end{equation} 

The graph $\Gamma(E_{n+1,\text{\rm rt}})=\Gamma(E_{n+1,\text{\rm rt}})=\Gamma(2)\oplus \Gamma(E_{n,\text{\rm rt}})$, where $\Gamma(2)$ consists of a single vertex with weight $2$. Therefore 
$$d_{\text{\rm rt},0}(\overline E_{n+1})=2d_{\text{\rm rt},0}(\overline E_{n})-d_{\text{\rm rt},1}(\overline E_{n})=
d_{\text{\rm rt},0}(\overline E_{n})+(d_{\text{\rm rt},0}(\overline E_{n})-d_{\text{\rm rt},1}(\overline E_{n})).$$
Consequently, 
$$d_{\text{\rm lt},0}(\overline E_{n+1})= d_{\text{\rm lt},0}(\overline E_n)+d_{\text{\rm lt},1}(\overline E_n)>
d_{\text{\rm rt},0}(\overline E_{n})+(d_{\text{\rm rt},0}(\overline E_{n})-d_{\text{\rm rt},1}(\overline E_{n}))=d_{\text{\rm rt},0}(\overline E_{n+1}),$$
since $d_{\text{\rm lt},0}(\overline E_{n})>d_{\text{\rm rt},0}(\overline E_{n})$ by assumption and 
$d_{\text{\rm rt},0}(\overline E_{n})\leq 
d_{\text{\rm lt},1}(\overline E_{n})+d_{\text{\rm rt},1}(\overline E_{n})$ by Lemma \ref{lem-ineq}.

In the second case similar calculations give the inequality 
$$d_{\text{\rm rt},0}(\overline E_{n+1})= d_{\text{\rm rt},0}(\overline E_n)+d_{\text{\rm rt},1}(\overline E_n)<
d_{\text{\rm lt},0}(\overline E_{n})+(d_{\text{\rm lt},0}(\overline E_{n})-d_{\text{\rm lt},1}(\overline E_{n}))=d_{\text{\rm lt},0}(\overline E_{n+1}),$$
since $d_{\text{\rm rt},0}(\overline E_{n})<d_{\text{\rm lt},0}(\overline E_{n})$ by assumption and 
$d_{\text{\rm lt},0}(\overline E_{n})> 
d_{\text{\rm rt},1}(\overline E_{n})+d_{\text{\rm lt},1}(\overline E_{n})$ by Lemma \ref{lem-ineq}. 

The converse statement directly follows from Statement \ref{state} and Mumford's criterion. \qed

\subsection{Chains of $\sigma$-processes resolving $\mathcal P$-simplest singularities } \label{sect-4.2}
By definition, a sequence $\overline{\sigma}_n$ of $\sigma$-processes resolves the singular point of a curve germ $(B_0,p_0)$ 
if $\overline{\sigma}_n^{-1}(B_0)$ is a divisor with normal crossing and the number $n>0$ of blowups  included in the sequence   
$\overline{\sigma}_n$ is minimal. 

\begin{rem} \label{rem-s} The reason why we consider a curve germ given by equation $x^k-y=0$, where $k\in \mathbb N$, to be singular  is 
that at the intermediate steps of blowups, the proper inverse image $B_m$ of the curve germ $B_0$ may no longer be singular, but $\overline{\sigma}_m^{-1}(B_0)$ is not a divisor with normal crossings and we must continue to blowup points until we get a divisor with normal crossings.
Therefore, in particular, in the case of singularity type $\mathcal P_{1,1}$ of the "singular"\, curve given by equation $x-y=0$, by definition we must perform one blowup. \end{rem} 

Denote by $\mathcal B$ the set of pairs $(B_{k_1,k_2},\overline{\sigma}_n)$ in which $\overline{\sigma}_n$ is the chain of $\sigma$-processes resolving the singular point of a curve $B_{k_1, k_2}$ given by equation $x^{k_1}-y^{k_2}=0$, where $\{ k_1,k_2\}\in\mathcal O$. We associate with a pair $(B_{k_1,k_2},\overline{\sigma}_n)$ the path $\lambda((B_{k_1,k_2},\overline{\sigma}_n)):=l(\overline s(\overline{\sigma}_n))$ 
in the graph $\Gamma_{\overline{\mathcal O}_{\mathcal D}}$ connecting the vertex $\overline s_{\overline{\sigma}_n}\in 
\Gamma_{\overline{\mathcal O}_{\mathcal D}}$ with the root of  $\Gamma_{\overline{\mathcal O}_{\mathcal D}}$.

\begin{thm}\label{sigma-p} Let a chain  $\overline{\sigma}_n=\sigma_1\circ\dots\circ\sigma_n$ of $\sigma$-processes resolves the singular point of the curve $B_{k_1,k_2}$ given by equation $x^{k_1}-y^{k_2}=0$, where $\{ k_1,k_2\}\in\mathcal O$.Then the correspondence $\lambda: (B_{k_1,k_2},\overline{\sigma}_n))\mapsto l(\overline s(\overline{\sigma}_n))$ has the following properties:
\begin{itemize} 
\item[$(\bf 1)$]  $\overline s_{\overline{\sigma}_n}=\text{\rm pr}^{-1}_{\Gamma_{\overline{\mathcal O}_{\mathcal D}},\Gamma_{\overline{\mathcal O}}}
(\{ k_1,k_2\})= \text{\rm pr}^{-1}(\{ k_1,k_2\})$,
\item[$(\bf 2)$] $\text{\rm ln}(\overline s_{\overline{\sigma}_n})=n_{\mathcal E}(k_1,k_2)+1=n$ and if we number sequentially the edges included in the path $l(\overline s_{\overline{\sigma}_n})$ connecting the vertex $\overline s_{\overline{\sigma}_n}$ with the root of the graph $\Gamma_{\overline{\mathcal O}_{\mathcal D}}$, in the order starting from the edge whose origin is the vertex $\overline s_{\overline{\sigma}_n}$, then for each $i=1,\dots,n$ 
    the $\sigma$-process $\sigma_i$ incoming in $\overline{\sigma}_n$ is in on-to-one correspondence with the edge incoming in the path  
    $l(\overline s_{\overline{\sigma}_n})$ with number $i$.  
\end{itemize}\end{thm}
\proof  
Let $(x_0,y_0)$ be coordinates in $V_0$ and a curve germ $B_0=B_{k_{1,0},k_{2,0}}\subset V_0$ is given by equation $x_0^{k_{1,0}}-y_0^{k_{2,0}}=0$,
where $\{ k_{1,0},k_{2,0}\}\in \mathcal O$ (and hence, by definition, $k_{1,0}\geq k_{2,0}$). 

By definition, the $\sigma$-process with center at $p_0=\{ x_0=y_0=0\}\in V_0$ is a map 
$$\sigma :V_1=\{ (v,(z_1:z_2))\in V_0\times\mathbb P^1\mid x_0(v)z_2=y_0(v)z_1\}\to V_0.$$ 
The inverse image $E_1=\sigma^{-1}(p_0)=\{ p_0\}\times \mathbb P^1\subset V_1$. The surface $V_1$ is covered by two neighborhoods $V_{1,z_1}=\{ z_1\neq 0\}$ and $V_{1,z_2}=\{ z_2\neq 0\}$.   Functions $x_{1,1}=x_0$ and $y_{1,1}=\frac {z_2}{z_1}$ are coordinates in $V_{1,z_1}$, the restriction of $\sigma$ to $V_{1,z_1}$ is given by functions $x_0=x_{1,1}$, $y_0=x_{1,1}y_{1,1}$, and the intersection $E_1\cap V_{1,z_1}$ is given by equation $x_{1,1}=0$. Similarly, functions $y_{2,1}=y_0$ and $x_{2,1}=\frac {z_1}{z_2}$ are coordinates in $V_{1,z_2}$, the restriction of $\sigma$ to $V_{1,z_2}$ is given by functions $x_0=x_{2,1}y_{2,1}$, $y_0=y_{2,1}$, and the intersection $E_1\cap V_{z_2}$ is given by equation $y_{2,1}=0$. 

The proper inverse image $B_1=\sigma_1^{-1}(B_0)$ lies in $V_{1,z_1}$ and it is given in $V_{1,z_1}$ by equation $x_{1,1}^{k_{1,0}-k_{2,0}}-y_{1,1}^{k_{2,0}}=0$.
 
Let's accept the agreement that at each blowup $\sigma_{i+1}: V_{i+1}\to V_i$ in the chain $\overline{\sigma}_n=\sigma_1\circ\dots\circ\sigma_n$ of $\sigma$-processes resolving the singular point of the curve germ $B_0$, the proper inverse image $B_{i+1}$ of the curve $B_i$  is given by equation $x_{i+1}^{k_{1,i+1}}-y_{i+1}^{k_{2,i+1}}=0$, where $\{ k_{1,i+1},k_{2,i+1}\}\in\mathcal O$. Therefore if $k_{1,0}<2k_{2,0}$, then before to continue to resolve the singular point of $B_1$, we must perform the coordinate change $x_{1,1}\leftrightarrow y_{1,1}$. 

Consequently, if we identify the $\sigma$-process $\sigma_1$ with the edge whose origin is the vertex $\overline s_{\overline{\sigma}_n}$, and use induction by  $n$, assuming that for $m\leq n-1$ the statements of Theorem hold, then we get that statements of Theorem holds for all $n\in\mathbb N$. Note also that an edge with number $i+1$ incoming in the path $l(\overline s_{\overline{\sigma}_n})$ is marked by transformation $\varepsilon_1$ if $B_i$ is given by equation $x_{1,i}^{k_{1,i}}-y_{1,i}^{k_{2,i}}=0$, where $k_{1,i}\geq 2k_{2,i}$, and it is marked by transformation $\varepsilon_2$ if $B_i$ is given by equation $x_{1,i}^{k_{1,i}}-y_{1,i}^{k_{2,i}}=0$, where $k_{1,i}<2k_{2,i}$.    \qed \\  

Let $\overline{\sigma}_n$ be a chain of $\sigma$-processes resolving the singular point of the curve $B_{k_1, k_2}$ given by equation $x^{k_1}-y^{k_2}=0$, 
where $\{ k_1,k_2\}\in\mathcal O$. We associate with a pair $(B_{k_1,k_2},\overline{\sigma}_n)$ a {\it partially weighted graph} $\Gamma(B_{k_1,k_2})$ depicted in Fig. 3, 

\begin{picture}(300,85)
\put(45,30){\circle*{3}}
\put(40,35){$\mbox{$w$}_{1}$}
\put(42,20){$\mbox{$v$}_{1}$}
\put(45,30){\line(1,0){40}}
\put(85,30){\circle*{3}}
\put(80,35){$\mbox{$w$}_{2}$}
\put(82,20){$\mbox{$v$}_{2}$}
\put(85,30){\line(1,0){40}}
\put(132,30){$\dots$}
\put(150,30){\line(1,0){40}}
\put(193,33){$\mbox{\tiny $-1$}$}
\put(190,70){\circle*{3}}
\put(180,67){$\mbox{$b$}$}
\put(190,70){\line(0,-1){40}}
\put(190,30){\circle*{3}}
\put(187,20){$\mbox{$v_{n_0}$}$}
\put(192,30){\line(1,0){40}}
\put(235,30){$\dots$}
\put(335,30){\circle*{2}}
\put(287,35){$\mbox{$w_{n-1}$}$}
\put(295,30){\circle*{3}}
\put(290,20){\mbox{$v_{n-1}$}}
\put(255,30){\line(1,0){40}}
\put(327,35){\mbox{$w_{n}$}}
\put(335,30){\circle*{3}}
\put(332,20){\mbox{$v_n$}}
\put(295,30){\line(1,0){40}}
\put(180,-5){$\text{Fig.}\, 3$}
\end{picture} \vspace{0.5cm}

\noindent   where the vertex $b$ corresponds to the proper inverse image $B_n\subset \overline{\sigma}_n^{-1}(B_{k_1,k_2})$ of the curve $B_{k_1,k_2}$, this vertex has not a weight, and $\Gamma(B_{k_1,k_2})\setminus \{ b\}=\Gamma(\overline E_n)$ is the oriented chain.

The following proposition is well known.
\begin{prop} {\rm (}\cite{W}{\rm )} Two curve germs $(B_1,o)$ and $(B_2,o)$ are equisingular deformation  equivalent if and only if their graphs $\Gamma(B_1)$ and $\Gamma(B_2)$ are isomorphic as partially weighted graphs.
\end{prop}

\noindent In particular, if the singularity of the branch curve germ $B_P$ of a map $f_P:U_P\to V_p$ is $\mathcal P_{k_1,k_2}$-simplest, then $\Gamma(B_P)=\Gamma(B_{k_1,k_2})$.

\subsection{Local fundamental  groups of curve germs having $\mathcal P$-simplest singula\-ri\-ties}\label{fund}
Consider a curve germ $B_P\subset V_p$ 
having at $p\in V_p$ a $\mathcal P_{k_1,k_2}$-simplest singularity the graph $\Gamma(B_P)$ of which is depicted in Fig. 3. 
The map $\overline{\sigma}:V_n\to V_p$ defines an isomorphism 
$$\overline{\sigma}_*: \pi_1(V_n\setminus \overline{\sigma}_n^{-1}(B_P))\to \pi_1(V_p\setminus B_P)=\pi_1^{loc}(B_P,p)$$
that allows us to identify $\pi_1^{loc}(B_P,p)$ with $\pi_1(V_n\setminus \overline{\sigma}_n^{-1}(B_P))$.

It follows from Theorem 3 in \cite{K-L} that the local fundamental group $\pi_1^{loc}(B_P,p)$ has the following presentation. 
The group $\pi_1^{loc}(B_P,p)$ is generated by $n+1$ elements $x_0,x_1,\dots,x_n$, which are geometric generators  that are in one-to-one correspondence with the vertices $b,v_1,\dots, v_n$ of the graph $\Gamma(B_P)$: $x_0$ is represented by a loop that is a simple circuit around the proper inverse image 
$B\subset V_n$ of the curve germ $B_P$ and $x_i$, $1\leq i\leq n$ are represented by  loops that are simple circuits around the curves $E_{\beta^{-1}(i)}\subset \overline E_n$. The  elements $x_0,x_1,\dots,x_n$ are subject to the following defining relations:
\begin{equation} \label{rel1} [x_1,x_2]=\dots =[x_i,x_{i+1}]=\dots =[x_{n-1},x_n]=[x_0x_{k_0}]=1, \phantom{aaaaaaa} \end{equation} 
\begin{equation}\label{rel2} x_{n_0}^{-1}x_{n_0-1}x_0x_{n_0+1}=1,\phantom{aaaaaaaaa} \end{equation}
\begin{equation}\label{rel3} x_1^{-w_1}x_2=x_{n-1}x_n^{-w_n}=1, \phantom{aaaaaaaa} \end{equation} 
\begin{equation}\label{rel4} \phantom{aaaaaaaaaaaaaaaaa} x_{i-1}x_i^{-w_i}x_{i+1}=1\qquad\quad \text{for}\,\, 2\leq i\leq n-1,\,\, i\neq n_0,  \end{equation}    
were $w_i$ are the weights of vertices $v_i$ of $\Gamma(B_P)$. 
 \begin{prop}\label{prop7} The group $\pi_1(V_p\setminus B_P)=\pi_1(V_n\setminus \overline{\sigma}_n^{-1}(B_P)))$ is generated by elements $x_0$, $x_{n_0-1}$, and $x_{n_0+1}$being subject to the following relations: 
\begin{equation} \label{7-0}
[x_{0},x_{n_0-1}x_0x_{n_0+1}]=[x_{n_0-1},x_{n_0-1}x_0x_{n_0+1}]=[x_{n_0+1},x_{n_0-1}x_0x_{n_0+1}]=1,
\end{equation}
\begin{equation} \label{7-1}
(x_{n_0-1}x_0x_{n_0+1})^{-d_{\text{\rm lt},1}(E_{\text{\rm lt}})}x_{n_0-1}^{d(E_{\text{\rm lt}})}=
(x_{n_0-1}x_0x_{n_0+1})^{-d_{\text{\rm lt},1}(E_{\text{\rm rt}})}x_{n_0+1}^{d(E_{\text{\rm rt}})}=1.
\end{equation}

The element $x_{n_0}=x_{n_0-1}x_0x_{n_0+1}$ belongs to the center of $\pi_1(V_n\setminus B_P)$ and
\begin{equation}\label{7-4}  x_{n_0-i} =x_{n_0-1}^{(-1)^{i-1}D_{i-1}(w_{n_0-i+1},\dots,w_{n_0-1})}x_{n_0}^{(-1)^{i-1}D_{i-2}(w_{n_0-i+1},\dots,w_{n_0-2})} 
\end{equation}
{for} $2\leq i\leq n_0-1$,
\begin{equation}\label{7-5} 
x_{n_0+i} =x_{n_0+1}^{(-1)^{i-1}D_{i-1}(w_{n_0-1},\dots,w_{n_0+1-i})}x_{n_0}^{(-1)^{i-1}D_{i-2}(w_{n_0-1},\dots,w_{n_0+1-i})} 
\end{equation}
{for} $2\leq i\leq n-n_0$. \end{prop} 
\proof By (\ref{rel2}), $x_{n_0}=x_{n_0-1}x_0x_{n_0+1}$ and relations (\ref{7-0}) follow from relations (\ref{rel1}). Therefore $x_{n_0}$ belongs to the center of the group generated by $x_{n_0-1},x_0$, and $x_{n_0+1}$. 

Equalities (\ref{7-4}) and (\ref{7-5}) are obtained recurrently using relations (\ref{rel4}). We obtain relations (\ref{7-1}) if we substitute in (\ref{rel3}) instead of the letters $x_1$, $x_2$ the words $W_1(x_{n_0-1},x_{n_0})$, $W_2(x_{n_0-1},x_{n_0})$ standing in the right side of equality (\ref{7-4}) ($i=n_0-1$ and $i=n_0-2$) and instead of the letters $x_{n-1}$, $x_n$ the words $W_{n-1}(x_{n_0+1},x_{n_0})$, $W_{n}(x_{n_0+1},x_{n_0})$ standing in the right side of equality (\ref{7-5}) ($i=n-n_0-1$ and $i=n-n_0$). 

\subsection{Cyclic quotient singularities} \label{cyclics}
In this subsection, properties of cyclic quotient singularities are reminded, proofs of which are contained, for example, in \cite{Ba}. 

Denote by $h:U=\mathbb C^2\to T=\mathbb C^2$ a finite morphism, $\deg h=k^2$, given by functions $t_i=u_i^k$, $i=1,2$, where $t_1,t_2$ (resp., $u_1,u_2$) are coordinates in $T$ (resp., $U$). The morphism $h$ is branched in the $t_i$-axises $C_{t_i}=\{t_i=0\}\subset T$ and it is ramified with multiplicity $k$ in  the $u_i$-axises $C_{u_i}=\{ u_i=0\}\subset U$, $i=1,2$. Denote by $o_{t}=\{ t_1=t_2=0\}\in T$ and $o:=o_u=\{u_1=u_2=0\}\in U$ the origins of coordinate systems in $T$ and $U$.  

The morphism $h$ defines the imbedding $h^*: \mathbb C[t_1,t_2]\hookrightarrow \mathbb C[u_1,u_2]$ of the polynomial ring $\mathbb C[t_1,t_2]$ and the imbedding $h^*: \mathbb C(t_1,t_2)\hookrightarrow \mathbb C(u_1,u_2)$ of the field of rational functions $\mathbb C(u_1,u_2)$ that is a Galois expansion of the field $C(t_1,t_2)$ whose  Galois group $\overline G=\text{\rm Gal}(U/T)=\mathbb Z_k\times\mathbb Z_k$ is an abelian group 
generated by $g_1=(1,0)$ and $g_2=(0,1)$ acting as follows:
$$ g_1: (u_1,u_2)\mapsto (e^{\frac{2\pi i}{k}}u_1,u_2),\quad  g_2: (u_1,u_2)\mapsto (u_1,e^{\frac{2\pi i}{k}}u_2).$$ 
We consider the surface $T=U/\overline G$ as the quotient of $U$ under the action of $\overline G$ on $U$.

Denote by $G_{a,b}$ a cyclic subgroup of $\overline G$ generated by $g_{a,b}=g_1^ag_2^b$  and by $V_{a,b}=U/G_{a,b}$ the quotient of $U$ under the action of $G_{a,b}$ on $U$. The surface $V_{a,b}$ is a normal variety and its ring of regular functions is the integral closure of $\mathbb C[t_1,t_2]$ in the field 
$$\mathcal K_{a,b}=\mathbb C(u_1,u_2)^{G_{a,b}}=\{ f\in \mathbb C(u_1,u_2)\mid g(f)=f\,\, \text{\rm for}\,\, g\in G_{a,b}\}.$$

In what follows (unless otherwise agreed) we will assume that either $ab=0$,  or $1\leq |a|,|b|\leq k-1$ and $a$, $b$, and $k$ are pairwise coprime integers.   

The surface $V_{a,b}$ is non-singular if and only if $ab=0$ and in this case, say $a=0$, the quotient map $h_{0,b}: U\to V_{0,b}=\mathbb C^2$ is given by functions  $v_1=u_1$, $v_2=u_2^{\frac{k}{m}}$, where $v_1,v_2$ are coordinates in $V_{0,b}$ and $m=\text{g.c.d.}(k,b)$. The degree of $h_{0,b}$ is equal to $\frac{k}{m}$ and $h_{0,b}$ is ramified in $C_{u_2}$ with multiplicity $\frac{k}{m}$.

Consider the case when $a=q\in\mathbb N$ and $b=1$, where $k$ and $q$ are  coprime. The cyclic  subgroup $G_k:=G_{q,1}\subset \overline G$ of order $k$ is generated by $g_{q,1}=g^q_1g_2\in \overline G$ that acts on $U$ as follows: $ g_q: (u_1,u_2)\mapsto (e^{\frac{2\pi qi}{k}}u_1,e^{\frac{2\pi i}{k}}u_2)$. It is easy to see that $\mathcal K_k=\mathbb C(u_1,u_2)^{G_{k}}\subset\mathbb C(u_1,u_2)$ is a Galois extension of $\mathbb C(t_1,t_2)$ and it is generated over $\mathbb C(t_1,t_2)$ by the function $v=u_1u^{k-q}_2$ that is a root of equation $v^k=t_1t_2^{k-q}$. In other words, $\mathcal K_k$ is the field of rational functions of a surface $\widetilde V\subset \mathbb C^3$ given by equation $v^k=t_1t_2^{k-q}$. The imbedding $\mathcal K_k\hookrightarrow\mathbb C(u_1,u_2)$ defines  a morphism $\widetilde h_{q,1}:U\to \widetilde V\subset \mathbb C^3$ given by functions $t_1=u_1^k$, $t_2=u_2^k$, $v=u_1u^{k-q}_2$.  

Let $\nu: V\to \widetilde V$ be the normalization of the surface $\widetilde V$. Then there is a factorization 
$$ h: U\stackrel{h_{q}}{\longrightarrow} V\stackrel{\nu}{\longrightarrow}\widetilde V \stackrel{\text{\rm pr}}{\longrightarrow} T,$$ 
where $\nu \circ h_q=\widetilde h_{q,1}$ and $\text{\rm pr}: \widetilde V\to T$ is the restriction to $\widetilde V$ of the projection 
$\text{\rm pr}:(t_1,t_2,v)\mapsto (t_1,t_2)$. 
The point $o_v=h_q(o)$ is a single singular point of the surface $V$ and $h_q$ is not ramified in $C_{u_1}$ and $C_{u_2}$, and it is ramified over $o_v$ only.
The surface $U\setminus \{ o\}$ is simply connected and $h_q:U\setminus \{ o\}\to V\setminus \{ o_v\}$ is an unramified cover. Therefore $h_q:U\setminus \{ o\}\to V\setminus \{ o_v\}$ is the universal cover of the surface $V\setminus \{ o_v\}$ and $G_k=G_{q,1}=\pi_1(V\setminus \{ o_v\})$. 

Let a quasi-projective normal variety $\overline V=\cup\overline V_i$ whose field of rational functions $\mathbb C(\overline V)$ is $\mathcal K_k$,  be covered by affine varieties $\overline V_i$. It is not difficult to show that the integral closures $R_i$ of the rings of regular functions 
$\mathbb C[\overline V_i]$ of varieties $\overline V_i$ in the field $\mathbb C(u_1,u_2)$ define a normal quasi-projective variety $\overline U$ covered by affine normal varieties $\overline U_i=\text{\rm Spec}\, R_i$ such that the group $G_k\subset \text{\rm Aut}(\mathbb C(u_1,u_2)/\mathbb C)$ defines an action on $\overline U$ such that $\overline V=\overline U/G_k$. The converse statement is also true. Namely, if a quasi-projective variety 
$\overline U =\cup \overline U_i$ is covered by affine normal varieties $\overline U_i$, such that the rings or regular functions 
$\mathbb C[\overline U_i]\subset \mathbb C(u_1,u_2)$ are invariant under the action of the group $G_k$, then the group $G_k$ acts on $\overline U$ and the quotient variety $\overline V=\overline U/G_k$ is well defined. 

\subsection{An important example}\label{example}
Let $\mathcal U=\mathbb P^1_{(z_1:z_2)}\times \mathbb P^1_{(z'_1:z'_2)}$, where $(z_1:z_2)$ and $(z'_1:z'_2)$ are homogeneous coordinates in the projective lines $\mathbb P^1_{(z_1:z_2)}$ and $\mathbb P^1_{(z'_1:z'_2)}$. The group $G_k$ acts on $\mathcal U$ if we put 
$$g_{q,1}: (z_1:z_2)\mapsto (e^{\frac{2\pi qi}{k}}z_1:z_2),\qquad g_{q,1}: (z'_1:z'_2)\mapsto (e^{\frac{2\pi i}{k}}z'_1:z'_2).$$
The quadric $\mathcal U$ is covered by four open subsets $U_{j_1,j_2}=\{z_{j_1}z'_{j_2}\neq 0\}$ isomorphic to $\mathbb C^2$, where $(j_1,j_2)\in \{ 1,2\}^2$.
The functions $u_1=\frac{z_1}{z_2}$ and $u_2=\frac{z'_1}{z'_2}$ are coordinates in the neighborhood $U=U_{2,2}= \mathbb C^2$ and the action of $G_k$ on $U_{2,2}$ is the considered above action of $G_k$ on $U$.  The functions $u_1=\frac{z_1}{z_2}$ and $\widetilde u_2=\frac{z'_2}{z'_1}$ are coordinates in the neighborhood $U_{2,1}= \mathbb C^2$ and the action of the generator $g_{q,1}$ of $G_k$ is the following one: 
$g_{q,1}(u_1,\widetilde u_2)=(e^{\frac{2\pi qi}{k}}u_1,e^{\frac{-2\pi i}{k}}\widetilde u_2)$. Note that an another generator, namely, $g_{k-q,1}=g_{q,1}^{k-1}$ of $G_k$ acts on $U_{2,1}$ as follows: $g_{k-q,1}(u_1,\widetilde u_2)=(e^{\frac{2\pi (k-q)i}{k}}u_1,e^{\frac{2\pi i}{k}}
\widetilde u_2)$. The description of the action of the group $G_k$ on the neighborhoods $U_{1,2}$ and $U_{1,1}$ is left to the reader. 

The quotient space $\mathcal V=\mathcal U/G_k$ is covered by four neighborhoods $V_{j_1,j_2}=U_{j_1,j_2}/G_k$ and it has four singular points 
$\widetilde o_{j_1,j_2}=h_q(o_{j_1,j_2})$, where $h_q:\mathcal U\to\mathcal V$ is the 
quotient mapping and $o_{j_1,j_2}=\{z_{|j_1-1|}=z'_{|j_2-1|}=0\}\in U_{j_1,j_2}$. 

Let $r:\overline{\mathcal V}\to\mathcal V$ be the minimal resolution of the singular points $\widetilde o_{j_1,j_2}\in \mathcal V$ (i.e., $\overline{\mathcal V}$ is a nonsingular surface and $r^{-1}(\widetilde o_{j_1,j_2})$ contains no $(-1)$-curves). 
Denote \begin{itemize} 
\item $C_1=\{z_1=0\}$, $C_{0}=\{ z'_1=0\}$, and $C_{\infty}=\{ z'_2=0\}$ the curves in $\mathcal U$, 
\item $o_0=C_1\cap C_0$ and $o_{\infty}=C_1\cap C_{\infty}$,
\item $\widetilde C_1=h_q(C_1)$, $\widetilde C_0=h_q(C_0)$, and $\widetilde C_{\infty}=h_q(C_{\infty})$ the curves in $\mathcal V$ 
(and, resp., their intersections with neighborhoods $V_{j_1,j_2}$ if it is not lead to misunderstanding),
\item $\widetilde o_0=h_q(o_0)$ and $\widetilde o_{\infty}=h_q(q_{\infty})$,
\item $\overline C_1$, $\overline C_0$, and $\overline C_{\infty} 
$ the proper inverse images of curves $\widetilde C_1$, $\widetilde C_0$, and $\widetilde C_{\infty}$ under the mapping $r:\overline{\mathcal V}\to\mathcal V$ 
(and, resp., their intersections with neighborhoods $\overline V_{j_1,j_2}=r^{-1}(V_{j_1,j_2})$ if it is not lead to misunderstanding),
\item $r^{-1}(\widetilde C_1\cup \widetilde C_0)=\overline C_1\cup E_{1,0}\cup\dots\cup E_{n(0),0}\cup \overline C_0$, where  
$$\overline E_{q,0}:=E_{1,0}\cup\dots\cup E_{n(0),0}=r^{-1}(\widetilde o_{0}),$$ 
\item $r^{-1}(\widetilde C_1\cup \widetilde C_{\infty})=\overline C_1\cup E_{1,\infty}\cup\dots\cup E_{n(\infty),\infty}\cup \overline C_{\infty}$, where  
$$\overline E_{k-q,\infty}:=E_{1,\infty}\cup\dots\cup E_{n(\infty),\infty}=r^{-1}(\widetilde o_{\infty}).$$
\end{itemize}

It follows from \cite{Ba}, Ch. III, Sect. 5, that  $E_{i,*}$ are nonsingular rational curves for $1\leq i\leq n(*)$ (here $*$ is either $0$, or $\infty$), and the graph $\Gamma(r^{-1}(\widetilde C_1\cup \widetilde C_*))=\Gamma_1\oplus \Gamma(\overline E_{n(*),*})\oplus \Gamma_*$ is an oriented chain depicted in Fig. 4 (the vertex $c_1$ is the origin), where $\Gamma_1$ is a graph consisting of a single vertex $c_1$ that corresponds to the curve $\overline C_1$, $\Gamma_*$ is a graph consisting of a single vertex $c_*$ that corresponds to the curve $\overline C_*$ and the graph $\Gamma(\overline E_{n(*),*})$ is the weighted oriented graph of the curve $\overline E_{n(*),*}=E_{1,*}\cup\dots\cup E_{n(*),*}$. The graph $\Gamma(\overline E_{n(*),*})$ is isomorphic to the graph $\Gamma_{\frac{k}{p(*)}}$  (here $p(0)=q$  and $p(\infty)=k-q$) that is defined in Subsection \ref{sect3.2}. Note that 
\begin{equation} \label{ch?} \begin{array}{c} d(\overline E_{n(*),*})=D_{\text{\rm rt},0}(\Gamma(\overline E_{n(*),*}))=\det M(\overline E_{n(*),*})=k, \\ d_{\text{\rm rt},1}(\overline E_{n(*),*})=D_{\text{\rm rt},1}(\Gamma(\overline E_{n(*),*}))=p(*),\end{array} \end{equation}
since the intersection matrix $I(\overline E_{n(*),*})$ of the components of $\overline E_{n(*),*}$ is negative definite. 

\begin{picture}(300,60)
\put(62,30){\circle{4}}
\put(59,18){$c_{1}$}
\put(64,30){\line(1,0){40}}
\put(99,37){$w_{1,*}$}
\put(105,30){\circle*{4}}
\put(100,18){$v_{1,*}$}
\put(145,30){\circle*{4}}
\put(141,37){$w_{2,*}$}
\put(140,19){$v_{2,*}$}\put(108,30){\line(1,0){40}}
\put(145,30){\line(1,0){40}}\put(215,30){\line(1,0){40}} \put(255,30){\circle*{4}}
\put(195,30){$.\, .\, .$}
\put(247,19){$v_{n(*)-1,*}$} \put(245,37){$w_{n(*)-1,*}$} \put(257,30){\line(1,0){40}} \put(299,30){\circle*{4}}
\put(294,37){$w_{n(*),*}$} \put(295,18){$v_{n(*),*}$}
\put(299,30){\line(1,0){40}} \put(341,30){\circle{4}}
\put(337,18){$c_*$}
\put(190,-5){$\text{Fig.}\, 4$}
\end{picture} \vspace{0.5cm}
 
The singularity type of the germ $(\mathcal V, o_{*})$ is denoted in \cite{Ba} by $A_{k,p(*)}$.

Before continuing the proof of Theorem  \ref{main0}, we make the following 
\begin{rem} {\rm Let $p\in C\subset X$ be a non-singular point of a curve $C$ lying in a normal surface $X$ and let $\widetilde h:X \to Y$ be a dominant rational mapping. In this case the restriction $\widetilde h_{\mid C}:C\to Y$ of $\widetilde h$ to $C$ is a regular mapping. This allow us to define the {\it image  of the point $p$ relative to the curve} $C$ as the image $\widetilde h_{\mid C}(p)$ even in the case when $p$ is a undeterminacy point of $\widetilde h$. In particular, in what follows, if we consider a rational dominant mapping $\widetilde h:\mathcal U\to Y$ we will talk about the image $\widetilde h(o_*):=\widetilde h_{\mid C_{1}}(o_*)\in Y$ of $o_*\in \mathcal U$ without specifying that it is the image relative to the curve $C_{1}$.}     
\end{rem}

Consider the following commutative diagram

\begin{picture}(400,70)
\put(152,51){$\widetilde U_*$}
\put(158,47){\vector(3,-2){45}}
\put(168,23){$\widetilde h_q$} 
\put(168,55){\vector(1,0){37}}
\put(180,57){$\widetilde r$}
\put(208,51){$\overline U_*$}
\put(222,55){\vector(1,0){37}}
\put(235,57){$\overline r$}
\put(260,51){$U_*$}
\put(263,47){\vector(0,-1){30}}
\put(265,30){$h_q$}
\put(258,5){$V_*$}
\put(222,8){\vector(1,0){36}}
\put(235,0){$r$}
\put(207,5){$\overline V_*$}
\put(211,47){\vector(0,-1){30}}
\put(198,30){$\overline h_q$} 
\put(190,-15){$\text{Diagram}\, 1.$}
\end{picture} \vspace{0.7cm}

\noindent where $U_*:=U_{j_1,j_2}$ in which the point $o_*$ lies, $V_* =h_q(U_*)$,  $\overline U_*$ is 
the normalization of the fibre product $U_*\times_{V_*} \overline V_*$ over $V_*$ of $U_*$ and $\overline V_*$, the mappings $\overline h_q$ and $\overline r$ are the projections to the factors, and $\widetilde r:\widetilde U_*\to \overline U_*$ is the minimal resolution of singularities of $\overline U_*$. The morphisms in this diagram have the following properties:  
\begin{itemize} 
\item[$(\text{d}_1)$]  $h_q$ and $\overline h_q$ are finite morphisms, $\deg h_q=\deg \overline h_q=k$.
\item[$(\text{d}_2)$] the branch locus of $\overline h_q$ is contained in $\overline E_{n(*),*}$; in particular, $\overline h_q$ is branched in $E_{1,*}$   with multiplicity $k$, since, by Proposition \ref{prop3}, the group $\pi_1(\overline V_*\setminus \overline E_{n(*),*})\simeq \mathbb Z_k$ is generated by the element $x_1$ represented by a simple loop around the curve $E_{1,*}$.  
\item[$(\text{d}_3)$] $r$, $\overline r$, and $\widetilde r$ are birational proper morphisms.
\item[$(\text{d}_4)$] $\text{\rm Sing}\, \overline U_*\subset \overline h_q^{-1}(\cup_{i=1}^{n(*)-1}(E_{i,*}\cap E_{i+1,*}))$.
\item[$(\text{d}_5)$] $\overline r^{-1}(o_*)=\overline h_q^{-1}(E_{1,*}\cap \overline C_1)=h_q^{-1}(E_{1,*})\cap \overline r_{\mid C_1}^{-1}(C_{1})$ is a nonsingular point of $\overline U_*$. 
\end{itemize}

Since $\widetilde U_*$ and $U_*$ are nonsingular surfaces, it follows from $(\text{d}_3)$ that the birational proper morphism 
$\overline r\circ \widetilde r:\widetilde U_*\to U_*$ is a composition $\widetilde{\sigma}_{m(*)}:\widetilde U_*\to U_*$  of $m(*)$ blowups $\sigma_i:U_{i,*}\to U_{i-1,*}$, $i=1,\dots, m(*)$,  of points $o_{i-1,*}\in U_{i-1,*}$, where $o_{0,*}=o\in U_{0,*}=U_*$ and $m(*)$ is the number of irreducible components of the curve $(\overline r\circ \widetilde r)^{-1}(o_*)$. Note that blowups of two different points commute. Therefore we can assume that the centers of the first $\Delta_{k,p(*)}$ blowups $\sigma_i:U_{i,*}\to U_{i-1,*}$ are the points $(\sigma_{1}\circ \dots \circ \sigma_{i-1})_{\mid C_1}^{-1}(o_*)$ for $i=1,\dots,\Delta_{k,p(*)}$ and for $i>\Delta_{k,p(*)}$ the centers of blowups $\sigma_i:U_{i,*}\to U_{i-1,*}$ do not coincide with the preimage of the point $o_*$. 
\begin{lem} \label{delta1} The following equalities hold: $\Delta_{k,p(*)}=p(*)$.
\end{lem}
\proof  We consider the case when $*=0$ only, since the case $*=\infty$ is completely similar.

In the case $*=0$, the functions $u_1,u_2$ defined above, are the coordinates in the neighborhood $U_0\simeq \mathbb C^2$. Then if 
$\sigma: U_{1}\to U_0$ is  the $\sigma$-process with center at $o_0$, then $U_1=U_{1,1}\cup U_{1,2}$,  
where $U_{1,i}\simeq \mathbb C^2$, the functions  $u_{1,1}=\frac{u_1}{u_2},u_{1,2}=u_2$ are coordinates in $U_{1,1}$, and 
$u_{2,1}=u_1,u_{2,2}=\frac{u_2}{u_1}$ are coordinates in $U_{1,2}$. The image $o_1=\sigma_{\mid C_1}^{-1}(o_0)$ of the point $o_0$ lies in $U_{1,1}$, $u_{1,1}=0$ is an equation of $\sigma_{\mid C_1}^{-1}(C_1)$, and $u_{1,2}=0$ is an equation of the exceptional curve of $\sigma$.
The generator $g=g_{q,1}$ of the group $G_k$ acts on $U_1$ and on $U_2$ as follows:
$$ g( (u_{1,1},u_{1,2}))=((e^{\frac{2\pi (q-1)i}{k}}u_{1,1},e^{\frac{2\pi i}{k}}u_{1,2})),\,\, 
g( (u_{2,1},u_{2,2}))=((e^{\frac{2\pi qi}{k}}u_{2,1},e^{\frac{2\pi (k-q+1)i}{k}}u_{2,2})).$$ 
The action of $G_k$ defines the quotient variety $V_1=U_1/G_k$ covered by two neighborhoods $V_{1,1}=U_{1,1}/G_k$ and $V_{1,2}=U_{1,2}/G_k$. The quotient mapping $h_1:U_{1}\to V_1=U_1/G_k$ is finite and it is ramified only at the two points $o_1=\{u_{1,1}=u_{1,2}=0\}\in U_{1,1}$ and $o'_1=\{u_{2,1}=u_{2,2}=0\}\in U_{1,2}$ with multiplicity $k$ and their images $h_{1}(o_1)$ and $h_1(o'_1)$ are singular points of $V_1$ if $q\neq 1$. If $q=1$, then $V_{1,1}$ is a non-singular surface and $h_1$ is ramified in the exceptional curve $\sigma^{-1}(o_0)\subset U_{1,1}$ with multiplicity $k$.      

Note that after $n=\Delta_{k,q}$ blowups of the point $o_0$ we must obtain a surface $U_{n,1}$ that contains the proper inverse image of the curve $C_1$ (and respectively, it contains the preimage $o_n$  of the point $o_0$) such that the induced map $\overline h_q : U_{n,1}\to \overline V_0$ is a regular morphism at  $o_n$ and $\overline h_q(o_0)\in E_{1,0}$.
Therefore it is easy to see (taking into account the considered above properties of action of the group $G_k$ on the blowup $\sigma_1:U_1\to U_0$ of the point $o_0$) that the inequality $\Delta_{k,q}\geq q$ follows from properties $(\text{d}_2)$, $(\text{d}_4)$, and $(\text{d}_5)$.

Similar arguments imply inequality $\Delta_{k,k-q}\geq k-q$. Therefore to prove Lemma, it suffices to prove that $\Delta_{k,q}+\Delta_{k,k-q}=k$. 
To prove this, note that the projection $\text{\rm pr}_1:\mathcal U\to \mathbb P_{(z_1:z_2)}^1$ defines a ruled structure on $\mathcal U$ that is invariant under the action of the group $G_k$ on $\mathcal U$. Therefore this ruled structure defines a ruled structure on $\mathcal V$ and on $\overline{\mathcal V}$. One of the fibres of this ruled structure is a fibre $F$ supported on the curve $\overline F=\overline E_{n(0),0}\cup \overline C_1\cup 
\overline E_{n(\infty),\infty}$ the oriented weighted graph $\Gamma(\overline F)$ of which is depicted in Fig. 5 in which the vertex $v_{n(0),0}$ is the origin and the vertex $v$ corresponding to the curve $\overline C_1$ is the center of the graph $\Gamma(\overline F) =(\Gamma(\overline E_{n(0),0})\setminus\{ c_0\})^{inv}\oplus (\Gamma(\overline E_{n(\infty),\infty})\setminus\{ c_1,c_0\})$ (see also Fig. 4). 

\begin{picture}(300,60)
\put(22,30){\circle*{4}}
\put(9,18){$v_{n(0),0}$}
\put(24,30){\line(1,0){40}}
\put(9,37){$w_{n(0),0}$}
\put(57,37){$w_{n(0)-1,0}$}
\put(69,30){\circle*{4}}
\put(57,18){$v_{n(0)-1,0}$}
\put(64,30){\line(1,0){40}}
\put(110,30){$.\, .\, .$}
\put(128,30){\line(1,0){40}}
\put(168,30){\circle*{4}}
\put(155,37){$w_{1,0}$}
\put(159,19){$v_{1,0}$}
\put(168,30){\line(1,0){40}}
\put(212,30){\line(1,0){40}} \put(210,30){\circle*{4}}
\put(207,19){$v$} \put(205,37){$w$}
\put(254,30){\line(1,0){40}} 
\put(302,30){$.\, .\, .$}
 \put(252,30){\circle*{4}}
\put(245,19){$v_{1,\infty}$} \put(245,37){$w_{1,\infty}$} 
\put(366,30){\circle*{4}}
\put(354,37){$w_{n(\infty),\infty}$} \put(355,18){$v_{n(\infty),\infty}$}
\put(324,30){\line(1,0){40}}
\put(197,0){$\text{Fig.}\, 5$}
\end{picture} \vspace{0.5cm} 

\begin{claim}\label{semidef} The weight $w$ of the vertex $v$ is equal to $w=-(\overline C_1^2)_{\overline{\mathcal V}}=1$.
 \end{claim}
\proof It is well known that the intersection matrix $I(\overline F)$ of the irreducible components of the fibre $F$ is negative semi-definite. Therefore,  by Lemma \ref{prop-1}, the matrix $M(\Gamma(F))$ is positive semi-definite and, in particular, $\det M(\Gamma(\overline F))=0$. Applying decomposition of  the determinant of $M(\Gamma(\overline F))$ by its  $(n(0)+1)$th row (see equation (\ref{eq0})),  it follows from equalities (\ref{ch?}) that
$\det M(\Gamma(\overline F))=w k^2-kq-k(k-q)=0$ and hence $w=1$. \qed \\ 

It follows from properties ($\text{d}_1$) -- ($\text{d}_5$) that  $\widetilde h_q=\overline h_q\circ \widetilde r:\widetilde{\mathcal U}\to \overline{\mathcal V}$ is a finite map  over the curve $\overline C_1$, $\deg \widetilde h_q=k$. Therefore $(\widetilde h_q^{-1}(\overline C_1),
\widetilde h_q^{-1}(\overline C_1))_{\widetilde{\mathcal U}}=\deg \widetilde h_q\cdot (\overline C_1^2)_{\overline{\mathcal V}}=-k$. 

On the other hand, by definition of the numbers $\Delta_{k,q}$ and $\Delta_{k,k-q}$, the sequence of blowups $\overline r\circ\widetilde r:\widetilde{\mathcal U}\to\mathcal U$ decreases the intersection number $(C_1,C_1)_{\mathcal U}=0$ by $\Delta_{k,q}+\Delta_{k,k-q}$. Therefore $\Delta_{k,q}+\Delta_{k,k-q}=k$. \qed \\ 

Consider the case when $k=k_1k_2$, where $k_1$ and $k_2$ are coprime integers. Let $q=m_1k_1+q_1$, where $0\leq m_1<k_2$ and $0<q_1<k_1$. Denote by $G_{k_1}$ a cyclic subgroup of $G_k$ of order $k_1$  generated by $g_{q_1,1}:=g_{q,1}^{k_2}$. It is easy to see that the group $G_{k_1}$ acts on $U$ as follows:
$g_{q_1,1}(u_1,u_2)=(e^{\frac{2\pi q_1i}{k_1}}u_1,e^{\frac{2\pi i}{k_1}}u_2)$. Denote by $\widetilde G_{k_2}=G_k/G_{k_1}$ the quotient group. Then the quotient mapping $h_q:\mathcal U\to \mathcal V_k=\mathcal U/G_k$ can be factorises: 
$$ h_q: \mathcal U\stackrel{h_{q_1}}{\longrightarrow}\mathcal V_{k_1}=\mathcal U/G_{k_1}\stackrel{h_{q,q_1}}{\longrightarrow}\mathcal V_{k}=
\mathcal V_{k_1}/\widetilde G_{k_2}.$$ 

In notation used above, we put $\mathcal V=\mathcal V_k$ and denote by $\overline{\mathcal V}_{k_1}$ the normalization of the fibre product 
$\mathcal V_{k_1}\times_{\mathcal V} \overline{\mathcal V}$ over $\mathcal V$ of $\mathcal V_{k_1}$ and $\overline{\mathcal V}$, denote by $\overline h_{q,q_1}:\overline{\mathcal V}_{k_1}\to \overline{\mathcal V}$ and $\overline r_{k_1}:\overline{\mathcal V}_{k_1}\to \mathcal V_{k_1}$ the projections to the factors. Obviously, there is a finite regular mapping $\overline h_{q_1}:\overline{\mathcal U}\to\overline{\mathcal V}_{k_1}$ such that 
$h_{q_1}\circ \overline r=\overline r_{k_1}\circ\overline h_{q_1}$ and $\overline h_q=\overline h_{q,q_1}\circ \overline h_{q_1}$.

Let $\widetilde r_{k_1}:\widetilde{\mathcal V}_{k_1}\to\overline{\mathcal V}_{k_1}$ be the minimal resolution of the singular points of 
$\overline{\mathcal V}_{k_1}$. Denote by $\overline C_{1,k_1}=\overline h_{q_1}(C_1)\subset \mathcal V_{k_1}$ and by $\widetilde C_{1,k_1}$ the image of $\overline C_{1,k_1}$ under the birational mapping $\widetilde r_{k_1}^{-1}$. It is easy to see that the morphisms $\overline h_{q_1}$, $\overline r_{k_1}$, and $\widetilde r_{k_1}$  have properties similar to properties $(\text{d}_1)$ -- $(\text{d}_5)$. In particular, 
\begin{itemize}
\item[$(\widetilde{\text{d}}_1)$] $\mathcal V_{k_1}=h_{q_1}(U_{2,2})\cup h_{q_1}(U_{2,1})\cup h_{q_1}(U_{1,2})\cup h_{q_1}(U_{1,1})$,  
\item[$(\widetilde{\text{d}}_2)$] $\text{Sing}\, \mathcal V_{k_1}=\{ \widetilde o_{k_1,j_1,j_2}=h_{q_1}(o_{j_1,j_2})\}_{(j_1,j_2)\in \{1,2\}^2}$, 
\item[$(\widetilde{\text{d}}_3)$] $\overline E_{k_1,j_1,j_2}=\overline r_{k_1}^{-1}(\widetilde o_{k_1,j_1,j_2})$\, is\, a chain\, of\, rational\, curves\, for\, each\, pair\, $(j_1,j_2)$, 
    $\text{Sing}\,\overline{\mathcal V}_{k_1}\subset \cup_{(j_1,j_2)}\text{Sing}\,\overline E_{k_1,j_1,j_2}$, and $\text{Sing}\, \overline{\mathcal V}_{k_1}\cap \overline C_{1,k_1}=\emptyset$, 
\item[$(\widetilde{\text{d}}_4)$] $(\cup_{(j_1,j_2)}\overline E_{k_1,j_1,j_2})\cap \overline C_{1,k_1}=
\{ \overline o_{k_1,0}=\overline h_{q_1\mid C_1}(o_0),
\overline o_{k_1,\infty}=\overline h_{q_1\mid C_1}(o_{\infty})\}$, 
\item[$(\widetilde{\text{d}}_5)$] $\overline h_{q_1}:\overline{\mathcal U}\to\overline{\mathcal V}_{k_1}$ is ramified over $\cup_{(j_1,j_2)}\overline E_{k_1,j_1,j_2}$ and $\overline h_{q,q_1}:\overline{\mathcal V}_{k_1}\to\overline{\mathcal V}$ is ramified in $\cup_{(j_1,j_2)}\overline E_{k_1,j_1,j_2}$. 
\end{itemize}  

Let $\widetilde r_{k_1}:\widetilde{\mathcal V}_{k_1}\to\overline{\mathcal V}_{k_1}$ be the minimal resolution of singular points of 
$\overline{\mathcal V}_{k_1}$ that lie in $\overline r_{k_1}^{-1}(\widetilde o_{2,2})$ and $\widetilde C_{1,k_1}$ the image of $\overline C_{1,k_1}$ under the birational mapping $\widetilde r_{k_1}^{-1}$. Note that 
$\overline r_{k_1}\circ\widetilde r_{k_1\mid \widetilde V_{k_1,2,2}}: \widetilde V_{k_1,2,2}\to h_{q_1}(U_{2,2})$ is a (not minimal) resolution of the singular point $\widetilde o_{k_1,0}\in h_{q_1}(U_{2,2})$, where $\widetilde{V}_{k_1,2,2}$ is the inverse image of $h_{q_1}(U_{2,2})$ under the mapping $\overline r_{k_1}\circ\widetilde r_{k_1}$. Denote by $\widetilde{\sigma}_{k_1,2,2}:\widetilde{\mathcal V}_{k_1}\to \overline{\mathcal V}_{k_1,2,2}$ the sequence of blow downs of the exceptional curves of the first kind ($(-1)$-curves) in order to obtain the minimal resolution of the singular point $\widetilde o_{k_1,0}\in h_{q_1}(U_{2,2})$. Let $\Delta_{\frac{k}{q},\frac{k_1}{q_1}}$ be the number of $\sigma$-processe incoming in the sequence  $\widetilde{\sigma}_{k_1,2,2}$ of $\sigma$-processes whose exceptional $(-1)$-curves intersect the images of the curve $\overline C_{1,k_1}$.  

Similarly denote by $\widetilde{\sigma}_{2,2}:\widetilde{\mathcal U}\to \overline{\mathcal U}_{2,2}$ the blow down of the curve 
$(\overline r\circ\widetilde r)^{-1}(o_0)$ to a non-singular point $o'_0$ of $\overline{\mathcal U}_{2,2}$. 

In what follows, by $(C_1,C_1)_Y$ we denote the self-intersection number on a surface $Y$ of the image of curve $C_1\subset \mathcal U$ under a rational mapping $\mathcal U\to Y$.
  
\begin{lem}\label{delta2} If $k=k_1k_2$, where $k_1$ and $k_2$ are coprime integers, and $q=m_1k_1+q_1$, where $0\leq m_1<k_2$ and $0<q_1<k_1$, then 
$\Delta_{\frac{k}{q},\frac{k_1}{q_1}}=m_1$.
\end{lem}
\proof It follows from properties $(\text{d}_1)$ -- $(\text{d}_5)$ and  $(\widetilde{\text{d}}_1)$ -- $(\widetilde{\text{d}}_5)$ that the restrictions of the morphisms $\overline h_{q_1}$ and $\overline h_q$ to a small neighborhood of the curve $C_1$ and the restriction of the morphism $\overline h_{q,q_1}$  to a small neighborhood of the curve $\overline C_{1,k_1}$ are finite mappings of smooth surfaces. Therefore the self-intersections numbers 
$(C_1,C_1)_{\overline{\mathcal U}}$, $(C_{1},C_{1})_{\overline{\mathcal V}_{k_1}}$, and $(C_{1},C_{1})_{\overline{\mathcal V}}$ are well defined and 
it follows from properties $(\text{d}_1)$ -- $(\text{d}_5)$ and  $(\widetilde{\text{d}}_1)$ -- $(\widetilde{\text{d}}_5)$ that  
$(C_{1},C_{1})_{\widetilde{\mathcal U}}=(C_{1},C_{1})_{\overline{\mathcal U}}$ and $(C_{1},C_{1})_{\widetilde{\mathcal V}_{k_1}}=(C_{1},C_{1})_{\overline{\mathcal V}_{k_1}}$. 

By Statement \ref{semidef}, $(C_{1},C_{1})_{\overline{\mathcal V}}=-1$. Therefore 
$(C_1,C_1)_{\overline{\mathcal U}}=-k_1k_2$, $(C_{1},C_{1})_{\overline{\mathcal V}_{k_1}}=-k_2$, since $\deg h_q=k=k_1k_2$ and $\deg h_{q,q_1}=k_2$.

By definition of $\Delta_{\frac{k}{q},\frac{k_1}{q_1}}$, we have
\begin{equation} \label{VkV} (C_{1},C_{1})_{\overline{\mathcal V}_{k_1,2,2}}=-k_2+\Delta_{\frac{k}{q},\frac{k_1}{q_1}}.
\end{equation}
Applying Lemma \ref{delta1} to the rational mapping $\overline{\mathcal U}_{2,2}\to \overline{\mathcal V}$, we obtain equality 
\begin{equation}\label{U} (C_1,C_1)_{\overline{\mathcal U}_{2,2}}=\deg \overline h_q\cdot (C_1C_1)_{\overline{\mathcal V}}+\Delta_{k,q}=
-k_1k_2+q=-k_1k_2+k_1m_1+q_1.
\end{equation}
Applying Lemma \ref{delta1} to the rational mapping $\overline{\mathcal U}_{2,2}\to \overline{\mathcal V}_{k_1,2,2}$, we obtain equality 
$(C_1,C_1)_{\overline{\mathcal U}_{2,2}}=\deg \overline h_{q_1}\cdot (C_1C_1)_{\overline{\mathcal V}_{k_1,2,2}}+\Delta_{\frac{k}{q},\frac{k_1}{q_1}}$ and by (\ref{VkV}), we obtain equality
\begin{equation}\label{U1} (C_1,C_1)_{\overline{\mathcal U}_{2,2}}=k_1(-k_2+\Delta_{\frac{k}{q},\frac{k_1}{q_1}})+\Delta_{k_1,q_1}=
-k_1k_2+k_1\Delta_{\frac{k}{q},\frac{k_1}{q_1}}+q_1.
\end{equation}
Now, Lemma follows from equalities (\ref{U}) and (\ref{U1}). \qed

\section{The end of the proof of Theorem \ref{main0}} \label{sect5} 

\subsection{Monodromy homomorphisms}\label{monod}
 Let $\{ f_P:U_P\to V_p\}\in \mathcal G_{(2),\mathcal P} $ be a germ of generic morphism branched in a curve germ $B_P\subset V_p$ having singular point of $\mathcal P$-simplest type and $f_{P*}:\pi_1(V_p\setminus B_P)\to \mathbb S_{\deg f_P}$ be the monodromy homomorphism of the germ $f_P$. 
\begin{prop} \label{prop10}  The monodromy group $G_{f_P}=\mathbb S_{\deg f_P}$ and $\deg f_P\leq \mu_p( B_P)+1$.
\end{prop}
\proof  By Zariski -- van Kampen Theorem (see, for example, \cite{K}), the fundamental group $\pi_1^{loc}(B_P,p)=\pi_1(V_p\setminus B_P,q)$  is generated by $\mu_p(B_P)$ geometric generators, i.e.,  elements of $\pi^{loc}_1(V_p\setminus B_P,q)$ represented by  simple loops around the curve germ $B_P$ near  non-singular points of $B_P$. By Riemann - Stein Theorem, the germ $f_P$ defines and uniquely defined by monodromy homomorphism $f_{P*}:\pi^{loc}_1(B_P,p)\to \mathbb S_{\deg f_P}$. It follows from generality conditions ({ii}) and ({iii}) that $f_{P*}(\gamma)\in \mathbb S_{\deg f_P}$ is a transposition for a geometric generator $\gamma\in\pi_1(V_p\setminus B_P,q)$. Therefore the monodromy group $G_{f_P}\subset\mathbb S_{\deg f_P}$ is generated by $\mu_p(B_P)$ transpositions   and it acts transitively on the fibre $f_P^{-1}(q)$, since $U_P$ is a connected germ of the non-singular surface $S$. Hence, $G_{f_P}=\mathbb S_{\deg f_P}$ and $\deg f_P\leq \mu_p(B_P)+1$. \qed \\

It is well known that if degree $\deg f_P=2$ of   a finite holomorphic mapping  $f_P$ from a connected germ $U_P$ of non-singular surface to $V_p$ bi-holomorphic to a bi-disc, then the branch curve germ $B_P$ of $f_P$ is  a non-singular curve germ. Therefore there are local coordinates $x,y$ in $V_p$ such that $B_P$ is given in $V_p$ by equation $x^2_1-y_1=0$, i.e., $B_P$ has a "singular"\, point of $\mathcal P_{2,1}$-simplest singularity type. 

Conversely, if $\mu_p(B_f)=1$, that is, the branch curve germ $B_P$ of $f_P$ is  a non-sigular germ, then it follows from Proposition \ref{prop10} that $\deg f_P\leq 2$. 

Further we assume that $\deg f_P\geq 3$.

\subsection{Contraction of the left and right parts of exceptional curve $\overline E_n$}\label{diagr} We use notations introduced above. In particular, 
$\overline{\sigma}_n:V_n\to V=V_p$ is a chain of $\sigma$-processes resolving the singular point $p$ of a curve germ  $B_{P}\subset V$ of $\mathcal P_{k_1,k_2}$-simplest singularity  type, $\overline E_n=\overline{\sigma}_n^{-1}(p)$ is a chain of rational curves, $E_n\subset \overline E_n$ is the exceptional curve of the last blowup $\sigma_n$ incoming in the chain $\overline{\sigma}_n$, $(E_n,E_n)_{V_n}=-1$. 
The intersection matrices of the irreducible components of the left and right parts $E_{n,\text{lt}}$  and $E_{n,\text{rt}}$ of the exceptional curve $\overline E_n\subset V_n$ are negative definite. Therefore, by Grauert's criterion (\cite{Gr}),  there is a holomorphic bimeromorphic contraction $\varphi :V_n\to \overline V$ contracting $E_{n,\text{lt}}$  and $E_{n,\text{rt}}$ to singular points 
$P_{\text{lt}}=\varphi(E_{\text{lt}})$  and $P_{\text{rt}}=\varphi(E_{\text{rt}})$ of $\overline V$. Note that there is a factorisation 
$\overline{\sigma}_n=\psi\circ\varphi$, where $\psi: \overline V\to V_p$ is a map contracting the curve $E_{0}=\varphi(E_{n})\subset\overline V$ to the point $p\in V_p$. 

Denote $\overline B=\varphi(B)\subset \overline V$, $P_B=\varphi(B\cap E_{n})=\overline B\cap E_0\in \overline V$  and let 
$\overline W_{\text{lt}}=\varphi( W_{\text{lt}})\subset \overline V$ (resp., $\overline W_{\text{rt}}=\varphi( W_{\text{rt}})\subset \overline V$ and 
$\overline W_{\overline B}=\varphi( W_{B})\subset \overline V$) be a "tubular"\, neighborhood of $P_{\text{lt}}$ (resp., $P_{\text{rt}}$ and $\overline B$), where $W_{\text{lt}}\subset V_n$ (resp., $W_{\text{rt}}\subset V_n$ and $W_{B}\subset V_n$) is a "tubular"\, neighborhood of the curve $E_{n,\text{lt}}$ (resp., $E_{n,\text{rt}}$ and the curve germ $B$).

The singular points  $P_{\text{lt}}$  and $P_{\text{rt}}$ of $\overline V$ are cyclic quotient singularities. 
The singularity of the germ $(\overline W_{\text{lt}}, P_{\text{lt}})$ at $P_{\text{lt}}$ is of type $A_{k_{\text{lt}},q_{\text{lt}}}$, where $k_{\text{lt}}=d_{\text{lt},0}(\overline E_{n})$ and $q_{\text{lt}}=d_{\text{lt},1}(\overline E_{n})$. 
Similarly, $A_{k_{\text{rt}},q_{\text{rt}}}$ is the singularity type of $(\overline W_{\text{rt}},P_{\text{rt}})$ at $P_{\text{rt}}$, where   
$k_{\text{rt}}=d_{\text{rt},0}(\overline E_{n})$ and $q_{\text{rt}}=d_{\text{rt},1}(\overline E_{n})$, and the singularity types $A_{k_{\text{lt}},q_{\text{lt}}}$ and $A_{k_{\text{rt}},q_{\text{rt}}}$ satisfy the equality 
\begin{equation} \label{fin} k_{\text{lt}}k_{\text{rt}}-k_{\text{lt}}q_{\text{rt}}-k_{\text{rt}}q_{\text{lt}}=1.\end{equation} 

Consider the following commutative diagram 2  

\begin{picture}(300,185)
\put(100,160){$\widetilde U_{\text{r}}$} 
\put(111,157){\vector(2,-3){39}}
\put(202,160){$\widetilde U$} 
\put(205,156){\vector(0,-1){58}}
\put(193,127){$\widetilde{\varphi}$}
\put(135,127){$\overline{\varphi}$}
\put(115,165){\vector(1,0){80}}
\put(215,165){\vector(1,0){80}}
\put(245,168){$\widetilde{\text{f}}$}
\put(155,168){$\widetilde{\text{r}}$}  
\put(100,10){$U_P$} 
\put(115,15){\vector(1,0){180}} 
\put(205,3){$\text{f}_\text{P}$}
\put(300,160){$V_n$} 
\put(298,158){\vector(-2,-3){40}} 
\put(267,127){$\varphi$}
\put(300,10){$V_p$} 
\put(150,85){$\overline U_{\text{r}}$} 
\put(151,82){\vector(-2,-3){40}} 
\put(170,46){$\overline{\psi}$}
\put(125,56){$\widetilde{\psi}$}
\put(163,90){\vector(1,0){35}}
\put(212,90){\vector(1,0){35}}
\put(201,85){$\overline U$} 
\put(204,81){\vector(-3,-2){90}}
\put(228,93){$\overline{\text{f}}$}
\put(177,93){$\overline{\text{r}}$} 
\put(250,85){$\overline V$} 
\put(258,83){\vector(2,-3){40}} 
\put(266,47){$\psi$}
\put(310,85){$\overline{\sigma}_n$}
\put(305,155){\vector(0,-1){133}}
\put(105,155){\vector(0,-1){133}} 
\put(95,85){$\widetilde{\sigma}$}
\put(185,-15){$\text{Diagram}\, 2$}
\end{picture} \vspace{0.7cm}

\noindent in which \begin{itemize}
\item[$1)$]  $f_P:U_P\to V_p$ is  an irreducible germ belonging to $\mathcal G_{(2),\mathcal P}$, $f_P$ is branched in $B_P\subset V_p$
\item[$2)$] $\overline U$ is the normalization of the fibre product $U_P\times_{V_p} \overline V$ over $V_p$ of $U_P$ and $\overline V$, the mappings $\overline{\psi}$ and $\overline f$ are the projections to the factors, and $\overline r:\overline U_r\to \overline U$ is the minimal resolution of singularities of $\overline U$,
\item[$3)$] 
$\widetilde U$  is  the normalization of the fibre product $\overline U\times_{\overline V} V_n$ over $\overline V$ of $\overline U$ and $V_n$ and the mappings $\widetilde{\varphi}$ and $\widetilde f$ are the projections to the factors, 
\item[$4)$] $\widetilde U_r$ is the minimal resolution of singularities of
the normalization of the fibre product $\overline U_r\times_{\overline U} \widetilde U$ over $\overline U$ of $\overline U_r$ and $\widetilde U$ and  the mappings $\overline{\varphi}$ and $\widetilde r$ are the projections to the factors.  
\end{itemize} 

The surfaces $U_P$, $\overline U_r$, and $\widetilde U_r$ are non-singular  and $\overline{\varphi} :\widetilde U\setminus \widetilde{\sigma}^{-1}(P)\to \overline U_r\setminus \widetilde{\psi}^{-1}(P)$ and $\widetilde{\psi} :\overline U_r\setminus \widetilde{\psi}^{-1}(P)\to  U_P\setminus f_P^{-1}(p)$ are bimeromorphic holomorphic maps (here  $f_P^{-1}(p)=P$). Therefore, as is well known, $\overline{\varphi}$ and $\widetilde{\psi}$ are compositions of $\sigma$-processes with centers at points and hence, $\widetilde{\sigma}^{-1}(P)$ and $\widetilde{\psi}^{-1}(P)=\overline r^{-1}(\overline E_0)$ are connected curves  whose irreducible components are all rational curves and where $\overline E_0=\overline f^{-1}(E_0)$. In addition, note that $\overline f^{-1}(P_{\text{lt}})$ (resp., $\overline f^{-1}(P_{\text{rt}})$) is a finite set of points and $\widetilde{\varphi}$ contracts the connected components of $\widetilde f^{-1}(E_{n,\text{lt}})$ (resp., $\widetilde f^{-1}(E_{n,\text{rt}})$) to the points lying in 
$\overline f^{-1}(P_{\text{lt}})$ (resp., $\overline f^{-1}(P_{\text{rt}})$). Consequently, 

$(*)$ {\it the curve $\overline E_0=\overline f^{-1}(E_0)$ is an irreducible rational curve and 

\phantom{aa}\, $\widetilde{\varphi}_{\mid E'_0}: E'_0=\widetilde f_{|E_n}^{-1}(E_{n})\to \overline E_0$ is a birational regular mapping.} 

\subsection{Monodromy homomorphisms (continuation)} \label{subsect-5.3} The mappings 
$$\varphi:V_n\setminus (\overline E_n\cup B)\to \overline V\setminus (E_0\cup \overline B)\,\, \text{and}\,\, \psi:\overline V\setminus (E_0\cup \overline B)\to V_p\setminus B_P$$  
are  biholomorphic. Therefore we can identify the fundamental groups $\pi_1(\overline V\setminus (E_0\cup \overline B))$ and  $\pi_1(V_p\setminus B_P)$ with 
$\pi_1(V_n\setminus (\overline E_n\cup B))$. This identification alow us to identify the  monodromy homomorphisms 
$f_{P*}:\pi_1(V_p\setminus B_P)\to \mathbb S_{\deg f_P}$, $\widetilde f_{*}:\pi_1(V_n\setminus (\overline E_n\cup B))\to \mathbb S_{\deg f_P}$, and $\overline f_{*}:\pi_1(\overline V\setminus (E_0\cup \overline B)\to \mathbb S_{\deg f_P}$.

By Proposition \ref{prop7}, the group $\pi_1(V_n\setminus (\overline E_n\cup B))$ is generated by elements $x_0$, $x_{n_0-1}$, $x_{n_0+1}$,  
where the elements $x_0$, $x_{n_0-1}$, and $x_{n_0+1}$ are  represented by loops that are simple circuits around the curves $B$, $E_{\beta^{-1}(n_0-1)}$, and $E_{\beta^{-1}(n_0+1)}$. Therefore the  monodromy group 
$$G_{f_P}=f_{P,*}(\pi_1(V_p\setminus B_P)=f_{P,*}(\pi_1(V_n\setminus (\overline E_n\cup B)) =\mathbb S_{\deg f_P}$$ 
is generated by permutations $f_{P*}(x_{n_0-1})$, $f_{P*}(x_0)$, and $f_{P*}(x_{n_0+1})$.

By Proposition \ref{prop10}, the monodromy homomorphism $f_{P*}:\pi_1(V_p\setminus B_P)\to \mathbb S_{\deg f_P}$ is an epimorphism. The center of the group $\mathbb S_{\deg f_P}$ is trivial if $\deg f_P\geq 3$. Therefore $x_{n_0}=x_{n_0-1}x_0x_{n_0+1}\in \ker f_{P*}=\ker \overline f_*=\ker \widetilde f_*$, 
since by Proposition \ref{prop10},  the element $x_{n_0}$ which is represented by  a simple loop around the curve $E_{n}$, belongs to the center of  
$\pi_1(V_n\setminus (\overline E_n\cup B))$. Consequently, the element 
$$f_{P*}(x_{n_0-1})f_{P*}(x_0)f_{P*}(x_{n_0+1})=\text{Id}$$
is the identical permutation and the group $\mathbb S_{\deg f_P}$ is generated by both two elements $f_{P*}(x_{n_0-1})$, $f_{P*}(x_0)$ and two elements $f_{P*}(x_0)$, $f_{P*}(x_{n_0+1})$. Therefore

$(**)$ {\it the mapping  $\widetilde f: \widetilde U\to V_n$ is\, branched\, in\, $E_{n,\text{lt}}\cup E_{n,\text{rt}}\cup B$ only\, and\, it\, is\, not 

\phantom{aaa}\, ramified over the curve $E_{n}$, the mapping $\overline f: \overline U\to \overline V$ is ramified over the curve 

\phantom{aaa}\, $B$ and over the points $P_{\text{lt}}$ and $P_{\text{rt}}$, and it  is not ramified in the curve $\overline R$, \newline 
$(***)$  
$\deg \overline f_{\mid \overline R}=\deg \widetilde f_{\mid R_0}=\deg f_P$,
the covering $\overline f_{\mid \overline R}:\overline R\to E_0$ is branched in $P_{\text{lt}}$, $P_{\text{rt}}$, 
\newline 
\phantom{aaaaa}\, and $P_B$ only.} 
 
Note that by Proposition \ref{prop3}, we have 
$$\pi_1(\overline W_{\text{lt}}\setminus P_{\text{lt}})=\pi_1(W_{\text{lt}}\setminus E_{n,\text{lt}})=
\mathbb Z/k_{\text{lt}}\mathbb Z=\langle x_{n_0-1}\,\,|\,\, x_{n_0-1}^{k_{\text{lt}}}=1\rangle,$$  
$$\pi_1(\overline W_{\text{rt}}\setminus P_{\text{rt}})=\pi_1(W_{\text{rt}}\setminus E_{n,\text{rt}})=
\mathbb Z/k_{\text{rt}}\mathbb Z=\langle x_{k_0+1}\,\,|\,\, x_{k_0+1}^{k_{\text{rt}}}=1\rangle,$$
$$\pi_1(\overline W_{\overline B}\setminus\overline B)=\pi_1(W_{B}\setminus B)=\mathbb Z=\langle x_{0}\,\,|\,\, \emptyset\rangle.$$  
Consequently, the  monodromy group 
$$G_{f_P}=f_{P,*}(\pi_1(V_p\setminus B_P)=f_{P,*}(\pi_1(V_n\setminus (\overline E_n\cup B)) =\mathbb S_{\deg f_P}$$ 
is generated by  both two permutations $f_{P*}(x_{n_0-1})$, $f_{P*}(x_0)$ and two permutations  $f_{P*}(x_0)$, $f_{P*}(x_{n_0+1})$, 
the order of $f_{P*}(x_{n_0-1})$  is a divisor of $k_{\text{lt}}$, the order of $f_{P*}(x_{n_0+1})$ is a divisor of $k_{\text{rt}}$, and $f_{P*}(x_{0})$ is a transposition, since $\{ f_P:U_P\to V_p\} \in \mathcal G_{(2),\mathcal P}$.

\begin{claim} \label{symm} Let $g_1=c_1\dots c_t\in \mathbb S_d$ be a decomposition of a permutation $g_1$ into a product of pairwise disjoint cycles 
$c_1,\dots, c_t$ of length $\text{ln}(c_i)=l_i\geq 1$, $\sum_{i=1}^tl_i=d$, and $g_2$ be a transposition. Assume that the elements $g_1$ and $g_2$ generate the group $\mathbb S_d$. Then  either $t=1$ and $l_1=d$, or $t=2$ and $l_1+l_2=d$. In addition, if $t=2$ then $g_1g_2$ is a cycle of length $\text{ln}(g_1g_2)=d$, and if $t=1$ then $g_2g_1=c_1'c'_2$ is a product of two disjoint cycles $c'_1$ and $c'_2$, $\text{ln}(c'_1)+\text{ln}(c'_2)=d$.
\end{claim}
\proof Let $g_2=(1,2)$ and $\langle g_1,g_2\rangle$ be a subgroup generated by $g_1$ and $g_2$. Then there are two possibilities: either the set $\{ 1,2\}$ is a subset of integers contained  in a cycle, say $c_1$, or $1$ is contained in a cycle, say, $c_1$ and $2$ is contained in a cycle, say, $c_2$.
 It is easy to check that in the first case the group $\langle g_1,g_2\rangle$ is a subgroup of the direct product 
$\mathbb S_{l_1}\times\mathbb S_{l_2}\times\dots\times\mathbb S_{l_t}\subset \mathbb S_d$, in the second case the group $\langle g_1,g_2\rangle$ is a subgroup of the direct product $\mathbb S_{l_1+l_2}\times\mathbb S_{l_3}\times\dots\times\mathbb S_{l_t}\subset \mathbb S_d$, and in both cases the restriction to $\langle g_1,g_2\rangle$ of the projection of the direct product to the first factor is an epimorphism. 
Therefore $l_1=d$ in the first case and $l_1+l_2=d$ in the second case, since $\langle g_1,g_2\rangle =\mathbb S_d$. \qed \\    

Denote $t_{\text{lt}}$ (resp., $t_{\text{rt}}$) the number of cycles in the decomposition of the permutation $f_{P*}(x_{n_0-1})$ (resp., $f_{P*}(x_{n_0+ 1})$) into a product of pairwise disjoint cycles. It follows from Statement \ref{symm} that $t_{\text{lt}}+t_{\text{rt}}=3$. Therefore there are two possibilities:
either $t_{\text{lt}}=2$ and $t_{\text{rt}}=1$ (case $(2,1)$), or $t_{\text{lt}}=1$ and $t_{\text{rt}}=2$ (case $(1,2)$). In what follows, we will consider case $(2,1)$ in more detail, because cases $(2,1)$ and $(1,2)$ are similar and in many ways, to consider case $(1,2)$, it is enough to replace $x_{n_0-1}$ with $x_{n_0+1}$ and to replace the indices $\text{lt}$ and $\text{rt}$ in case $(2,1)$ with opposite indices, i.e. replace $\text{lt}$ with $\text{rt}$ and $\text{rt}$ with $\text{lt}$.

In case $(2,1)$,  $f_{P*}(x_{n_0+1})=c_{\text{rt},1}$ is a cycle of length $l_{\text{rt},1}=d$ and $f_{P*}(x_{n_0-1})=c_{\text{lt},1}c_{\text{lt},2}$ is the product of two cycles of lengths $l_{\text{lt},i}=\text{ln}(c_{\text{lt},i})$, $i=1,2$, the sum of which: $l_{\text{lt},1}+l_{\text{lt},2}=d$. Note that $l_{\text{rt},1}$ is a divisor of $k_{\text{rt}}$ (let $k_{\text{rt}}=k_{\text{rt},1}l_{\text{rt},1}$) and $l_{\text{lt},i}$ are  divisors of 
$k_{\text{lt}}$ (let $k_{\text{lt}}=k_{\text{lt},i}l_{\text{lt},i}$).  Let $q_{\text{rt}}=m_{\text{rt},1}k_{\text{rt},1}+q_{\text{rt},1}$, where $0<q_{\text{rt},1}<k_{\text{rt},1}$, and $q_{\text{lt}}=m_{\text{lt},i}k_{\text{lt},i}+q_{\text{lt},i}$, where $0<q_{\text{lt},i}<k_{\text{it},i}$ for $i=1,2$. 

A priori, case $(2,1)$ breaks down into the following six sub-cases: sub-case $(2,1)_{2_1}$ $k_{\text{rt},1}>1$ and $k_{\text{lt},i}>1$ for $i=1,2$,
sub-case $(2,1)_{2_0}$ $k_{\text{rt},1}=1$ and $k_{\text{lt},i}>1$ for $i=1,2$, sub-case $(2,1)_{1_1}$ $k_{\text{rt},1}>1$ and $k_{\text{lt},1}>k_{\text{lt},2}=1$, sub-case $(2,1)_{1_0}$ $k_{\text{rt},1}=1$ and $k_{\text{lt},1}>k_{\text{lt},2}=1$, 
sub-case $(2,1)_{0_1}$ $k_{\text{rt},1}>1$ and $k_{\text{lt},1}=k_{\text{lt},2}=1$, 
sub-case $(2,1)_{0_0}$ $k_{\text{rt},1}=k_{\text{lt},1}=k_{\text{lt},2}=1$.

In case $(2,1)$,  the inverse image $\overline f^{-1}(\overline W_{\text{rt}})=\overline W_{\text{rt},1}\subset \overline U$ and the point 
$\overline P_{\text{rt},1}= \overline f^{-1}(P_{\text{rt}})$  is a singular point of the singularity type $A_{k_{\text{rt},1},q_{\text{rt},1}}$ if $k_{\text{rt},1}>1$. The inverse image $\overline f^{-1}(\overline W_{\text{lt}})=\overline W_{\text{1t},1}\coprod \overline W_{\text{1t},2}\subset 
\overline U$ is a disjoint union of two neighborhoods $\overline W_{\text{1t},1}$ and $\overline W_{\text{1t},2}$ and the point 
$\overline P_{\text{lt},i}= \overline f^{-1}(P_{\text{lt}})\cap \overline W_{\text{1t},i}$, $i=1,2$,  is a singular point of the singularity type $A_{k_{\text{lt},i},q_{\text{lt},i}}$ if $k_{\text{lt},i}>1$. 

Similar to case $(2,1)$, case $(1,2)$ breaks down into six sub-cases also.

Denote by $\widetilde E_0=(\overline r_{|\overline E_0})^{-1}(\overline E_0)$ the proper inverse image of the rational curve 
$\overline E_0\subset \overline U$ and by $\widetilde E_{\text{rt},1}=\overline r^{-1}(\overline P_{\text{rt},1})$ the chain of rational curves if $k_{\text{rt},1}>1$, and by $\widetilde E_{\text{lt},i}=\overline r^{-1}(\overline P_{\text{lt},1})$, $i=1,2$, the chain of rational curves if $k_{\text{lt},i}>1$. Note that the self-intersection numbers of the irreducible components of $\widetilde E_{\text{rt},1}$ (resp, of
$\widetilde E_{\text{lt},i}$) are less than $-1$, since $\overline r:\overline U_r\to\overline U$ is the minimal resolution of the singular points of 
$\overline U$. Therefore $(\widetilde E_0,\widetilde E_0)_{\overline U_r}=-1$, since $\widetilde{\psi}:\overline U_r\to U_P$ is a composition of $\sigma$-processes and $\widetilde{\psi}^{-1}(P)=\widetilde E_0\cup (\cup \widetilde E_{*,*})$, where $\cup \widetilde E_{*,*}$ is the union of the exceptional curves of $\overline r$.

\subsection{Subcase $(2,1)_{2_0}$} In this case $k_{\text{lt},1}$ and $k_{\text{lt},2}$ must be coprime integers, since 
$\overline r^{-1}(\overline E_0)=\widetilde E_{\text{lt},1}\cup \widetilde E_0\cup \widetilde E_{\text{lt},2}$ must be a chain of exceptional curves of the chain of $\sigma$-processes $\psi:\overline U_r\to U_P$, $(\widetilde E_0,\widetilde E_0)_{\overline U_r}=-1$. 
Therefore we must have equality
\begin{equation} \label{fin2} k_{\text{lt},1}k_{\text{lt},2}-k_{\text{lt},1}q_{\text{lt},2}-k_{\text{lt},2}q_{\text{lt},1}=1,\end{equation} 
since $\Gamma(\widetilde E_{\text{lt},i})=\Gamma_{\frac{k_{\text{lt},i}}{q_{\text{lt},i}}}$. In addition, the integers $l_{\text{lt},1}$ and $l_{\text{lt},2}$ must be coprime also. Indeed, if $m=\text{c.g.d.}(l_{\text{lt},1},l_{\text{lt},2})$ 
and $l_{\text{lt},i}=ml'_{\text{lt},i}$,  then $k_{\text{lt}}=ml'_{\text{lt},i}k_{\text{lt},i}$ and $d=l_{\text{rt},1}=k_{\text{rt}}=m(l'_{\text{lt},1}+l'_{\text{lt},2})$. But, $k_{\text{rt}}$ and $k_{\text{lt}}$ are coprime and therefore $m$ must be equal to $1$, and hence, we must have equalities $l_{\text{lt},1}=k_{\text{lt},2}$ and $l_{\text{lt},2}=k_{\text{lt},1}$, i.e., we must have equalities $k_{\text{lt}}=k_{\text{lt},1}k_{\text{lt},2}$ and $k_{\text{rt}}=k_{\text{lt},1}+k_{\text{lt},2}$. Moreover, 
$\widetilde s=(k_{\text{lt},1},k_{\text{lt},2},q_{\text{lt},1},q_{\text{lt},2},q_{\text{rt}},q_{\text{lt}},m_{\text{lt},1},m_{\text{lt},2})\in D_{\mathcal P}$ must be a solution of system of equations (\ref{sys1}), (\ref{syst2}), (\ref{syst3}).  

Conversely, let $\widetilde s=(k_{1},k_{2},q_{1},q_{2},q_{3},q_{4},m_{1},m_{2})\in D_{\mathcal P}$  and $p$ be the singular point of a curve germ $B_P=B_{k_1k_2,k_1+k_2}\subset V_p$ of $\mathcal P_{k_1,k_2,k_1+k_2}$-simplest singularity type, $\overline{\sigma}_n:V_n\to V_p$ be the minimal resolution of the singular point of $B_{P}$, $\overline E_n=\overline{\sigma}_n^{-1}(p)$, $\Gamma(E_{n,\text{lt}})=\Gamma_{\frac{k_1k_2}{q_4}}$, and 
$\Gamma(E_{n,\text{rt}})=\Gamma_{\frac{k_1+k_2}{q_3}}$. 

In notation defined in Subsection \ref{subsect-5.3}, it follows from Proposition \ref{prop7} and Statement \ref{symm}that there is an epimorphism 
$f_{P*}: \pi_1(V_p\setminus B_{P})=\pi_1(V_n\setminus \overline{\sigma}_n^{-1}(B_{P}))\to \mathbb S_{k_1+k_2}$ such that  
$f_{P*}(x_{n_0-1})=c_{\text{lt},1}c_{\text{lt},2}$ is the product of two disjoint cycles of length $\text{ln}(c_{\text{ln},i})=k_i$, $i=1,2$, 
$f_{P*}(x_{n_0+1})=c_{\text{rt},1}$ is a cycle of length $\text{ln}(c_{\text{rn},1})=k_1+k_2$, $f_{P*}(x_{0})$ is a transposition, and $f_{P*}(x_{n_0})=f_{P*}(x_{n_0-1})f_{P*}(x_0)f_{P*}(x_{n_0+1})=\text{Id}$. By  Riemann-Stein Theorem, $f_{P*}$ defines a finite covering $f_P:U_P\to V_p$ of degree $d=\deg f_{P}=k_1+k_2$, where, a priori, $U_P$ is a germ of normal surface and $P=f_P^{-1}(p)$ is, probably, its singular point. 

As above, we can define a diagram of holomorphic mappings similar to Diagram 2. In particular, in this diagram in notation used in Subsection \ref{diagr}, the singular points  $P_{\text{lt}}$  and $P_{\text{rt}}$ of $\overline V$ are cyclic quotient singularities.  
The singularity of the germ $(\overline W_{\text{lt}}, P_{\text{lt}})$ at $P_{\text{lt}}$ is of type $A_{k_1k_2,q_4}$ and 
$A_{k_{1}+k_2,q_{3}}$ is the singularity type of $(\overline W_{\text{rt}},P_{\text{rt}})$ at $P_{\text{rt}}$. The finite mapping $\overline f:\overline U\to\overline V$ is ramified over the points $P_{\text{lt}}$, $P_{\text{rt}}$, and over $\overline B=\varphi(B)$ only and it is not ramified over the rational curve $E_0=\varphi(E_n)$. The inverse image $\overline f^{-1}(\overline W_{\text{rt}})=\overline W_{\text{rt},1}\subset \overline U$ is connected, the point 
$\overline P_{\text{rt},1}= \overline f^{-1}(P_{\text{rt}})$  is a non-singular point, and $\overline f$ is ramified at $\overline P_{\text{rt},1}$ with multiplicity $k_1+k_2$. The inverse image $\overline f^{-1}(\overline W_{\text{lt}})=\overline W_{\text{1t},1}\coprod \overline W_{\text{1t},2}\subset \overline U$ is a disjoint union of two neighborhoods $\overline W_{\text{1t},1}$ and $\overline W_{\text{1t},2}$, the point 
$\overline P_{\text{lt},1}= \overline f^{-1}(P_{\text{lt}})\cap \overline W_{\text{1t},1}$  is a singular point of the singularity type $A_{k_{2},q_{2}}$, and $\overline f$ is ramified at $\overline P_{\text{lt},1}$ with multiplicity $k_1$. Similarly, the point 
$\overline P_{\text{lt},2}= \overline f^{-1}(P_{\text{lt}})\cap \overline W_{\text{1t},2}$  is a singular point of the singularity type $A_{k_{1},q_{1}}$, and $\overline f$ is ramified at $\overline P_{\text{lt},2}$ with multiplicity $k_2$. Note that it easily follows from Hurwitz formula applied to 
$\overline f_{\mid \overline E_0}:\overline E_0\to E_0\subset \overline V$ that $\overline E_0$ is a rational curve. 

The holomorphic bimeromorphic mapping $\overline r: \overline U_r\to \overline U$ is the minimal resolution of the singular points of $\overline U$. 
The curve $\overline r^{-1}(\overline E_0)=\widetilde E_{\text{lt}1}\cup\widetilde E_0\cup \widetilde E_{\text{lt},2}$, where $\widetilde E_0=\overline r_{\mid \overline E_0}^{-1}(\overline E_0)$ and $\widetilde E_{\text{lt},i}=\overline r^{-1}(P_{\text{lt},i})$, $i=1,2$. The weighted graph $$\Gamma(\overline r^{-1}(\overline E_0))=\Gamma(\widetilde E_{\text{lt},1})\oplus \Gamma(\widetilde E_{0})\oplus \Gamma(\widetilde E_{\text{lt},2})$$ is isomorphic to the graph $\Gamma_{\frac{k_2}{q_2}}^{inv}\oplus \Gamma(\widetilde E_{0})\oplus \Gamma_{\frac{k_1}{q_1}}$ and if we apply decomposition of  the determinant of 
$M(\Gamma(\overline r^{-1}(\overline E_0))$ by its  row coorsponding to the vertex $\Gamma(\widetilde E_{0})$ we obtain equality 
$\det M(\Gamma(\overline r^{-1}(\overline E_0))=-(\widetilde E_0,\widetilde E_0)_{\overline U_r}k_1k_2-k_1q_2-k_2q_1$.

To show that $U_P$ is non-singular, let us apply Lemmas \ref{delta1} and \ref{delta2} to mappings $\varphi :\overline V\to V_n$ and $\overline f\circ\overline r:\overline U_r\to \overline V$ to obtain equality 
$$ (\widetilde E_0,\widetilde E_0)_{\overline U_r}=\deg f_P\cdot (E_n,E_n)_{V_n}+\Delta_{\frac{k_1+k_2}{q_3}}+\Delta_{\frac{k_1k_2}{q_4},\frac{k_1}{q_1}}+  
\Delta_{\frac{k_1k_2}{q_4},\frac{k_2}{q_2}}= -(k_1+k_2)+q_3+m_1+m_2.$$
By equality (\ref{Delta}), $m_1+m_2+q_3=k_1+k_2-1$. Therefore $ (\widetilde E_0,\widetilde E_0)_{\overline U_r}=-1$ and hence, 
$\det M(\Gamma(\overline r^{-1}(\overline E_0))=k_1k_2-k_1q_2-k_2q_1$. But, $\widetilde s=(k_{\text{lt},1},k_{\text{lt},2},q_{\text{lt},1},q_{\text{lt},2},q_{\text{rt}},q_{\text{lt}},m_{\text{lt},1},m_{\text{lt},2})$  
is a solution of system of equations (\ref{sys1}), (\ref{syst2}), (\ref{syst3}).  Therefore $s=(k_1,k_2,q_1,q_1)$ is a solution of equation (\ref{sys1}), and hence, $\det M(\Gamma(\overline r^{-1}(\overline E_0))=1$. Consequently, $\overline r^{-1}(\overline E_0)$ is the chain of exceptional curves of a chain $\overline{\sigma}_{\text{ln}(\Gamma(\overline r^{-1}(\overline E_0))}$ of $\sigma$-processes. Finally, we obtain that 
$\widetilde{\psi}=\overline{\sigma}_{\text{ln}(\Gamma(\overline r^{-1}(\overline E_0))}$, i.e., $U_P$ is a germ of non-singular surface.

As a corollary, we obtain Statement (1) of Theorem \ref{main0}.      

\subsection{Subcases $(1,2)_{1_1}$ and $(1,2)_{0_1}$} In these subcases, $k_{\text{rt},1}>k_{\text{rt},2}=1$ and $k_{\text{lt},1}> 1$ 
in subcase $(1,2)_{1_1}$ and $k_{\text{lt},1}= 1$ in subcase $(1,2)_{0_1}$.

We have equalities: $k_{\text{lt}}=k_{\text{lt},1}l_{\text{lt},1}$, $k_{\text{rt}}=k_{\text{rt},1}l_{\text{rt},1}=l_{\text{rt},2}$, and 
$d=\deg f_P=l_{\text{lt},1}=l_{\text{rt},1}+l_{\text{rt},2}=l_{\text{rt},1}+k_{\text{rt}}$. Therefore $l_{\text{rt},1}=1$, because otherwise $l_{\text{rt},1}$ is a common divisor of $k_{\text{lt}}$ and $k_{\text{rt}}$, and hence $l_{\text{lt},1}=1+k_{\text{rt}}$,
\begin{equation} \label{dkq} \begin{array}{c} d=\deg f_P=k_{\text{rt},1}+1,\,\, k_{\text{rt}}=k_{\text{rt},1},\,\,\,  q_{\text{rt}}=q_{\text{rt},1},\\  
k_{\text{lt}}=k_{\text{1t},1}(1+k_{\text{rt},1}),\,\,\, q_{\text{lt}}=k_{\text{lt},1}m_{\text{lt},1}+q_{\text{lt},1},\end{array}\end{equation}
where $0<q_{\text{lt},1}<k_{\text{lt},1}$ in  subcase $(1,2)_{1_1}$ and $m_{\text{lt},1}=0$ in subcase $(1,2)_{0_1}$. 
If we substitute $k_{\text{rt},1}$ instead of $k_{\text{rt}}$, $q_{\text{rt},1}$ instead of $q_{\text{rt}}$, $k_{\text{1t},1}(1+k_{\text{rt}})$ instead of $k_{\text{lt}}$ and $k_{\text{lt},1}m_{\text{lt},1}+q_{\text{lt},1}$ instead of $q_{\text{lt}}$ in equality 
$k_{\text{lt}}k_{\text{rt}}-k_{\text{lt}}q_{\text{rt}}-k_{\text{rt}}q_{\text{lt}}=1$ that must hold, then we obtain  equality 
\begin{equation} \label{k1q1} k_{\text{lt},1}(1+k_{\text{rt}})k_{\text{rt},1}-k_{\text{lt},1}(1+k_{\text{rt}})q_{\text{rt},1}-k_{\text{rt},1}k_{\text{lt},1}m_{\text{rt},1}
-k_{\text{rt},1}q_{\text{lt},1}=1.
\end{equation} 
If $k_{\text{lt},1}>1$, then the integers $k_{\text{lt},1}$ and $k_{\text{rt},1}$ must be coprime also, since 
$\overline r^{-1}(\overline E_0)=\widetilde E_{\text{lt},1}\cup \widetilde E_0\cup \widetilde E_{\text{rt},1}$ must be a chain of exceptional curves of the chain of $\sigma$-processes $\psi:\overline U_r\to U_P$, and equality $(\widetilde E_0,\widetilde E_0)_{\overline U_r}=-1$ must hold. 
Therefore in subcase $(1,2)_{1_1}$, we must have equality
\begin{equation} \label{fin12} k_{\text{lt},1}k_{\text{rt},1}-k_{\text{lt},1}q_{\text{rt},1}-k_{\text{rt},1}q_{\text{lt},1}=1,\end{equation} 
since $\Gamma(\widetilde E_{\text{lt},1})=\Gamma_{\frac{k_{\text{lt},1}}{q_{\text{lt},1}}}$ and 
$\Gamma(\widetilde E_{\text{rt},1})=\Gamma_{\frac{k_{\text{rt},1}}{q_{\text{rt},1}}}$. If we substitute $k_{\text{lt},1}k_{\text{rt},1}-k_{\text{lt},1}q_{\text{rt},1}-1$ in equality (\ref{k1q1}) instead of $k_{\text{rt},1}q_{\text{lt},1}$, 
subtract 1 from the left and right sides of the resulting equality, and then divide the left side of the resulting equality by $k_{\text{lt},1}>1$, as a result, we obtain equality 
\begin{equation} \label{fin12m} m_{\text{lt},1}= k_{\text{rt},1}-q_{\text{rt},1}.\end{equation} 
Summarizing the above, we get that subcase $(1,2)_{1_1}$ occurs only if  

{\it there is $\overline s=(k_{\text{lt},1},k_{\text{rt},1},q_{\text{lt},1},q_{\text{rt},1})\in \mathcal D_{\mathcal P}$, $q_{\text{lt},1}q_{\text{rt},1}\neq 0$, the\, coordinates\, of which are

connected by equality $(\ref{fin12})$, and $d$, $k_{\text{lt}}$, $k_{\text{rt}}$, $q_{\text{lt}}$, $q_{\text{rt}}$ are defined by equalities} $(\ref{dkq})$, $(\ref{fin12m})$.

Conversely in subcase $(1,2)_{1_1}$, if $\overline s\in \mathcal D_{\mathcal P,0}$, where $\overline s=(k_{\text{lt},1},k_{\text{rt},1},q_{\text{lt},1},q_{\text{rt},1})$ are such that $q_{\text{lt},1}q_{\text{rt},1}\neq 0$, then the coordinates of $\overline s$ are connected by equality (\ref{fin12}). Define $k_{\text{lt}}$, $q_{\text{lt}}$, $k_{\text{rt}}$, and $q_{\text{lt}}$ using equalities (\ref{dkq}) and (\ref{fin12m}). The proof that  $\widetilde s=(k_{\text{lt}},k_{\text{rt}},q_{\text{lt}},q_{\text{rt}})$ is a solution of equation  (\ref{sys1}), is the same as the proof of equality (\ref{fin12m}).  Therefore $k_{\text{lt}}$ and $k_{\text{rt}}$ are coprime integers and hence the singularity type of the singular point $p$ of a curve $B_{k_{\text{lt}},k_{\text{rt}}}\subset V_p$ given by equation $x^{k_{\text{lt}}}-y^{k_{\text{rt}}}=0$ is 
$\mathcal P_{\{ k_{\text{lt}},k_{\text{rt}}\}}$-simplest, where $k_{\text{lt}}=k_{\text{lt},1}(k_{\text{rt}}+1)$.  

Let $g_{1}$, $g_{2}$, and $g_3$ be permutations generating the group $\mathbb S_{k_{\text{rt}}+1}$ such that $g_1g_2g_3=\text{Id}$ and $g_1$ is a cycle of length $l_{\text{lt},1}=k_{\text{rt}}+1$, $g_2$ is a transposition, and $g_3=c_{\text{rt},1}c_{\text{rt},2}$ is a product of two disjoint cycles of lengths, resp, $l_{\text{rt},1}=1$ and $l_{\text{rt},2}=k_{\text{rt}}$. In notation defined in Subsection \ref{subsect-5.3}, it easily follows from Proposition \ref{prop7} that there is an epimorphism $f_{P*}: \pi_1(V_p\setminus B_{k_{\text{lt}},k_{\text{rt}}})=\pi_1(V_n\setminus \overline{\sigma}_n^{-1}(B_{k_{\text{lt}},k_{\text{rt}}}))\to \mathbb S_{k_{\text{rt}}+1}$ such that  $f_{P*}(x_{n_0-1})=g_1$, $f_{P*}(x_{0})=g_2$, and 
$f_{P*}(x_{n_0+1})=g_{3}$.  By  Riemann-Stein Theorem, the  epimorphism $f_{P*}$ defines a finite covering $f_P:U_P\to V_p$ of degree 
$d=\deg f_{P}=k_{\text{rt}}+1$, where, a priori, $U_P$ is a germ of normal surface such that $P=f_P^{-1}(p)$ is probably a singular point of $U_P$. 

As above, we can define a diagram of mappings similar to Diagram 2. In particular, in this diagram in notation used above, the singular points  $P_{\text{lt}}$  and $P_{\text{rt}}$ of $\overline V$ are cyclic quotient singularities.  
The singularity of the germ $(\overline W_{\text{lt}}, P_{\text{lt}})$ at $P_{\text{lt}}$ is of type $A_{k_{\text{lt}},q_{\text{lt}}}$ and 
$A_{k_{\text{rt}},q_{\text{rt}}}$ is the singularity type of $(\overline W_{\text{rt}},P_{\text{rt}})$ at $P_{\text{rt}}$. The finite mapping $\overline f:\overline U\to\overline V$ is ramified over the points $P_{\text{lt}}$, $P_{\text{rt}}$, and over $\overline B=\varphi(B)$ only and it is not ramified over the rational curve $E_0=\varphi(E_n)$. The inverse image $\overline f^{-1}(\overline W_{\text{lt}})=\overline W_{\text{lt},1}\subset \overline U$ is connected and $\overline f$ is ramified at $\overline P_{\text{lt},1}$ with multiplicity $k_{\text{rt}}+1$. The point 
$\overline P_{\text{lt},1}= \overline f^{-1}(P_{\text{lt}})$  is a singular point of the singularity type $A_{k_{\text{lt},1},q_{\text{lt},1}}$ in subcase $(1,2)_{1_1}$ and $\overline W_{\text{lt},1}$  is a  non-singular surface in subcase $(1,2)_{0_1}$. The inverse image $\overline f^{-1}(\overline W_{\text{rt}})=\overline W_{\text{rt},1}\coprod \overline W_{\text{rt},2}\subset \overline U$ is a disjoint union of two neighborhoods $\overline W_{\text{rt},1}$ and $\overline W_{\text{rt},2}$, the point $\overline P_{\text{rt},1}= \overline f^{-1}(P_{\text{rt}})\cap \overline W_{\text{rt},1}$  is a singular point of the singularity type $A_{k_{\text{rt},1},q_{\text{rt},1}}$, and $\overline f$ is not ramified at $\overline P_{\text{rt},1}$. The point 
$\overline P_{\text{rt},2}= \overline f^{-1}(P_{\text{rt}})\cap \overline W_{\text{rt},2}$  is a non-singular point and $\overline f$ is ramified at $\overline P_{\text{rt},2}$ with multiplicity $k_{\text{rt}}$. Note that it easily follows from Hurwitz formula applied to 
$\overline f_{\mid \overline E_0}:\overline E_0\to E_0\subset \overline V$ that $\overline E_0$ is a rational curve. 

In subcase $(1,2)_{1_1}$, the holomorphic bimeromorphic mapping $\overline r: \overline U_r\to \overline U$ is the minimal resolution of the singular points of $\overline U$. The curve $\overline r^{-1}(\overline E_0)=\widetilde E_{\text{lt},1}\cup\widetilde E_0\cup \widetilde E_{\text{rt},1}$, where 
$\widetilde E_0=\overline r_{\mid \overline E_0}^{-1}(\overline E_0)$ and $\widetilde E_{\text{lt},1}=\overline r^{-1}(P_{\text{lt},1})$ and 
$\widetilde E_{\text{rt},1}=\overline r^{-1}(P_{\text{rt},1})$.  

The weighted graph 
$\Gamma(\overline r^{-1}(\overline E_0))=\Gamma(\widetilde E_{\text{lt},1})\oplus \Gamma(\widetilde E_{0})\oplus \Gamma(\widetilde E_{\text{rt},1})$ is isomorphic to the graph $\Gamma_{\frac{k_{\text{lt},1}}{q_{\text{lt},1}}}^{inv}\oplus \Gamma(\widetilde E_{0})\oplus \Gamma_{\frac{k_{\text{rt},1}}{q_{\text{rt},1}}}$ and if we apply decomposition of  the determinant of $M(\Gamma(\overline r^{-1}(\overline E_0))$ by its  row corresponding to the vertex $\Gamma(\widetilde E_{0})$ we obtain equality 
$$\det M(\Gamma(\overline r^{-1}(\overline E_0))=-(\widetilde E_0,\widetilde E_0)_{\overline U_r}
k_{\text{lt},1}k_{\text{rt},1}-k_{\text{lt},1}q_{\text{rt},1}-k_{\text{rt},1}q_{\text{rt},1}.$$

To show that $U_P$ is non-singular, let us apply Lemmas \ref{delta1} and \ref{delta2} to mappings $\varphi :\overline V\to V_n$ and $\overline f\circ\overline r:\overline U_r\to \overline V$ to obtain equality 
$$(\widetilde E_0,\widetilde E_0)_{\overline U_r}=\deg f_P\cdot (E_n,E_n)_{V_n}+\Delta_{\frac{k_{\text{lt}}}{q_{\text{lt}}}, \frac{k_{\text{lt},1}}{q_{\text{lt},1}}}+ \Delta_{\frac{k_{\text{rt}}}{q_{\text{rt}}}}= -(k_{\text{rt},1}+1)+m_{\text{lt},1}+q_{\text{rt},1}.$$
By equality (\ref{fin12m}), $m_{\text{lt},1}+q_{\text{rt},1}=k_{\text{rt},1}$. Therefore $ (\widetilde E_0,\widetilde E_0)_{\overline U_r}=-1$ and hence, 
$$\det M(\Gamma(\overline r^{-1}(\overline E_0))=k_{\text{lt},1}1k_{\text{rt},1}-k_{\text{lt},1}1q_{\text{rt},1}-k_{\text{rt},1}1q_{\text{rt},1}=1.$$

Consequently, $\overline r^{-1}(\overline E_0)$ is the chain of exceptional curves of a chain 
$\overline{\sigma}_{\text{ln}(\Gamma(\overline r^{-1}(\overline E_0))}$ of $\sigma$-processes, i.e., 
$\widetilde{\psi}=\overline{\sigma}_{\text{ln}(\Gamma(\overline r^{-1}(\overline E_0))}$, and hence, $U_P$ is a germ of non-singular surface.

In subcase $(1,2)_{0_1}$, the holomorphic bimeromorphic mapping $\overline r: \overline U_r\to \overline U$ is the minimal resolution of the singular points of $\overline U$. The curve $\overline r^{-1}(\overline E_0)=\widetilde E_0\cup \widetilde E_{\text{rt},1}$, where $\widetilde E_0=
\overline r_{\mid \overline E_0}^{-1}(\overline E_0)$ and  $\widetilde E_{\text{rt},1}=\overline r^{-1}(P_{\text{rt},1})$. The weighted graph 
$\Gamma(\overline r^{-1}(\overline E_0))=\Gamma(\widetilde E_{0})\oplus \Gamma(\widetilde E_{\text{rt},1})$ is isomorphic to the graph 
$\Gamma(\widetilde E_{0})\oplus \Gamma_{\frac{k_{\text{rt},1}}{q_{\text{rt},1}}}$. 

Applying Lemmas \ref{delta1} and \ref{delta2} to mappings $\varphi :\overline V\to V_n$ and $\overline f\circ\overline r:\overline U_r\to \overline V$ to obtain equality 
\begin{equation}\label{case3} (\widetilde E_0,\widetilde E_0)_{\overline U_r}=\deg f_P\cdot (E_n,E_n)_{V_n}+\Delta_{\frac{k_{\text{lt}}}{q_{\text{lt}}}}+ \Delta_{\frac{k_{\text{rt}}}{q_{\text{rt}}}}= -(k_{\text{rt}}+1)+q_{\text{lt}}+q_{\text{rt}}.\end{equation}

The surface $U_P$ is non-singular if and only if  $\widetilde{\psi} :\widetilde U_r\to U_P$ is a composition of $\sigma$-processes. Since the self-intersections numbers of the irreducible components of the curve $E_{\text{rt},1}$ are less than $-1$, then we must have equality $(\widetilde E_0,\widetilde E_0)_{\overline U_r}=-1$ and all irreducible components of the curve $E_{\text{rt},1}$ must be $(-2)$-curves. Therefore the singularity type of $\overline P_{\text{rt},1}$ must be $A_{k_{\text{rt}},k_{\text{rt}}-1}$, i.e., $q_{\text{rt}}=k_{\text{rt}}-1$, and it follows from equality (\ref{case3}) that $q_{\text{lt}}=1$.  Consequently, the curve germ $B_P$ must be equisingular equivalent to a curve germ given by equation $x^{k_{\text{rt}}+1}-y^{k_{\text{rt}}}=0$ and $\mu_p(B_P)=k_{\text{rt}}$.

As a corollary, we obtain Statement (2) of Theorem \ref{main0}.  

\subsection{Remaining subcases} To complete the proof of Theorem \ref{main0}, let us show that remaining cubcases are impossible.

Subcase $(2,1)_{2_1}$ and similar subcase $(1,2)_{1_2}$ are impossible, since $(-1)$-curve $\widetilde E_0$ can not intersect three other components of the exceptional curve of a composition of $\sigma$-processes. 

Subcase $(2,1)_{0_0}$ and similar subcase $(1,2)_{0_0}$ are impossible, since in this case $d=2k_{lt}=k_{rt}$ ($d=k_{lt}=2k_{rt}$ in subcase $(1,2)_{0_0}$). On the over hand, it follows from equality (\ref{fin}) that $k_{lt}$ and $k_{rt}$ are coprime integers. 

Subcase $(2,1)_{1_1}$ is impossible, since in this case inequality $d=l_{\text{rt,1}}<k_{\text{rt}}$ contradicts to inequality $k_{\text{lt}}>k_{\text{rt}}$ and equality $d=k_{\text{lt}}+l_{\text{lt},1}$.

Subcase $(2,1)_{0_1}$ is impossible, since in this case inequality $d=l_{\text{rt,1}}<k_{\text{rt}}$ contradicts to inequality $k_{\text{lt}}>k_{\text{rt}}$ and  equality $d=l_{\text{lt},1}+l_{\text{lt},2}=2k_{\text{lt}}$.

Subcase $(2,1)_{1_0}$ is impossible, since in this case equality $d=l_{\text{rt,1}}=k_{\text{rt}}$ contradicts to inequality $k_{\text{lt}}>k_{\text{rt}}$ and  equality $d=l_{\text{lt},1}+l_{\text{lt},2}=l_{\text{lt},1}+k_{\text{lt}}$.

Subcase $(1,2)_{1_0}$ when $k_{\text{lt},1}>1$ and $k_{\text{rt},1}=k_{\text{rt},2}=1$, is impossible, since in this case $d=l_{\text{lt},1}=l_{\text{rt},1}+l_{\text{rt},2}=2k_{\text{rt}}$, but $k_{\text{lt}}=k_{\text{lt},1}l_{\text{lt},1}$ and $k_{\text{rt}}$ must be coprime integers. 

Subcase $(1,2)_{0_2}$ when $k_{\text{lt},1}=1$ and $k_{\text{rt},1}>1$ and $k_{\text{rt},2}>1$ is impossible. Indeed, in this case we have 
$k_{\text{lt}}= l_{\text{lt},1}$, $k_{\text{rt}}= l_{\text{rt},1}k_{\text{rt},1}=k_{\text{rt},2}l_{\text{rt},2}$, and 
$d= l_{\text{lt},1}=l_{\text{rt},1}+l_{\text{rt},2}$. The integers $l_{\text{rt},1}$ and $l_{\text{rt},2}$ are coprime, since $k_{\text{lt}}$ and $k_{\text{rt}}$ must be coprime. Therefore $l_{\text{rt},1}=k_{\text{rt},2}$ and  $l_{\text{rt},2}=k_{\text{rt},1}$, and hence 
$k_{\text{rt}}=k_{\text{rt},1}k_{\text{rt},2}$ and $k_{\text{lt}}=k_{\text{rt},1}+k_{\text{rt},2}$. Note also that the inequality $k_{\text{lt}}=k_{\text{rt},1}+k_{\text{rt},2}>k_{\text{rt}}=k_{\text{rt},1}k_{\text{rt},2}$ must hold. Let $k_{\text{rt},1}\geq k_{\text{rt},2}$. Then
$2k_{\text{rt},1}\geq k_{\text{rt},1}k_{\text{rt},2}$. But, this is possible only if $k_{\text{rt},2}=2$ and $k_{\text{rt},1}+2\geq 2k_{\text{rt},1}$.
Therefore $k_{\text{rt},1}$ must be equal to $2$ also and hence, we have equality $k_{\text{lt}}=k_{\text{rt}}=4$ that is forbidden, since inequality $k_{\text{lt}}>k_{\text{rt}}$ must hold.  \vspace{0.4cm}

\subsection{Proof of Theorem \ref{main00}} The curve germs $B_P$ having $\mathcal P_{k+1,k}$-simplest singularity type  
are dual to smooth curve germs $(C,p)\subset (\mathbb C^2,p)$ for which $p$ is a $(k-1)$-fold inflection point, i.e., $(L,C)_p=k+1$, where $(L,C)_p$ is the intersection number of $C$ and its tangent line $L$ at the point $p$. In the proof of Theorem 5 in \cite{K-Q}, it was shown that  the irreducible germs of dualizing covers associated with curves immersed in the projective plane have extra property. Therefore Theorem \ref{main00} follows from Theorem \ref{main0} and Theorem 1, since $$\deg f_P= \mu_p(B_P)+1=k+1$$ for  germs $\{ f_P:U_P\to V_p\}\in \mathcal G_{(2),\mathbb N}$ branched in the curve germs $B_P$ having $\mathcal P_{k+1,k}$-simplest singularity type. 

\ifx\undefined\bysame
\newcommand{\bysame}{\leavevmode\hbox to3em{\hrulefill}\,}
\fi

\end{document}